\documentclass[onefignum,onetabnum]{siamart190516}

\usepackage{lipsum}
\usepackage{amsfonts,amsmath,amssymb}
\usepackage{graphicx}
\usepackage{epstopdf}
\usepackage{algorithmic}
\usepackage{booktabs}
\usepackage{multirow}
\usepackage{subfig}
\usepackage{tikz}
\usetikzlibrary{plotmarks,arrows.meta,positioning,matrix}
\usepackage{pgfplots}
\pgfplotsset{compat=1.16}
\ifpdf
  \DeclareGraphicsExtensions{.eps,.pdf,.png,.jpg}
\else
  \DeclareGraphicsExtensions{.eps}
\fi
\newcommand{\R}{{\mathbb R}}
\newcommand{\V}{{\mathcal V}}
\newcommand{\E}{{\mathcal E}}
\newcommand{\M}{{\mathcal M}}
\newcommand{\G}{{\mathcal G}}
\newcommand{\I}{{\mathcal I}}
\newcommand{\Range}{{\mathcal Range}}

\DeclareMathOperator\blockdiag{blockdiag}
\usepackage[inline]{enumitem}
\newsiamremark{remark}{Remark}
\newsiamremark{hypothesis}{Hypothesis}
\crefname{hypothesis}{Hypothesis}{Hypotheses}
\newsiamthm{claim}{Claim}

\headers{AMG preconditioners}{P. D'Ambra, F. Durastante, and S. Filippone}

\title{AMG preconditioners for Linear Solvers towards Extreme Scale\thanks{Submitted to the editors DATE.
\funding{This work is supported by the EU under the Horizon 2020 Project Energy oriented Centre of Excellence: toward exascale for energy  (EoCoE-II), Project ID: 824158; and partially supported by the INDAM-GNCS project ``Optimization and Advanced Linear Algebra for problems driven by PDEs''.}}}

\author{Pasqua D'Ambra\thanks{Institute for Applied Computing “Mauro Picone” (IAC), National Research Council (CNR),
Napoli, Italy. 
  (\email{pasqua.dambra@cnr.it}).}
\and Fabio Durastante\thanks{Department  of  Mathematics,  University  of  Pisa, Pisa, Italy, and IAC-CNR. (\email{fabio.durastante@unipi.it})}
\and Salvatore Filippone\thanks{Dept. of Civil Engineering and Computer Engineering, University of Rome Tor-Vergata, Rome, Italy, and IAC-CNR. E-mail \email{salvatore.filippone@uniroma2.it}}}

\usepackage{amsopn}
\DeclareMathOperator{\diag}{diag}

\ifpdf
\hypersetup{
  pdftitle={AMG preconditioners for Linear Solvers towards Extreme Scale},
  pdfauthor={P. D'Ambra, F. Durastante, and S. Filippone}
}
\fi

\begin{document}

\maketitle

\begin{abstract}
  Linear solvers for large and sparse systems are a key element of scientific applications, and their efficient implementation is necessary to harness the computational power of current computers. Algebraic MultiGrid (AMG) preconditioners are a popular ingredient of such linear solvers; this is the motivation for the present work where we examine some recent developments in a package of AMG preconditioners to improve efficiency, scalability and robustness on extreme scale problems. The main novelty is the design and implementation of a parallel coarsening algorithm based on aggregation of unknowns employing weighted graph matching techniques; this is a completely automated procedure, requiring no information from the user, and applicable to general symmetric positive definite (s.p.d.) matrices. The new coarsening algorithm improves in terms of numerical scalability at low operator complexity over decoupled aggregation algorithms available in previous releases of the package. The preconditioners package is built on the parallel software framework \texttt{PSBLAS}, which has also been updated to progress towards exascale. We present weak scalability results  on one of the most powerful supercomputer in Europe for linear systems with sizes up to $O(10^{10})$ unknowns. 
\end{abstract}

\begin{keywords}
  Algebraic MultiGrid, preconditioners, parallel scalability
\end{keywords}

\begin{AMS}
  65F08, 65F10, 65N55, 65Y05
\end{AMS}

\section{Introduction}
\label{sec:intro}

Solving algebraic linear systems of the form:
\begin{equation}
\label{eq:sys}
A x= b, 
\end{equation}
where $A \in \R^{n \times n}$ is a very large and sparse matrix and $x, b \in \R^n$ are vectors, is the most time consuming computational kernel in many areas of computational science and engineering, including more recent fields such as data and network analysis. Efficient and scalable methods and software libraries for solving such systems on the current generation of pre-exascale parallel computers are a key technology for high-resolution and high-fidelity simulations and analysis~\cite{PTRSA2020}. 

The Energy oriented Center of Excellence (EoCoE-II) is a project funded by the EU to enhance computing applications in the field of renewable energy generation and deployment to run on near future exascale parallel computers. In EoCoE-II, four out of the five flagship codes require efficient solution of systems of type \eqref{eq:sys} with sizes larger than $10^{10}$ to have reliable and predictive simulations of multi-physics phenomena; the size rapidly grows towards exascale as the models include more physics and ever more reliable mathematical formulations. 
Often the linear systems originate from time-dependent and/or non-linear models, so that the system~\eqref{eq:sys} has  to be solved multiple times within an iterative procedure; in this case the re-use of operators, such as preconditioners and/or Jacobian matrices, becomes crucial to reduce the overall time to solution. 

In this paper we will discuss our work on extensions of a  package of Algebraic MultiGrid (AMG) preconditioners, motivated by EoCoE-II applications and built on top of the parallel software framework \texttt{PSBLAS} (Parallel Sparse Basic Linear Algebra Subroutines), which implements parallel basic linear algebra operations tailored for iterative sparse linear solvers~\cite{FB2012,FC2000}. The main objective of  the extensions was to improve efficiency, scalability and robustness of the preconditioners package, by employing novel AMG algorithms, to work better at extreme scale. The extensions' implementation  also required a reworking of some basic kernels for parallel distributed matrix and data types management, for data communication, as well as the implementation of some additional Krylov and stationary iterative methods in the \texttt{PSBLAS} framework.

AMG methods are widely used as preconditioners in the iterative solution of large and sparse linear systems, in particular when the systems originate from the discretization of elliptic Partial Differential Equations (PDEs), on structured or unstructured meshes. Their success stems from their potential to achieve optimal complexity, i.e., a computational complexity that is linear in the size of the  system, and on optimal convergence properties which ensure, in some cases, that the number of linear iterations to obtain a given accuracy remains almost constant for increasing system size~\cite{S2001,V2008,XZ2017}.
The latter optimality property is often defined as \emph{algorithmic or numerical scalability}~\cite{CFHJMMMR2006}, and is a necessary condition to have truly scalable AMG.
Many more technological aspects have to be taken into account for efficient implementations on pre-exascale parallel machines, where many thousands of computational cores are distributed among parallel nodes connected by high-speed networks; often some nodes are equipped with accelerators exploiting fine-grained data parallelism and limited size high-speed memories, such as Graphics Processing Units (GPUs).     

Since the introduction of AMG methods, substantial efforts have been put in the design and implementation of efficient packages for their use as parallel preconditioners. {For a detailed description of the most important software projects aimed at preparing  parallel linear solvers and preconditioners for the exascale challenge we refer the reader to~\cite{PTRSA2020}. For comparison purposes we consider here the Hypre library~\cite{FJY2006}; Hypre is very well known for its BoomerAMG multigrid preconditioner package~\cite{MR1908755},  which includes many {\em de facto} standard AMG preconditioners, and is also interfaced by more comprehensive software frameworks, such as PETSc~\cite{petsc} and Trilinos~\cite{trilinos}.} 
In the early 2000s we proposed a package of parallel preconditioners, written in modern Fortran and exploiting basic concepts of object-oriented programming~\cite{BDDF2007,DDF2010}; previous versions of our package  include AMG preconditioners based on the smoothed aggregation method introduced in~\cite{VMB1996}. In this paper we present the most recent activities for extending the package with new aggregation methods, smoothers and coarsest solvers; all these developments have been driven by the desire to achieve extreme scalability and flexibility with respect to the problem size and to the number and features of parallel processors.

The paper is organized as follows. In section~\ref{sec:psblas} we describe the improvements in  the \texttt{PSBLAS} framework for moving towards exascale. Section~\ref{sec:amg4psblas} introduces AMG preconditioners and presents the main features and computational kernels of the new parallel algorithm for coarsening; we focus on the convergence properties of the resulting AMG preconditioner for the different parallel smoothers available in the new version of the package. In section~\ref{sec:testcase} we present the model problem and the parallel machines used in our study; section~\ref{sec:performance_results} provides  an exhaustive analysis of the performance results. Finally section~\ref{sec:end} addresses some concluding remarks.

\section{\texttt{PSBLAS} towards exascale}
\label{sec:psblas}
Our development is  based on the \texttt{PSBLAS} framework~\cite{FC2000,FB2012,CaFiRo:14,Filippone:2017}. Originally introduced for clusters that at the time were large-scale, it has gone through a number of revisions to keep up with the technology development of the past two decades, and the movement towards exascale is no exception. The software framework contains the computational building blocks for Krylov-type linear solvers on parallel computers, as well as support infrastructure to ease the writing of a parallel application using them. In particular, we introduced:
\begin{enumerate}
    \item a framework for handling the mapping between the global index space of the problem and the local portions of the data structures~\cite{FC2000};
    \item the handling and optimization of the halo data exchange, also known as nearest-neighbour data exchange, the essential communication kernel;
    \item an object-oriented architecture that enables choosing storage formats for sparse matrices and switching them at runtime to adapt to the application needs~\cite{FB2012,CaFiRo:14};
    \item a plugin for seamless integration of GPUs~\cite{CaFiRo:14,Filippone:2017}.
\end{enumerate}
During the development of the EoCoE project we have improved the handling of large index spaces requiring 8-byte integers, streamlined the process of setting up the data structures for halo data exchange, and also implemented some new computational kernels prompted by the extension of the preconditioners package. 

In applications dealing with a large and sparse linear system, the system matrix is typically associated with a graph, examples being the discretization mesh of a PDE and the graph representing a complex network. 
All such applications handle the global numbering of the graph, which induces the global numbering of unknowns and matrix indices. In normal practice the global graph/matrix is partitioned and split among processes, and each portion local to a process is  handled through  a local numbering scheme. The solution in \texttt{PSBLAS} is to have an \emph{index map} object contained in the \emph{communication descriptor} to keep track of the correspondence between local and global indices. With the target of handling more than $10^{10}$ degrees of freedom (dofs), it is clear that global indices require 8-byte integers, but that does not necessarily mean that any individual portion will require the same; indeed, having the local portions of the matrices run over 4-byte integers enables  memory savings that can be quite significant, especially when we consider accelerators such as  NVIDIA GPUs which do not support virtual memory and for which memory management is a major concern for the developer. In the current development version of \texttt{PSBLAS} we can choose at configuration time the number of bytes for local and global numbering separately, with the default of using 4 bytes for local   and 8 bytes for global indices. 

One of the main design points of \texttt{PSBLAS} was to make it as easy as possible for the application developer to specify the distribution of the index space, with the only constraint  that each global index/dof point is owned by one process; this is done at the time the descriptor for the index space is created. After this step, all processes need to figure out with whom they need to exchange data. In general there will be some mesh points whose value is needed to carry out the local part of the computation but are not locally owned, and are known as the \emph{halo}; for each halo index, we need to know the owner process. This question would be easy to answer if we had available a vector mapping each index to a process; indeed, that is one of the possible ways to partition an index space, but for very large index spaces this would imply an excessive  memory footprint. Instead, we normally  keep an amount of auxiliary memory that is proportional to the number of local and halo indices on the current process, a solution scalable for increasing number of computational cores; this can be done in two main variants, with a set of hash tables, or by imposing the constraint that the global indices owned by a process must be contiguous. Finding the process owning halo indices is equivalent to establishing a process topology mesh. To help with the construction of the data exchange lists:
\begin{enumerate}
    \item we have devised a new iterative algorithm to identify the owner process for a given non-local index; 
    \item we defined a new interface for the user to provide additional information about the process topology, if available; 
    \item we create, when necessary, a copy of an existing index map employing a renumbering into a block-contiguous format, so as to speed up subsequent halo ownership identification. 
\end{enumerate}

The algorithm to identify index owners is based on the concept of \emph{neighbouring processes}, i.e. processes that own indices needed by each other. The iterative algorithm alternates between probing, in which a global communication among processes establishes the ownership of a subset of the indices, and neighbour communication, in which the newly acquired neighbouring information is used to sweep through the remaining indices; the process is iterated until all index owners have been identified. Each of the communication steps has a limit on the amount of data to be exchanged, and potentially can be split into multiple sub-iterations; the data/memory limit can be adjusted by the user. 
In a block-contiguous numbering it is possible to determine index ownership by knowing the size of each contiguous block; this simplifies the operation, and therefore the library uses this numbering for  auxiliary index maps internal to the preconditioner objects.
All these aspects are handled internally by the software with minimal input by the user;  they influence the setup time of the linear system and of the preconditioner, but have essentially no impact on the runtime of the solver methods. 

\section{AMG for \texttt{PSBLAS}}
\label{sec:amg4psblas}

In~\cite{BDDF2007,DDF2010,ADDDF2010} we proposed a package of AMG preconditioners built on 
top of the \texttt{PSBLAS} framework; the  first version of the package implemented a multilevel version of some domain decomposition preconditioners of additive-Schwarz type and was based on a parallel decoupled version of the smoothed aggregation method described in \cite{VMB1996,TT2000} to generate the multilevel hierarchy of coarser matrices. In this paper we will present a new version of the package, which inherits all the new features of the \texttt{PSBLAS} infrastructure and significantly extends the previous version in terms of algorithms and software modules. We are therefore improving  flexibility, robusteness and computational complexity, but we preserve  the numerical scalability and concurrency of the preconditioners when tens of thousands cores are used, whilst at the same time including support for  GPU accelerators. 
Given the significant number of  changes and the increase in scope, we also decided to change the name of the package to \texttt{AMG4PSBLAS}. In the following we review the main features of the AMG methods, so that we can better describe the new parallel algorithms proposed in this paper and included in \texttt{AMG4PSBLAS}.

\subsection{Introduction to AMG}
\label{sec:introamg}

In the sequel we will consider systems of the type~\eqref{eq:sys}  where $A$ is symmetric and  positive-definite (s.p.d.). 
An AMG method applied to such systems can be viewed as a particular
instance of a general stationary iterative method:
\[
x^{(k)}=x^{(k-1)}+B(b-Ax^{(k-1)}), \ \ k=1, 2, \ldots \ \ \text{given} \ x^{(0)} \in \R^n,
\]
where the matrix $B$ is defined by recursion, as described in the following.
Let $A_l$ be a sequence of coarse matrices computed by the usual triple-matrix Galerkin product:
\[
A_{l+1}=(P_l)^TA_lP_l, \ \ l=0, \ldots, nl-1,
\]
with $A_0=A$ and $P_l$ a sequence of prolongation matrices of size $n_l \times n_{l+1}$, with $n_{l+1} < n_l$ and $n_0=n$. Let $M_l$ be an $A_l$-convergent smoother for $A_l$, i.e., $\|I-M_l^{-1}A_l\|_{A_l} < 1$, where $I$ is the identity matrix of size $n_l$ and $\| \cdot \|_{A_l}$ indicates the $A_l$-energy-norm\footnote{{Given a vector $\mathbf{v}$ and an s.p.d. matrix $A$, the $A$-energy-norm is defined as $\| \mathbf{v} \|_{A} = (\mathbf{v}^T A \mathbf{v})^{\frac{1}{2}}$.}}. The preconditioner matrix $B$ for the well known $V-$cycle, where 1 sweep of pre- and post-smoothing step is applied, is the linear operator corresponding to the multiplicative composition of the following error propagation matrices:
\begin{equation}
\label{eq:amg}
I-B_lA_l= (I-(M_l)^{-T}A_l)(I-P_lB_{l+1}(P_l)^TA_l)(I-M_l^{-1}A_l) \ \ \forall l<nl, 
\end{equation}
assuming that $B_{nl} \approx A_{nl}^{-1}$ is an approximation of the inverse of the coarsest-level matrix \cite{V2008}. The main feature of AMG methods is that the preconditioner setup is completely algebraic, i.e. the prolongation matrices $P_l$ and the smoothers $M_l$ are defined using only information extracted from the corresponding matrices $A_l$. This is the main strength of AMG methods with respect to their geometric counterparts; indeed, in principle they may be applied to general linear systems not  originating from PDEs. The AMG methods' convergence properties depend critically on the ability to define coarsening procedures and prolongation matrices which represent well at each new level the  lower end of the spectrum  of $A$. We refer the reader to~\cite{V2008,XZ2017} for the general theory of AMG methods.  
Two main approaches are used to setup the hierarchy of prolongation matrices in a completely algebraic setting:  {\em classical coarsening} and  {\em coarsening by aggregation}~\cite{S2001,XZ2017}. Classical coarsening separates the original index set into either coarse indices (C-indices), which form the coarse level, and fine indices (F-indices), whose unknowns will be interpolated by the C-indices values. 
The choice of interpolation 
formulas to transfer unknown values from coarser to finer levels defines the corresponding prolongation matrices. Coarsening by aggregation uses disjoint  aggregates of fine unknowns to form the coarse unknowns and, in general, the prolongation matrices are piecewise-constant interpolation matrices (unsmoothed aggregation~\cite{N2010,DV2013,BCHZ2013}) or a smoothed variant of  them (smoothed aggregation~\cite{VMB1996}). In both cases, the way to select aggregates of fine level  variables normally exploits heuristics measures of affinity (also known as {\em strength of connection}) among the variables; these measures have been constructed for systems arising from scalar elliptic PDEs, and often loose their robustness for more general systems. This paper deals with aggregation-based approaches; in particular, we present a parallel aggregation scheme exploiting maximum weight matching in the weighted adjacency graph of the sparse matrices at each level. 

\subsection{Parallel aggregation based on weighted graph matching}
\label{sec:parmatch}

We present the first version of a parallel aggregation scheme for coarsening on large distributed-memory architectures. The method, named {\em coarsening based on compatible weighted matching} was first introduced in~\cite{DV2013} and is already available in the sequential package described in~\cite{DFV2018}. A first parallel version of the method, exploiting fine-grained parallelism and specifically tailored for single GPU device is described in~\cite{BDP2020,BDP2020100041}. 
The method  is independent of any heuristics or a priori information on the {\em near kernel} of $A$, i.e., the lower part of the range of eigenvalues of the system matrix $A$ which is generally used to obtain good-quality aggregates, and it is a completely automatic procedure applicable to general s.p.d. systems. Furthermore, the coupled coarsening based on compatible weighted matching has the  advantage of building  coarse matrices which are well-balanced among parallel processes;  there is  no need for special treatment of process-boundary dofs accounting for inter-processes coupling, as often happens in the  coarsening procedures available in widely used software libraries. Finally, there is a significant flexibility in the choice of the size of aggregates: it is possible to have an  almost arbitrarily aggressive coarsening.

The {\em coarsening based on compatible weighted matching} is a recursive procedure which starts from the adjacency graph $G = (\V,\E)$ associated with the sparse matrix $A$, where the vertex set $\V$ consists of the row/column indices of $A$ and the edge set $\E$ corresponds to the index pairs $(i,j)$ of the nonzero entries in $A$. A \emph{matching} $\M$ in the graph $G$ is a subset of edges such that no two edges are incident on the same vertex.
Graph matching is a general way to obtain partitions of graphs and was proposed for aggregation-based AMG methods in previous works~\cite{BCZ2012,BCKZ2013}. In our method we associate to the graph $G$ a suitable edge weight matrix $C$, computed from the matrix $A$ and an arbitrary vector $w$, and exploit the maximum product matching as a tool to obtain good quality aggregates for fast convergent AMG preconditioners. 

Let $C$ be the following weight matrix:
\begin{equation}
\label{eq:weights}
(C)_{i,j} = c_{i,j} = 1 - \frac{2 a_{i,j} w_i w_j}{a_{i,i} w_i^2 + a_{j,j}w_j^2},
\end{equation}
where $a_{i,j}$ are the entries of $A$ and $w=(w_i)_{i=1}^n$ is a given vector, and let $\M$ be a {\em maximum product matching} in the graph $G$ with edge weight matrix $C$, i.e. $ \M= \arg \max_{\M'} \prod_{ (i,j) \in \M'} c_{ij}$. By applying a maximum product matching we can define the aggregates $\{\G_p\}_{p=1}^{n_p}$ for the row/column indices $\I$ of matrix $A$, consisting of pairs of indices, where $n_p = |\M|$ is the cardinality of the graph matching~$\M$. Equivalently, we are decomposing the index set as
\[
\I = \bigcup_{p=1}^{n_p} \G_p, \quad \G_p \cap \G_r = \emptyset \text{ if } p\neq r.
\] 
We observe that in case of sub-optimal matching, not all vertices will be endpoints of matched edges: we may have \emph{unmatched} vertices. In this case, each unmatched vertex corresponds to a singleton $G_s$, and $n_s$ is the total number of singletons. 
For each edge $e_{i \mapsto j} \in \mathcal{M}$  we can identify the vectors $w_e$ as the orthonormal projection of the vector $w$ on the aggregate $G_p$ defined by the pair $\{ i,j \}$.
Given the $w_e$ vectors, and an ordering of the indices  which  moves all the unknowns corresponding to unmatched vertices at the bottom, we can define a {\em prolongator}:
\begin{equation}\label{eq:prolongator}
P = \left (
\begin{array}{cc}
\tilde{P} & 0\\
0 & W
\end{array}
\right ) \in \R^{n \times n_c}, 
\end{equation}
where:
\[
\tilde{P}=\blockdiag(w_{e_1}, \ldots, w_{e_{np}} ),
\]
$W=\diag(w_s/|w_s|), \; s=1, \ldots, n_s$, corresponds to unmatched vertices and $n_c=n_p+n_s$.
The matrix $P$ we have just built is a piecewise constant interpolation operator whose range includes, by construction, the vector $w$. The recursive application of the above procedure defines an unsmoothed-type aggregation coarsening whose quality and convergence analysis have been discussed in~\cite{DDF2020}; the analysis enables an a posteriori evaluation of the quality of the aggregates and highlights the connection between the choice of the aggregates and the so-called compatible relaxation principle, originally introduced in~\cite{B2000} as a general way to obtain good quality coarsening in AMG.
Good choices for the vector $w$ are obtained by using {\em algebraically smooth vectors}  with respect to $A$, i.e. performing a few iterations of an $A$-convergent smoother $M$ on an arbitrary sample vector (see~\cite{DDF2020} for details). 
In the recursive application of the basic pairwise aggregation method, at each new level the input weighted graph $G=(\V,\E,C)$ corresponds to the adjacency graph of the computed coarse matrix whose weights are computed based on the restriction of the fine-level vector $w$ onto the coarse space.
It is possible to combine multiple steps of the basic pairwise aggregation by computing the product of $m$ consecutive pairwise prolongators, to obtain a more  aggressive coarsening; the resulting  aggregates merge multiple pairs and have an almost arbitrary large size $n_c=2^m$.

As already observed, the prolongation operator $P$ we used for moving between coarse and fine levels corresponds to a piecewise constant interpolation operator. In such a case the V-cycle proves inadequate to obtain an optimal AMG; it is then necessary to employ more robust cycles such as general Algebraic Multilevel Iteration (AMLI)~\cite{V2008}. In this study we employ a Krylov-based MG cycle called the K-cycle, where at each level except the fine and the coarsest ones we apply two iterations of a Flexible Conjugate Gradient (FCG) method with the already defined AMG method starting on the current level as preconditioner~\cite{NV2008}.

An alternative to improve convergence while still employing  a single V-cycle is to  consider the use of a more regular interpolation operator obtained by applying one step of a weighted-Jacobi smoother to the basic piecewise constant interpolation,
as in the smoothed aggregation AMG introduced in~\cite{VMB1996}. The actual prolongator $\overline{P}$ is then obtained from $P$ as $\overline{P}=(I-\omega D^{-1}A)P$, where $D=\diag(A)$ and $\omega = 1/\|D^{-1}A\|_{\infty} \approx 1/\rho(D^{-1}A)$, with $\rho(D^{-1}A)$ the spectral radius of $D^{-1}A$. We will provide results for this choice in the experiments section.

The main kernels in the coarsening based on compatible weighted matching are therefore the computation of a maximum product matching in an edge-weighted graph, the sparse matrix-sparse matrix products needed to apply multiple steps of pairwise aggregation and to compute the  Galerkin product, and, finally, the sparse matrix-vector product with  matrix $P^T$ needed to restrict the vector $w$ from finer to coarser level.

The computation of a maximum product matching in an edge-weighted graph can be transformed into the classic {\em maximum weight matching} by maximizing the following additive weight function:
\[
-\sum_{ (i,j) \in \M} (\log \max_i | c_{ij} | - \log |c_{ij} |) , \ c_{ij} \neq 0.
\]
The computation of an optimal maximum weight matching in a parallel setting is not a trivial task because the corresponding algorithm is intrinsically sequential. Indeed, practical approaches are based on a relaxation of the optimality requirement and look for solutions which approximate the optimal weight.   
In our parallel coarsening we use the \texttt{MatchBox-P} software library, which implements the parallel algorithm for the computation of half-approximate maximum weight matching described in~\cite{CDGHP2011}. 
This algorithm has a complexity $\mathcal{O}(|E| \Delta)$, where $|E|$ is the cardinality of the graph edge set and $\Delta$ is the maximum vertex degree\footnote{The degree of a vertex is the number of edges that are incident on it.}, and guarantees a solution that is at least half of the optimal weight. The \texttt{MatchBox-P} algorithm is based on the idea of identifying {\em locally dominant} edges, i.e., edges with largest weight for both  end-vertices. 
\texttt{MatchBox-P} supports data distributions arising from general row-block sparse matrix parallel distribution and implements  message aggregation and overlapping between communication and computation to reduce the impact of the data communication on the parallel efficiency. It uses MPI asynchronous communication functionalities and was demonstrated to have good weak scalability properties for large graphs resulting from 2D~PDE discretization on up to tens of thousands of parallel cores.

The construction of the multilevel hierarchy requires the computation of a number of parallel matrix-matrix products where the matrices are sparse. The implementation of this kernel can logically be split in two parts: the data exchange and the local computations. If we are computing the product $A \times B$, any given process will own a subset of the rows of $A$; to compute the product it will need not only the same subset of rows of $B$, but also the rows of $B$ corresponding to the column indices of non-zero entries in $A$, including \emph{halo} indices. The data exchange is thus essentially the same as the halo data exchange for vectors in a matrix-vector product, keeping in mind that items being exchanged are rows of sparse matrices, and in general the rows will be of differing lengths, and this must be accounted for in the  use of the basic \texttt{PSBLAS} data structures for halo data exchange. An additional twist applies here: during the process of collecting remote matrix rows, the set of column indices on each process will in general expand, and hence there will be an impact on the support data structures. Once the matrix pieces have been converted to local numbering we can use any serial sparse-matrix by sparse-matrix product code to complete the computation. 

\subsection{Parallel Smoothers and Coarsest Solvers}
\label{sec:smoothers}

Scalable AMG relies on smoothers that are both  highly parallel  and robust. A smoother is a  convergent iterative method which is able to swiftly reduce the error components associated with the large eigenvalues of the system matrix. For s.p.d. matrices, we consider A-convergent smoothers, represented by a matrix operator $M$, such that $\|I-M^{-1}A\|_A < 1$.
A common choice with a good smoothing factor is the Gauss-Seidel (GS) method, which is based on the convergent splitting of the system matrix $A=M-N$, with $M=L+D$ and $N=-L^T$, where $D=\operatorname{diag}(A)$ and $L$ is the lower triangular part of $A$. However, the GS method is intrinsically sequential; in parallel we typically use an inexact block-Jacobi version of it, where in the portion of the row-block local to each process the method acts as the GS method. This version is commonly known as {\em Hybrid GS (HGS)} and has been demonstrated to be a convergent smoother, with better smoothing properties than the block-Jacobi method when the local diagonal block of the matrix is sufficiently large with respect to the off-diagonal portion~\cite{BFKY2011}.  
This dependence of the smoothing properties of HGS on the size of local diagonal blocks of the distributed matrix may have an important impact for AMG when very large numbers of parallel cores have to be employed. Indeed, for coarser-level matrices, where the size of local diagonal block is reduced and the matrix becomes denser, the off-diagonal part becomes heavier and the smoothing properties of HGS may degrade.
In this setting a weighted version of the method, named $\ell_1-$HGS, can be useful. $\ell_1-$HGS is $A$-convergent and its smoothing properties appear less sensitive to the number of parallel cores and to the size of local diagonal blocks of the matrix, as discussed in~\cite{BFKY2011}. Let $A$ be divided in $np$ blocks of size $nb \times n$, as in a general row-block parallel distribution involving $np$ processes, and let $A_{pp}$ be the corresponding diagonal block of $A$. Starting from the usual decomposition of the block as $A_{pp}=L_{pp}+D_{pp}+L_{pp}^T$, where $D_{pp}=\operatorname{diag}(A_{pp})$ and $L_{pp}$ is the lower triangular part of $A_{pp}$, the $\ell_1-$HGS method is an inexact block-Jacobi method defined by the following diagonal block matrix:
\begin{equation}
\label{eq:ml1hgs}
M_{\ell_1-HGS}=\operatorname{diag}((M_{\ell_1-HGS})_p)_{p=1, \ldots np},
\end{equation}
where $(M_{\ell_1-HGS})_p= L_{pp}+D_{pp}+D_{\ell_1p}$, with $D_{\ell_1p}$ the diagonal matrix whose entries are:
\begin{equation}\label{eq:l1weights}
((d_{\ell_1})_i)_{i=1}^{nb}= \sum_{j \in \Omega_p^{nb}} |a_{ij}|.
\end{equation}
In the above definition $\Omega_p^{nb}$ is the set of indices of the $p$-th block of $A$ outside the diagonal block $A_{pp}$. Note that the matrix blocks do not necessarily have the same size, i.e. $nb$ can different for each of the $np$ blocks.   
For $nb=1$, $\ell_1-$HGS becomes the so-called $\ell_1-$Jacobi smoother which is $A-$convergent and generally has good smoothing properties. It is of particular interest for our aims, since its application allows to exploit the high level of parallelism of GPU accelerators, whereas the application of the $\ell_1$-HGS method is inefficient due to need of apply the intrinsically sequential triangular solve involved in the inversion of the diagonal blocks in \eqref{eq:ml1hgs}. 

In the present study we also experiment with sparse approximate inverses since they do not require solving an auxiliary linear systems and have  as their main kernel a sparse matrix-vector product, whose efficient implementation is available in the \texttt{PSBLAS}  GPU plugin. There exist several different algorithms for computing a
sparse approximate inverse; we focus here on the inversion and sparsification of an incomplete factorization introduced in~\cite{va-1999}. This strategy is based on the application of a sparse inversion technique for the triangular factors of an existing incomplete factorization in the form $M=L D L^T$, where, as usual, $D$ is a diagonal matrix and $L$ is lower triangular with an all ones main diagonal. In this way an expression for 
\[
M^{-1}=L^{-T}D^{-1}L^{-1}=Z D^{-1} Z^T
\]
is obtained, and the application of the smoother is reduced to the computation of a matrix-vector product. To have sparse expressions for the incomplete factorization of $A^{-1}$ it is necessary to employ a sparsification process during the computation of the matrix $M^{-1}$, i.e., a sparsification process for the matrix $Z$. As discussed in~\cite{va-1999}, and analyzed in detail in~\cite{bertaccini2016sparse}, the sparsification can be based on either a thresholding procedure or a positional dropping; we refer to~\cite[Chapter~3.5]{MR3793630} for a complete discussion. As for the other smoothers, we will always consider the case of computing the decomposition $M$ for the diagonal blocks $A_{pp}$ of $A$ in a parallel block-Jacobi setting, by means of the positional sparsification technique, i.e., the INVK variant from~\cite[Algorithm~2]{bertaccini2016sparse}. 

We can also  consider  the $\ell_1$-modification of this smoother by using again the modified block decomposition with the $D_{\ell_1,p}$ matrix in~\eqref{eq:l1weights}, that is, by computing the approximate inverse of the $A_{pp} + D_{\ell_1,p}$ blocks.

It is  possible to define damped or weighted versions of all these  methods, that is, to introduce the modified iteration
\begin{equation*}
x^{(k+1)} = (I - \omega M^{-1}A)x^{(i)} + \omega M^{-1}b;
\end{equation*}
a sufficient condition for convergence is that $\omega \in (0,2/\rho(M^{-1}A))$. When all the eigenvalues of $M^{-1}A$ are real and positive the best convergence ratio is obtained at the midpoint of the extreme eigenvalues $\omega_{\text{opt}} = 2/(\lambda_{\text{min}}(M^{-1}A) + \lambda_{\text{max}}(M^{-1}A))$.

Following the line of reasoning in~\cite{BFKY2011}, to evaluate the smoothing properties of the selected smoothers we can perform a smoothing analysis, for a fixed projector $P$, by looking at the spectral radius of the error propagation matrix in~\eqref{eq:amg} for the two grid method using only a pre-smoothing step: 
\begin{equation*}
\|I-BA\|_A^2 = \| (I - P(P^TAP)^{-1}P^TA)(I-M^{-1}A) \|_A^2 \leq 1 - \frac{1}{K};
\end{equation*}
for the convergence constant $K$ it holds
\begin{equation*}
K = \sup_{e} \frac{\|(I-PR) e \|_{\overline{M}}^2}{\|e\|_A^2} \geq 1, \qquad \overline{M} = M^T(M^T + M - A)^{-1} M,
\end{equation*}
where $R$ is any matrix such that $PR$ is a projection onto $\Range(P)$.
If we use the prolongator $P$ in~\eqref{eq:prolongator}, we have that $R = (P^TDP)^{-1}P^TD$, with $D=\operatorname{diag}(A)$, defines the coarse grid variables, and is such that $R P = I_{n_c}$, with $I_{n_c}$  the identity matrix of size $n_c$. We can then write the associated full-rank matrix $S \in \mathbb{R}^{(n-n_c) \times n}$ such that $RS = 0$ by the D-orthogonal vectors to the diagonal blocks of $P$.
Finally, we can explicitly compute $K$ from~\cite{FV2004}:
\begin{equation}\label{eq:smoothingfactor}
K = \frac{1}{\lambda_{\text{min}}\left[(S^T \overline{M} S)^{-1}(S^T A S)\right] },
\end{equation}
and use it to compare the smoothing properties for different choices of $M$ while keeping fixed the remaining parts of the hierarchy.
In Table~\ref{tab:smoothingfactor} we report the constant $K$ relative to the application of the  parallel smoothers we have discussed in this section in a two-level method built by our matching-based parallel coarsening, for a 3D Laplace equation with homogeneous Dirichlet boundary conditions on a unit cube.
\begin{table}[htbp]
	\centering
	\caption{Convergence constant $K$ from~\eqref{eq:smoothingfactor} for various smoothers for the homogeneous 3D Laplacian problem with parallel coarsening based on compatible weighted matching. The global size of the problem is $m = 2^4 \times 2^4 \times 2^4$ distributed over $np$ processes using a 3D block distribution. By INVK we denote here the approximate inverse computed from the incomplete LU factorization with $0$ levels of fill-in, and admitting a single level of fill-in in the inversion procedure.}
	\begin{tabular}{rrcccccccc}
		\multicolumn{8}{c}{Unsmoothed  prolongator with 3 sweeps of pairwise aggregation} \\
		\midrule
		m    & np & HGS  & $\ell_1-$HGS & HINVK & $\ell_1-$INVK & HINVK & $\ell_1-$INVK\\
		\cmidrule(l{2pt}r{2pt}){5-6} \cmidrule(l{2pt}r{2pt}){7-8}
		&   &      &        &       \multicolumn{2}{c}{$\omega = 1$}         & \multicolumn{2}{c}{$\omega_{\text{opt}}$}\\
		\midrule
		4096 &  1   & 1.3766 & 1.3766 & 1.5562 & 1.5562 & 1.0024 & 1.0024 \\
		2048 &  2   & 1.4194 & 1.5270 & 1.5273 & 1.7196 & 1.0356 & 1.0812 \\
		1024 &  4   & 1.4587 & 1.6621 & 1.6093 & 2.2149 & 1.2593 & 1.3638 \\
		512  &  8   & 1.4744 & 1.7803 & 1.8284 & 2.6713 & 1.5345 & 1.6147 \\
		256  &  16  & 1.4945 & 1.8230 & 1.8608 & 2.7307 & 1.5589 & 1.5816 \\
		128  &  32  & 1.5149 & 1.8682 & 1.8977 & 2.7972 & 1.5875 & 1.5669 \\
		64   &  64  & 1.5335 & 1.9162 & 1.9390 & 2.8715 & 1.6198 & 1.5736 \\
		32   & 128  & 1.5880 & 2.0343 & 2.0272 & 3.0707 & 1.7469 & 1.6461 \\
		16   & 256  & 1.6406 & 2.1594 & 2.1440 & 3.3688 & 1.9108 & 1.7678 \\
		8    & 512  & 1.6665 & 2.3088 & 2.3137 & 3.7280 & 2.1317 & 1.9254 \\
		\midrule
		\multicolumn{2}{c}{$l_1$-Jacobi} & \multicolumn{6}{c}{5.6220} \\
		\bottomrule
	\end{tabular}
	\label{tab:smoothingfactor}
\end{table}

We can observe that all the smoothers considered here have good smoothing properties, leading to good convergent two-level AMG. HGS appears to have the best behavior, with no significant deterioration for decreasing block size.
HINVK also produces good convergence constants for decreasing block size. Using the optimal value for $\omega$, as expected, improves the convergence, however using $\omega=1$ does not degrade $K$  in a significant way. Finally, we note that $\ell_1$ versions of the smoothers do not seem to improve $K$:  the constant value is generally worse than that of the corresponding standard smoother. 

In \texttt{AMG4PSBLAS} we provide interfaces to some of the widely used parallel direct solvers, such as \texttt{SuperLU}~\cite{L2005} and \texttt{MUMPS}~\cite{ADL2000}. However,  using direct solvers at the coarsest level of an AMG method on many thousands of parallel cores can be very expensive, because the  coarsest-level matrix tends to have a small size, therefore the cost of data communication dominates the local arithmetic computations. Sometimes, as in the \texttt{Hypre} library, the coarsest-level system is solved by a subset of parallel cores, leaving many processes idle, then the solution is redistributed to all the processes in order to proceed.
In this study we focus on a distributed coarsest solver involving all the parallel cores and limiting the maximum size of the coarsest-level matrix to a sufficiently large size. In order to have an effective parallel solution of the corresponding coarsest system we exploit a preconditioned version of the Conjugate Gradient (CG) Method, coupled with block-Jacobi preconditioners, as available in \texttt{PSBLAS}.

\section{Model Problem and Computational Environments}
\label{sec:testcase}

To  benchmark our code, we consider the solution of the Laplace equation in 3D, i.e., the self-adjoint elliptic partial differential equation of the form
\begin{equation}\label{eq:bench}
- \nabla \cdot ( K \nabla u ) = \mathbf{f}, \quad\text{ in } [0,1]^3,
\end{equation}
with homogeneous Dirichlet boundary conditions, $K= 1$ and unitary right-hand side. 
The discretization of this  problem is obtained by the usual 7-points finite difference stencil.

The choice of the above model problem satisfies two requirements: on the one hand it allows for comparisons with results available in the literature~\cite{NN2015,Baker2012}, and on the other hand it is of interest because it is representative of  the computational kernel at the core of two of the flagship applications included in the EoCoE-II project. Indeed, the Large Eddy Simulation of wind-turbines models in Alya~\cite{VAZQUEZ201615}, and the simulation of 3D hydrologic models, simulating subsurface flows based on the solution of Darcy's equation in Parflow~\cite{KOLLET2006945}, both reduce to elliptic equations similar in nature to~\eqref{eq:bench}.

For the performance results in section~\ref{sec:performance_results} we investigate weak scalability: we  aim  to  show that our solvers scale (almost) linearly with the amount of resources. For this task we consider two computational frameworks. The first one is a purely distributed context in which the resources to be considered for the weak scaling are the number of MPI ranks. The second one is a hybrid context in which we exploit both the distributed MPI architecture and the use of GPU accelerators; therefore the resources in this case are CUDA devices, i.e., the number of accelerators. In both cases, the resulting linear system is distributed by means of a block 3D distribution:  the MPI ranks are  arranged in a Cartesian grid similar to the cubic grid induced by the discretization, and  the  unknowns in each discretized  sub-block are then assigned to the corresponding  process. The weak scalability tests is then performed by  a number of subsequent refinements of the mesh in such a way that, as the number of processes is increased, the number of unknowns per process is constant.

The machine we used for our experiments is the \textbf{Piz Daint} supercomputer from the Partnership for Advanced Computing in Europe (PRACE). Piz Daint is operated by the Swiss National Supercomputing Center. It is based on the Cray Model XC40/Cray XC50 architecture with 5704 hybrid compute nodes (Intel Xeon E5-2690 v3 with Nvidia Tesla P100 accelerator), and 1431 multicore compute nodes (Intel Xeon E5-2695 v4), using the Cray Aries routing and communications ASIC with Dragonfly network topology. The gross total of cores is 387,872 and it is ranked at the twelfth position in the November~2020 \emph{TOP500}\footnote{Refer to~\url{https://www.top500.org/lists/top500/2020/11/}.}.

\section{Performance Results}
\label{sec:performance_results}

To measure the overall performance of the parallel AMG preconditioners introduced in this paper we focus on both \emph{algorithmic} and \emph{implementation} scalability. Perfect algorithmic scalability is achieved when $\rho(B^{-1}A) \approx 1$ independently of the global size $n$ of the linear system;  implementation scalability on the other hand corresponds to having an optimal application cost for $B^{-1}$ that is $\mathcal{O}(n)$ flops per iteration while achieving parallel speedup proportional to the number of processes employed. In all cases we use the AMG methods as preconditioners for a Flexible Conjugate Gradient (FCG) algorithm. Algorithmic scalability can  be analyzed by looking at the number of iterations needed by FCG to achieve a relative  residual norm of $10^{-6}$ as the size $n$ grows; implementation scalability is analyzed by considering both the total solve time for the procedure as well as the average time per iteration. We also analyze the timings for the preconditioner setup with increasing number of processes; all the timings reported in the figures are in seconds.

Regarding the memory footprint of the proposed multigrid hierarchies  {as well as the cost of the application of a V-cycle}, a quantitative measure is given by the operator complexity
\[\text{opc}=\frac{\sum_{l=0}^{\text{nl}-1}\operatorname{nnz}(A_{l})}{\operatorname{nnz}(A_0)} > 1.\]

 Note that for both AMG cycles we employed (K- and V-cycle), we always apply 1 pre-smoothing and 1 post-smoothing step. As coarsest solver we use a version of the CG method coupled with a block-Jacobi preconditioner, with incomplete LU factorization (ILU(0)) on the diagonal blocks; the iterative solution of the coarsest system is stopped when the relative residual is less than $10^{-4}$ or the number of iterations is larger than $30$. In all  experiments discussed in this paper we applied the parallel coarsening based on compatible weighted matching starting from a vector $w$ of all ones, since this vector is in the near kernel of our model test case. Furthermore, we used the library default approach which sets a target size of the coarsest-level matrix and stops the coarsening procedure as soon as the coarse matrix size is less than or equal to the target; in the experiments with the pure MPI version of the library we set $maxsize=200 \times np$, where $np$ is the number of cores.
 
 To identify the different algorithmic variants implemented in \texttt{AMG4PSBLAS} and analyzed in our experiments, we use the labeling convention described in Table~\ref{tab:tested_methods}: the overall name reported in the figures and in the analysis is obtained by combining together the labels of the various components{, e.g., KA1S1CS1 is the preconditioner employing a K-cycle, a hierarchy built with three sweeps of parallel coarsening based on matching and unsmoothed prolongator (unsmoothed parallel matching, for brevity), one sweep of the forward/backward Hybrid Gauss-Seidel smoother and the preconditioned CG method as coarsest solver}.
\begin{table}[h]
    \centering
    \caption{{Strings identifying each preconditioner are built by combination of the strings identifying the various algorithmic variants.}}
    \label{tab:tested_methods}
    \begin{tikzpicture}[node distance=0.25cm]
    \matrix (label) [matrix of nodes,draw,sharp corners]{%
        K & A1 & S2 & CS1 \\
    };
    \matrix (cycle) [matrix of nodes,below=of label-1-1,xshift=-3.8cm]{%
        \underline{Cycle}\\
        K\\
        V\\
    };
    \matrix (aggregation) [matrix of nodes,below=of label-1-2,xshift=-1.7cm,text width=4cm]{%
        \underline{\textbf{A}ggregation}\\
        1) Unsmoothed Parallel Matching: 3 sweeps \\
        2) Unsmoothed Parallel Matching: 4 sweeps \\
        3) Smoothed Parallel Matching: 3 sweeps \\
        4) Smoothed Parallel Matching: 4 sweeps \\
        5) VBM\\
        \underline{Hypre \textbf{C}oarsening}\\
        1) Falgout \\
        2) HMIS \\
        3) HMIS1 \\
    };
    \matrix (smoothers) [matrix of nodes,below=of label-1-3,xshift=1.1cm,text width=3cm]{%
        \underline{\textbf{S}moother}\\
        1) Hybrid Gauss-Seidel\\
        2) $l_1$--Hybrid Gauss-Seidel\\
        3) INVK\\
        4) $l_1$--INVK \\
        5) $l_1$--Jacobi \\
    };
    \matrix (coarsesolver) [matrix of nodes,below=of label-1-4,xshift=3.4cm,text width=3.5cm]{%
        \underline{\textbf{C}oarsest \textbf{S}olver}\\
        1) Preconditioned CG\\
        \underline{Hypre \textbf{C}oarsest \textbf{S}olver}\\
        2) Direct Solver \\
    };
    \draw[->] (cycle.north) |- (label-1-1.west);
    \draw[->] (aggregation.north) -| (label-1-2.south);
    \draw[->] (smoothers.north) -| (label-1-3.south);
    \draw[->] (coarsesolver.north) |- (label-1-4.east);
    \end{tikzpicture}
\end{table}

\subsection{Performance of the smoothers}
\label{res:smoothers}

In this section we analyze performance results obtained by using the smoothers discussed in Table~\ref{tab:smoothingfactor} on CPU cores of Piz Daint. We compare the different versions of the HGS smoothers and of HINVK smoothers (with $\omega=1$) within the K-cycle coupled with A1 aggregation (see Table~\ref{tab:tested_methods}); we refer to this type of preconditioners as KA1 types.      
We consider weak scalability with $256 \times 10^3$ unknowns (dofs: degrees of freedom) per core and up to $16384$ cores, i.e., we reach an overall number of $\sim 4.2 \times 10^9$ dofs. 

\begin{figure}[htbp]
    \centering
    \subfloat[\label{fig:kpmc3smoothit_pizdaint}]{
%
%
\definecolor{mycolor1}{rgb}{0.00000,0.44700,0.74100}%
\definecolor{mycolor2}{rgb}{0.85000,0.32500,0.09800}%
\definecolor{mycolor3}{rgb}{0.92900,0.69400,0.12500}%
\definecolor{mycolor4}{rgb}{0.49400,0.18400,0.55600}%
\begin{tikzpicture}

\begin{axis}[%
width=0.38\textwidth,
height=2.1in,
at={(2.509in,2.689in)},
scale only axis,
xmode=log,
xmin=1,
xmax=16384,
xtick={    1,     2,     4,     8,    16,    32,    64,   128,   256,   512,  1024,  2048,  4096,  8192, 16384},
xticklabels={    1,     2,     4,     8,    16,    32,    64,   128,   256,   512,  1024,  2048,  4096,  8192, 16384},
x tick label style={rotate=60,anchor=east,font=\footnotesize},
xminorticks=true,
log basis x=2,
xlabel style={font=\color{white!15!black},at={(0.5,-0.12)}},
xlabel={number of MPI cores},
ymin=11,
ymax=18,
yminorticks=true,
ylabel style={font=\color{white!15!black}},
title={Iterations},
axis background/.style={fill=white},
legend columns=2,
legend style={at={(0.5,-0.25)}, anchor=north, legend cell align=left, align=left, draw=none,font=\footnotesize}
]
\addplot [color=mycolor1, line width=4.0pt, mark=o, mark options={solid, mycolor1}]
  table[row sep=crcr]{%
1	12\\
2	12\\
4	12\\
8	12\\
16	13\\
32	13\\
64	12\\
128	14\\
256	14\\
512	13\\
1024	15\\
2048	15\\
4096	16\\
8192	17\\
16384	16\\
};
\addlegendentry{KA1S1CS1}

\addplot [color=mycolor2, line width=4.0pt, mark=asterisk, mark options={solid, mycolor2}]
  table[row sep=crcr]{%
1	15\\
2	16\\
4	12\\
8	15\\
16	17\\
32	13\\
64	15\\
128	14\\
256	13\\
512	15\\
1024	17\\
2048	13\\
4096	15\\
8192	18\\
16384	17\\
};
\addlegendentry{KA1S3CS1}

\addplot [color=mycolor3, line width=4.0pt, mark=square, mark options={solid, mycolor3}]
  table[row sep=crcr]{%
1	12\\
2	12\\
4	12\\
8	12\\
16	13\\
32	14\\
64	13\\
128	15\\
256	15\\
512	14\\
1024	15\\
2048	17\\
4096	17\\
8192	15\\
16384	17\\
};
\addlegendentry{KA1S2CS1}

\addplot [color=mycolor4, line width=4.0pt, mark=diamond, mark options={solid, mycolor4}]
  table[row sep=crcr]{%
1	15\\
2	13\\
4	11\\
8	12\\
16	13\\
32	13\\
64	13\\
128	13\\
256	13\\
512	14\\
1024	14\\
2048	14\\
4096	15\\
8192	14\\
16384	17\\
};
\addlegendentry{KA1S4CS1}

\end{axis}
\end{tikzpicture}%
}
    \subfloat[\label{fig:kpmc3smoothtimexiter_pizdaint}]{
%
%
\definecolor{mycolor1}{rgb}{0.00000,0.44700,0.74100}%
\definecolor{mycolor2}{rgb}{0.85000,0.32500,0.09800}%
\definecolor{mycolor3}{rgb}{0.92900,0.69400,0.12500}%
\definecolor{mycolor4}{rgb}{0.49400,0.18400,0.55600}%
\begin{tikzpicture}

\begin{axis}[%
width=0.38\textwidth,
height=2.1in,
at={(2.509in,2.689in)},
scale only axis,
xmode=log,
xmin=1,
xmax=16384,
xtick={    1,     2,     4,     8,    16,    32,    64,   128,   256,   512,  1024,  2048,  4096,  8192, 16384},
xticklabels={    1,     2,     4,     8,    16,    32,    64,   128,   256,   512,  1024,  2048,  4096,  8192, 16384},
x tick label style={rotate=60,anchor=east,font=\footnotesize},
xminorticks=true,
log basis x=2,
xlabel style={font=\color{white!15!black},at={(0.5,-0.12)}},
xlabel={number of MPI cores},
ymin=0.03,
ymax=0.23,
yminorticks=true,
ylabel style={font=\color{white!15!black}},
y tick label style={
    /pgf/number format/.cd,
        fixed,
        fixed zerofill,
        precision=2,
    /tikz/.cd
},
title={Time per iteration},
axis background/.style={fill=white},
legend columns=2,
legend style={at={(0.5,-0.25)}, anchor=north, legend cell align=left, align=left, draw=none,font=\footnotesize}
]
\addplot [color=mycolor1, line width=4.0pt, mark=o, mark options={solid, mycolor1}]
  table[row sep=crcr]{%
1	0.03\\
2	0.03\\
4	0.04\\
8	0.05\\
16	0.06\\
32	0.08\\
64	0.12\\
128	0.11\\
256	0.11\\
512	0.1\\
1024	0.1\\
2048	0.18\\
4096	0.15\\
8192	0.19\\
16384	0.18\\
};
\addlegendentry{KA1S1CS1}

\addplot [color=mycolor2, line width=4.0pt, mark=asterisk, mark options={solid, mycolor2}]
  table[row sep=crcr]{%
1	0.03\\
2	0.03\\
4	0.04\\
8	0.06\\
16	0.06\\
32	0.09\\
64	0.09\\
128	0.12\\
256	0.12\\
512	0.1\\
1024	0.1\\
2048	0.17\\
4096	0.21\\
8192	0.17\\
16384	0.21\\
};
\addlegendentry{KA1S3CS1}

\addplot [color=mycolor3, line width=4.0pt, mark=square, mark options={solid, mycolor3}]
  table[row sep=crcr]{%
1	0.03\\
2	0.03\\
4	0.04\\
8	0.05\\
16	0.06\\
32	0.08\\
64	0.09\\
128	0.11\\
256	0.11\\
512	0.12\\
1024	0.11\\
2048	0.12\\
4096	0.18\\
8192	0.21\\
16384	0.22\\
};
\addlegendentry{KA1S2CS1}

\addplot [color=mycolor4, line width=4.0pt, mark=diamond, mark options={solid, mycolor4}]
  table[row sep=crcr]{%
1	0.03\\
2	0.03\\
4	0.04\\
8	0.06\\
16	0.07\\
32	0.1\\
64	0.09\\
128	0.12\\
256	0.13\\
512	0.12\\
1024	0.12\\
2048	0.15\\
4096	0.2\\
8192	0.23\\
16384	0.23\\
};
\addlegendentry{KA1S4CS1}

\end{axis}
\end{tikzpicture}%
}
    \caption{Weak scaling results 256k dofs per core. Number of iterations and time per iteration for different smoothers when KA1-type preconditioners is used.}
    \centering
    \subfloat[\label{fig:kpmc3smoothsolvetime_pizdaint}]{
%
%
\definecolor{mycolor1}{rgb}{0.00000,0.44700,0.74100}%
\definecolor{mycolor2}{rgb}{0.85000,0.32500,0.09800}%
\definecolor{mycolor3}{rgb}{0.92900,0.69400,0.12500}%
\definecolor{mycolor4}{rgb}{0.49400,0.18400,0.55600}%
\begin{tikzpicture}

\begin{axis}[%
width=0.41\textwidth,
height=2.1in,
at={(2.509in,2.689in)},
scale only axis,
xmode=log,
xmin=1,
xmax=16384,
xtick={    1,     2,     4,     8,    16,    32,    64,   128,   256,   512,  1024,  2048,  4096,  8192, 16384},
xticklabels={    1,     2,     4,     8,    16,    32,    64,   128,   256,   512,  1024,  2048,  4096,  8192, 16384},
x tick label style={rotate=60,anchor=east,font=\footnotesize},
xminorticks=true,
log basis x=2,
xminorticks=true,
xlabel style={font=\color{white!15!black},at={(0.5,-0.12)}},
xlabel={number of MPI cores},
ymin=0.32,
ymax=3.98,
ylabel style={font=\color{white!15!black}},
axis background/.style={fill=white},
title={Execution Time for Solve (sec.)},
legend columns=2,
legend style={at={(0.5,-0.25)}, anchor=north, legend cell align=left, align=left, draw=none,font=\footnotesize}
]
\addplot [color=mycolor1, line width=4.0pt, mark=o, mark options={solid, mycolor1}]
  table[row sep=crcr]{%
1	0.35\\
2	0.34\\
4	0.45\\
8	0.64\\
16	0.75\\
32	1.06\\
64	1.39\\
128	1.48\\
256	1.49\\
512	1.29\\
1024	1.54\\
2048	2.75\\
4096	2.38\\
8192	3.26\\
16384	2.91\\
};
\addlegendentry{KA1S1CS1}

\addplot [color=mycolor2, line width=4.0pt, mark=asterisk, mark options={solid, mycolor2}]
  table[row sep=crcr]{%
1	0.46\\
2	0.47\\
4	0.51\\
8	0.87\\
16	1.03\\
32	1.2\\
64	1.36\\
128	1.62\\
256	1.6\\
512	1.51\\
1024	1.75\\
2048	2.19\\
4096	3.08\\
8192	3.04\\
16384	3.55\\
};
\addlegendentry{KA1S3CS1}

\addplot [color=mycolor3, line width=4.0pt, mark=square, mark options={solid, mycolor3}]
  table[row sep=crcr]{%
1	0.32\\
2	0.36\\
4	0.47\\
8	0.65\\
16	0.76\\
32	1.16\\
64	1.12\\
128	1.61\\
256	1.59\\
512	1.75\\
1024	1.72\\
2048	2.07\\
4096	3.07\\
8192	3.11\\
16384	3.68\\
};
\addlegendentry{KA1S2CS1}

\addplot [color=mycolor4, line width=4.0pt, mark=diamond, mark options={solid, mycolor4}]
  table[row sep=crcr]{%
1	0.4\\
2	0.4\\
4	0.47\\
8	0.74\\
16	0.85\\
32	1.27\\
64	1.19\\
128	1.58\\
256	1.63\\
512	1.66\\
1024	1.67\\
2048	2.09\\
4096	3.02\\
8192	3.2\\
16384	3.98\\
};
\addlegendentry{KA1S4CS1}

\end{axis}
\end{tikzpicture}%
}
    \subfloat[\label{fig:kpmc3smoothsetuptime_pizdaint}]{
%
%
\definecolor{mycolor1}{rgb}{0.00000,0.44700,0.74100}%
\definecolor{mycolor2}{rgb}{0.85000,0.32500,0.09800}%
\definecolor{mycolor3}{rgb}{0.92900,0.69400,0.12500}%
\definecolor{mycolor4}{rgb}{0.49400,0.18400,0.55600}%
\begin{tikzpicture}

\begin{axis}[%
width=0.41\textwidth,
height=2.1in,
at={(2.209in,2.689in)},
scale only axis,
xmode=log,
xmin=1,
xmax=16384,
xtick={    1,     2,     4,     8,    16,    32,    64,   128,   256,   512,  1024,  2048,  4096,  8192, 16384},
xticklabels={    1,     2,     4,     8,    16,    32,    64,   128,   256,   512,  1024,  2048,  4096,  8192, 16384},
x tick label style={rotate=60,anchor=east,font=\footnotesize},
xminorticks=true,
log basis x=2,
xminorticks=true,
xlabel style={font=\color{white!15!black},at={(0.5,-0.12)}},
xlabel={number of MPI cores},
ymode=log,
ymin=0.69,
ymax=12.9,
yminorticks=true,
ylabel style={font=\color{white!15!black}},
axis background/.style={fill=white},
title={Execution Time for Setup (sec.)},
legend columns = 2,
legend style={at={(0.5,-0.25)}, anchor=north, legend cell align=left, align=left, draw=none,font=\footnotesize}
]
\addplot [color=mycolor1, line width=4.0pt, mark=o, mark options={solid, mycolor1}]
  table[row sep=crcr]{%
1	0.69\\
2	0.83\\
4	1.1\\
8	1.36\\
16	1.66\\
32	2.07\\
64	2.34\\
128	2.68\\
256	2.86\\
512	3.54\\
1024	3.88\\
2048	4.77\\
4096	6.78\\
8192	8.62\\
16384	11.57\\
};
\addlegendentry{KA1S1CS1}

\addplot [color=mycolor2, line width=4.0pt, mark=asterisk, mark options={solid, mycolor2}]
  table[row sep=crcr]{%
1	1.28\\
2	1.42\\
4	1.76\\
8	2.04\\
16	2.35\\
32	2.75\\
64	3.02\\
128	3.39\\
256	3.55\\
512	4.28\\
1024	4.59\\
2048	5.53\\
4096	7.38\\
8192	8.91\\
16384	12.9\\
};
\addlegendentry{KA1S3CS1}

\addplot [color=mycolor3, line width=4.0pt, mark=square, mark options={solid, mycolor3}]
  table[row sep=crcr]{%
1	0.73\\
2	0.89\\
4	1.18\\
8	1.46\\
16	1.76\\
32	2.14\\
64	2.43\\
128	2.83\\
256	3.04\\
512	3.8\\
1024	4.11\\
2048	4.95\\
4096	7.02\\
8192	8.55\\
16384	11.81\\
};
\addlegendentry{KA1S2CS1}

\addplot [color=mycolor4, line width=4.0pt, mark=diamond, mark options={solid, mycolor4}]
  table[row sep=crcr]{%
1	1.29\\
2	1.49\\
4	1.81\\
8	2.09\\
16	2.39\\
32	2.81\\
64	3.12\\
128	3.46\\
256	3.62\\
512	4.42\\
1024	4.71\\
2048	5.51\\
4096	7.5\\
8192	9.14\\
16384	12.61\\
};
\addlegendentry{KA1S4CS1}

\end{axis}
\end{tikzpicture}%
}
    \caption{Weak scaling results 256k dofs per core. Execution times for the solve and setup for different smoothers when KA1-type preconditioners are used.}
\end{figure}

In Fig.~\ref{fig:kpmc3smoothit_pizdaint} we observe that for all smoothers we obtain a good algorithmic scalability, with the number of iterations ranging from $11$ to $18$. The HGS smoother is generally a little better than the $\ell1$-HGS, which sometimes requires one or two additional iterations. For a wide range of number of cores, $\ell1$-HINVK requires fewer iterations than HGS, but when the largest core count is used HGS is better.
In Fig.~\ref{fig:kpmc3smoothtimexiter_pizdaint} we report the solve time per FCG iteration for the different preconditioners. We observe that all preconditioners show a very good implementation scalability with a maximum time per iteration of $0.22$  seconds on $16384$ cores, corresponding to a solve time of $5 \times 10^{-11}$ seconds per dof. 

HGS shows a consistently good behavior with  increasing number of cores which is reflected in the total solve time, as shown in Fig.~\ref{fig:kpmc3smoothsolvetime_pizdaint}. 
In Fig.~\ref{fig:kpmc3smoothsetuptime_pizdaint} we show the setup time for the different preconditioners, including the setup of the AMG hierarchy, the setup of the smoother operators at each level of the hierarchy, and the setup of the coarsest-level solver. The setup of the preconditoners with the HINVK smoothers  requires a somewhat larger time; this is due to the additional (local) computations needed to perform the approximate inversion, whereas the HGS smoothers only require memory copies to generate the matrix splitting. The smoother setup time is completely flat, since it only depends on the size of the local matrix, whereas the hierarchy setup time tends to grow with the number of cores/size of the systems, hence the gap between the curves for the total setup time in Fig.~\ref{fig:kpmc3smoothsetuptime_pizdaint}  tends to close.

In all  cases the setup times show good scalability with a sub-linear increase for increasing number of cores and dofs. We did the same analysis with a more aggressive coarsening using $4$ sweeps of basic pairwise aggregation, i.e., when the size of aggregates is at most $16$, which we denote as KA2-type preconditioners. 
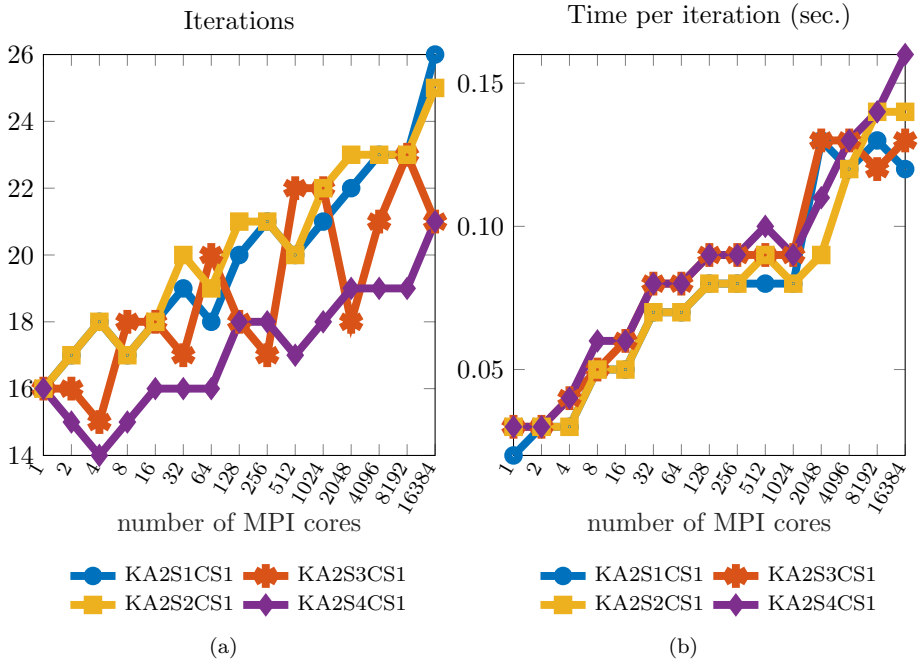
\begin{figure}[htbp]
    \centering
    \subfloat[\label{fig:kpmc4smoothit_pizdaint}]{
%
%
\definecolor{mycolor1}{rgb}{0.00000,0.44700,0.74100}%
\definecolor{mycolor2}{rgb}{0.85000,0.32500,0.09800}%
\definecolor{mycolor3}{rgb}{0.92900,0.69400,0.12500}%
\definecolor{mycolor4}{rgb}{0.49400,0.18400,0.55600}%
\begin{tikzpicture}

\begin{axis}[%
width=0.40\textwidth,
height=2.1in,
at={(2.509in,2.689in)},
scale only axis,
xmode=log,
xmin=1,
xmax=16384,
xtick={    1,     2,     4,     8,    16,    32,    64,   128,   256,   512,  1024,  2048,  4096,  8192, 16384},
xticklabels={    1,     2,     4,     8,    16,    32,    64,   128,   256,   512,  1024,  2048,  4096,  8192, 16384},
x tick label style={rotate=60,anchor=east,font=\footnotesize},
xminorticks=true,
log basis x=2,
xlabel style={font=\color{white!15!black},at={(0.5,-0.12)}},
xlabel={number of MPI cores},
ymin=14,
ymax=26,
yminorticks=true,
axis background/.style={fill=white},
title={Iterations},
legend columns = 2,
legend style={at={(0.5,-0.25)}, anchor=north, legend cell align=left, align=left, draw=none,font=\footnotesize}
]
\addplot [color=mycolor1, line width=3.0pt, mark=o, mark options={solid, mycolor1}]
  table[row sep=crcr]{%
1	16\\
2	17\\
4	18\\
8	17\\
16	18\\
32	19\\
64	18\\
128	20\\
256	21\\
512	20\\
1024	21\\
2048	22\\
4096	23\\
8192	23\\
16384	26\\
};
\addlegendentry{KA2S1CS1}

\addplot [color=mycolor2, line width=3.0pt, mark=asterisk, mark size=4.5pt, mark options={solid, mycolor2}]
  table[row sep=crcr]{%
1	16\\
2	16\\
4	15\\
8	18\\
16	18\\
32	17\\
64	20\\
128	18\\
256	17\\
512	22\\
1024	22\\
2048	18\\
4096	21\\
8192	23\\
16384	21\\
};
\addlegendentry{KA2S3CS1}

\addplot [color=mycolor3, line width=3.0pt, mark=square, mark options={solid, mycolor3}]
  table[row sep=crcr]{%
1	16\\
2	17\\
4	18\\
8	17\\
16	18\\
32	20\\
64	19\\
128	21\\
256	21\\
512	20\\
1024	22\\
2048	23\\
4096	23\\
8192	23\\
16384	25\\
};
\addlegendentry{KA2S2CS1}

\addplot [color=mycolor4, line width=3.0pt, mark=diamond, mark options={solid, mycolor4}]
  table[row sep=crcr]{%
1	16\\
2	15\\
4	14\\
8	15\\
16	16\\
32	16\\
64	16\\
128	18\\
256	18\\
512	17\\
1024	18\\
2048	19\\
4096	19\\
8192	19\\
16384	21\\
};
\addlegendentry{KA2S4CS1}

\end{axis}
\end{tikzpicture}%
}
    \subfloat[\label{fig:kpmc4smoothtimexiter_pizdaint}]{
%
%
\definecolor{mycolor1}{rgb}{0.00000,0.44700,0.74100}%
\definecolor{mycolor2}{rgb}{0.85000,0.32500,0.09800}%
\definecolor{mycolor3}{rgb}{0.92900,0.69400,0.12500}%
\definecolor{mycolor4}{rgb}{0.49400,0.18400,0.55600}%
\begin{tikzpicture}

\begin{axis}[%
width=0.40\textwidth,
height=2.1in,
at={(2.509in,2.689in)},
scale only axis,
xmode=log,
xmin=1,
xmax=16384,
xtick={    1,     2,     4,     8,    16,    32,    64,   128,   256,   512,  1024,  2048,  4096,  8192, 16384},
xticklabels={    1,     2,     4,     8,    16,    32,    64,   128,   256,   512,  1024,  2048,  4096,  8192, 16384},
x tick label style={rotate=60,anchor=east,font=\footnotesize},
xminorticks=true,
log basis x=2,
xlabel style={font=\color{white!15!black},at={(0.5,-0.12)}},
xlabel={number of MPI cores},
ymin=0.02,
ymax=0.16,
yminorticks=true,
ylabel style={font=\color{white!15!black}},
y tick label style={
    /pgf/number format/.cd,
        fixed,
        fixed zerofill,
        precision=2,
    /tikz/.cd
},
axis background/.style={fill=white},
title={Time per iteration (sec.)},
legend columns = 2,
legend style={at={(0.5,-0.25)}, anchor=north, legend cell align=left, align=left, draw=none,font=\footnotesize}
]
\addplot [color=mycolor1, line width=3.0pt, mark=o, mark options={solid, mycolor1}]
  table[row sep=crcr]{%
1	0.02\\
2	0.03\\
4	0.03\\
8	0.05\\
16	0.05\\
32	0.07\\
64	0.07\\
128	0.08\\
256	0.08\\
512	0.08\\
1024	0.08\\
2048	0.13\\
4096	0.12\\
8192	0.13\\
16384	0.12\\
};
\addlegendentry{KA2S1CS1}

\addplot [color=mycolor2, line width=3.0pt, mark=asterisk, mark size=4.5pt, mark options={solid, mycolor2}]
  table[row sep=crcr]{%
1	0.03\\
2	0.03\\
4	0.04\\
8	0.05\\
16	0.06\\
32	0.08\\
64	0.08\\
128	0.09\\
256	0.09\\
512	0.09\\
1024	0.09\\
2048	0.13\\
4096	0.13\\
8192	0.12\\
16384	0.13\\
};
\addlegendentry{KA2S3CS1}

\addplot [color=mycolor3, line width=3.0pt, mark=square, mark options={solid, mycolor3}]
  table[row sep=crcr]{%
1	0.03\\
2	0.03\\
4	0.03\\
8	0.05\\
16	0.05\\
32	0.07\\
64	0.07\\
128	0.08\\
256	0.08\\
512	0.09\\
1024	0.08\\
2048	0.09\\
4096	0.12\\
8192	0.14\\
16384	0.14\\
};
\addlegendentry{KA2S2CS1}

\addplot [color=mycolor4, line width=3.0pt, mark=diamond, mark options={solid, mycolor4}]
  table[row sep=crcr]{%
1	0.03\\
2	0.03\\
4	0.04\\
8	0.06\\
16	0.06\\
32	0.08\\
64	0.08\\
128	0.09\\
256	0.09\\
512	0.1\\
1024	0.09\\
2048	0.11\\
4096	0.13\\
8192	0.14\\
16384	0.16\\
};
\addlegendentry{KA2S4CS1}

\end{axis}
\end{tikzpicture}%
}
    \caption{Weak scaling results 256k dofs per core. Number of iterations and time per iteration for different smoothers when KA2-type preconditioners are used.}
\end{figure}
\begin{figure}[htbp]
    \centering
    \subfloat[\label{fig:kpmc4smoothsolvetime_pizdaint}]{
%
%
\definecolor{mycolor1}{rgb}{0.00000,0.44700,0.74100}%
\definecolor{mycolor2}{rgb}{0.85000,0.32500,0.09800}%
\definecolor{mycolor3}{rgb}{0.92900,0.69400,0.12500}%
\definecolor{mycolor4}{rgb}{0.49400,0.18400,0.55600}%
\begin{tikzpicture}

\begin{axis}[%
width=0.41\textwidth,
height=2.1in,
at={(2.509in,2.689in)},
scale only axis,
xmode=log,
xmin=1,
xmax=16384,
xtick={    1,     2,     4,     8,    16,    32,    64,   128,   256,   512,  1024,  2048,  4096,  8192, 16384},
xticklabels={    1,     2,     4,     8,    16,    32,    64,   128,   256,   512,  1024,  2048,  4096,  8192, 16384},
x tick label style={rotate=60,anchor=east,font=\footnotesize},
xminorticks=true,
log basis x=2,
xminorticks=true,
xlabel style={font=\color{white!15!black},at={(0.5,-0.12)}},
xlabel={number of MPI cores},
ymin=0,
ymax=3.38,
axis background/.style={fill=white},
title={Execution Time for Solve (sec.)},
legend columns = 2,
legend style={at={(0.5,-0.25)}, anchor=north, legend cell align=left, align=left, draw=none,font=\footnotesize}
]
\addplot [color=mycolor1, line width=3.0pt, mark=o, mark options={solid, mycolor1}]
  table[row sep=crcr]{%
1	0.37\\
2	0.43\\
4	0.58\\
8	0.79\\
16	0.88\\
32	1.25\\
64	1.32\\
128	1.6\\
256	1.67\\
512	1.62\\
1024	1.71\\
2048	2.77\\
4096	2.83\\
8192	3.1\\
16384	3.12\\
};
\addlegendentry{KA2S1CS1}

\addplot [color=mycolor2, line width=3.0pt, mark=asterisk, mark size=4.5pt, mark options={solid, mycolor2}]
  table[row sep=crcr]{%
1	0.4\\
2	0.44\\
4	0.53\\
8	0.95\\
16	0.99\\
32	1.32\\
64	1.58\\
128	1.7\\
256	1.59\\
512	1.91\\
1024	1.94\\
2048	2.32\\
4096	2.66\\
8192	2.84\\
16384	2.82\\
};
\addlegendentry{KA2S3CS1}

\addplot [color=mycolor3, line width=3.0pt, mark=square, mark options={solid, mycolor3}]
  table[row sep=crcr]{%
1	0.4\\
2	0.44\\
4	0.59\\
8	0.81\\
16	0.89\\
32	1.33\\
64	1.32\\
128	1.69\\
256	1.74\\
512	1.82\\
1024	1.83\\
2048	2.07\\
4096	2.66\\
8192	3.18\\
16384	3.38\\
};
\addlegendentry{KA2S2CS1}

\addplot [color=mycolor4, line width=3.0pt, mark=diamond, mark options={solid, mycolor4}]
  table[row sep=crcr]{%
1	0.41\\
2	0.43\\
4	0.51\\
8	0.83\\
16	0.93\\
32	1.26\\
64	1.29\\
128	1.71\\
256	1.69\\
512	1.73\\
1024	1.7\\
2048	2.02\\
4096	2.43\\
8192	2.59\\
16384	3.27\\
};
\addlegendentry{KA2S4CS1}

\end{axis}
\end{tikzpicture}%
}
    \subfloat[\label{fig:kpmc4smoothsetuptime_pizdaint}]{
%
%
\definecolor{mycolor1}{rgb}{0.00000,0.44700,0.74100}%
\definecolor{mycolor2}{rgb}{0.85000,0.32500,0.09800}%
\definecolor{mycolor3}{rgb}{0.92900,0.69400,0.12500}%
\definecolor{mycolor4}{rgb}{0.49400,0.18400,0.55600}%
\begin{tikzpicture}

\begin{axis}[%
width=0.41\textwidth,
height=2.1in,
at={(2.509in,2.689in)},
scale only axis,
xmode=log,
xmin=1,
xmax=16384,
xtick={    1,     2,     4,     8,    16,    32,    64,   128,   256,   512,  1024,  2048,  4096,  8192, 16384},
xticklabels={    1,     2,     4,     8,    16,    32,    64,   128,   256,   512,  1024,  2048,  4096,  8192, 16384},
x tick label style={rotate=60,anchor=east,font=\footnotesize},
xminorticks=true,
log basis x=2,
xminorticks=true,
xlabel style={font=\color{white!15!black},at={(0.5,-0.12)}},
xlabel={number of MPI cores},
ymode=log,
ymin=0.67,
ymax=12.08,
yminorticks=true,
axis background/.style={fill=white},
title={Execution Time for Setup (sec.)},
legend columns = 2,
legend style={at={(0.5,-0.25)}, anchor=north, legend cell align=left, align=left, draw=none,font=\footnotesize}
]
\addplot [color=mycolor1, line width=3.0pt, mark=o, mark options={solid, mycolor1}]
  table[row sep=crcr]{%
1	0.67\\
2	0.8\\
4	1.09\\
8	1.33\\
16	1.63\\
32	2.02\\
64	2.34\\
128	2.63\\
256	2.79\\
512	3.48\\
1024	3.81\\
2048	4.64\\
4096	6.34\\
8192	8.21\\
16384	11.14\\
};
\addlegendentry{KA2S1CS1}

\addplot [color=mycolor2, line width=3.0pt, mark=asterisk, mark size=4.5pt, mark options={solid, mycolor2}]
  table[row sep=crcr]{%
1	1.23\\
2	1.38\\
4	1.69\\
8	1.97\\
16	2.28\\
32	2.67\\
64	3\\
128	3.33\\
256	3.43\\
512	4.16\\
1024	4.47\\
2048	5.42\\
4096	7.16\\
8192	8.92\\
16384	11.97\\
};
\addlegendentry{KA2S3CS1}

\addplot [color=mycolor3, line width=3.0pt, mark=square, mark options={solid, mycolor3}]
  table[row sep=crcr]{%
1	0.72\\
2	0.86\\
4	1.14\\
8	1.42\\
16	1.72\\
32	2.13\\
64	2.43\\
128	2.76\\
256	2.92\\
512	3.74\\
1024	3.97\\
2048	4.63\\
4096	6.69\\
8192	8.01\\
16384	11.46\\
};
\addlegendentry{KA2S2CS1}

\addplot [color=mycolor4, line width=3.0pt, mark=diamond, mark options={solid, mycolor4}]
  table[row sep=crcr]{%
1	1.24\\
2	1.39\\
4	1.73\\
8	2.01\\
16	2.32\\
32	2.73\\
64	2.99\\
128	3.36\\
256	3.52\\
512	4.31\\
1024	4.64\\
2048	5.29\\
4096	7.18\\
8192	8.73\\
16384	12.08\\
};
\addlegendentry{KA2S4CS1}

\end{axis}
\end{tikzpicture}%
}
    \caption{Weak scaling results 256k dofs per core. Execution times for the solve and setup for different smoothers when KA2-type preconditioners are used.}
\end{figure}
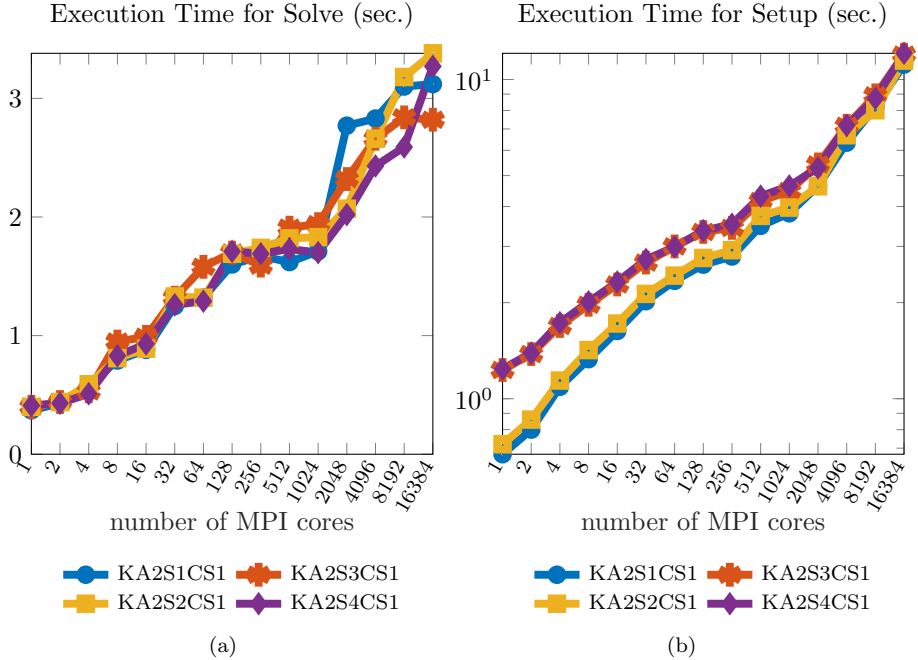
In Fig.~\ref{fig:kpmc4smoothit_pizdaint}-\ref{fig:kpmc4smoothtimexiter_pizdaint} we report number of iterations and solve time per iterations, respectively, when the different preconditioners are applied. As expected, in this case the number of iterations is generally larger than the case of KA1-type preconditioners, however we still get a fairly good numerical scalability, with a moderate increase in the number of iterations for increasing number of cores. The best convergence behavior is obtained when HINVK smoothers are employed, while the minimum number of iterations is generally obtained by $\ell1$-HINVK. This good convergence behavior balances the small increase in solve time per iteration required by the HINVK smoothers, resulting in a total solve time which is generally better when $\ell1$-HINVK is employed, except on $16384$ cores, where HINVK and HGS obtain a slightly better solve time (see Fig.~\ref{fig:kpmc4smoothsolvetime_pizdaint}). For the sake of completeness we report in Fig.~\ref{fig:kpmc4smoothsetuptime_pizdaint} the setup time for the different preconditioners based on the more aggressive coarsening, showing a behavior similar to that  of the KA1-type preconditioners.

\subsection{Comparison with \texttt{Hypre}}
\label{res:hypre}

We now  compare our preconditioners with some of those available in the \texttt{Hypre} library~\cite{FJY2006}, using default algorithmic parameters for best practice. Specifically, we compare with three different coarsening approaches, i.e., {Falgout}~\cite{MR1908755},~\cite[section~3.2]{MR2267933}, {HMIS}~\cite{MR2205610}, and {HMIS1}~\cite[section~3]{MR2650221}; we refer to section~\ref{sec:hyprecoarsening} in the supplementary materials for a more extensive discussion of their features. All these coarsening schemes are used to generate an AMG hierarchy which is applied as a V-cycle preconditioner for a CG method, with one sweep of forward/backward HGS being applied as pre/post-smoother. Default choices are used for the coarsest system, where a direct method is employed; using the same notation applied for our preconditioners, the three \texttt{Hypre} preconditioners are denoted, respectively, VC1S1CS2, VC2S1CS2, and VC3S1CS2.

{We compare these preconditioners from \texttt{Hypre} with the KA1S1CS1 and VA3S1CS1 methods in \texttt{AMG4PSBLAS}
and with the VA5S1CS1 preconditioner, based on the parallel decoupled version of the smoothed aggregation strategy from \cite{VMB1996,TT2000}, implemented in the previous version of the  library~\cite{BDDF2007,DDF2010}, i.e., \texttt{MLD2P4} v. 2.2.}

{For these experiments we used the Piz Daint machine up to 8192 CPU cores.}
\begin{figure}[htbp]
    \centering
%
%
\definecolor{mycolor1}{rgb}{0.00000,0.44700,0.74100}%
\definecolor{mycolor2}{rgb}{0.85000,0.32500,0.09800}%
\definecolor{mycolor3}{rgb}{0.92900,0.69400,0.12500}%
\definecolor{mycolor4}{rgb}{0.49400,0.18400,0.55600}%
\definecolor{mycolor5}{rgb}{0.46600,0.67400,0.18800}%
\definecolor{mycolor6}{rgb}{0.30100,0.74500,0.93300}%
\begin{tikzpicture}

\begin{axis}[%
width=0.41\textwidth,
height=2.1in,
at={(0in,0in)},
scale only axis,
xmode=log,
xmin=1,
xmax=8192,
xtick={    1,     2,     4,     8,    16,    32,    64,   128,   256,   512,  1024,  2048,  4096,  8192},
xticklabels={    1,     2,     4,     8,    16,    32,    64,   128,   256,   512,  1024,  2048,  4096,  8192},
x tick label style={rotate=60,anchor=east,font=\footnotesize},
xminorticks=true,
log basis x=2,
xminorticks=true,
xlabel style={font=\color{white!15!black},at={(0.5,-0.12)}},
xlabel={number of MPI cores},
ymin=0,
ymax=5.186356,
yminorticks=true,
ylabel style={font=\color{white!15!black}\small},
title={Operator Complexity},
axis background/.style={fill=white},
legend columns=3, 
legend style={at={(1.0,-0.2)}, anchor=north, legend cell align=left, align=left, draw=none,font=\small}
]
\addplot [color=mycolor1, dashdotted, line width=3.0pt]
  table[row sep=crcr]{%
1	3.574353\\
2	3.868418\\
4	4.537781\\
8	4.10222\\
16	4.776941\\
32	4.958702\\
64	4.322591\\
128	4.977869\\
256	5.120193\\
512	4.449998\\
1024	4.519012\\
2048	5.184155\\
4096	5.13489\\
8192	4.566168\\
};
\addlegendentry{VC1S1CS2}

\addplot [color=mycolor2, dashed, line width=3.0pt]
  table[row sep=crcr]{%
1	3.118791\\
2	3.104061\\
4	3.078508\\
8	3.11332\\
16	3.079452\\
32	3.048086\\
64	3.118852\\
128	3.082773\\
256	3.075108\\
512	3.121879\\
1024	3.107453\\
2048	3.074346\\
4096	3.089975\\
8192	3.106731\\
};
\addlegendentry{VC2S1CS2}

\addplot [color=mycolor3, dotted, line width=3.0pt, mark size=2.5pt, mark=+, mark options={solid, mycolor3}]
  table[row sep=crcr]{%
1	1.350973\\
2	1.334193\\
4	1.324674\\
8	1.344212\\
16	1.320715\\
32	1.31974\\
64	1.343717\\
128	1.319893\\
256	1.321543\\
512	1.343195\\
1024	1.331307\\
2048	1.321115\\
4096	1.328164\\
8192	1.330828\\
};
\addlegendentry{VC3S1CS2}

\addplot [color=mycolor4, dotted, line width=3.0pt, mark size=2.5pt, mark=+, mark options={solid, mycolor4}]
  table[row sep=crcr]{%
1	1.14\\
2	1.141\\
4	1.142\\
8	1.142\\
16	1.142\\
32	1.143\\
64	1.142\\
128	1.143\\
256	1.143\\
512	1.143\\
1024	1.143\\
2048	1.144\\
4096	1.143\\
8192	1.143\\
};
\addlegendentry{KA1S1CS1}

\addplot [color=mycolor5, dotted, line width=3.0pt, mark size=2.5pt, mark=+, mark options={solid, mycolor5}]
  table[row sep=crcr]{%
1	1.82\\
2	1.846\\
4	1.873\\
8	1.889\\
16	1.905\\
32	1.92\\
64	1.93\\
128	1.938\\
256	1.947\\
512	1.951\\
1024	1.956\\
2048	1.961\\
4096	1.962\\
8192	1.965\\
};
\addlegendentry{VA3S1CS1}

\addplot [color=mycolor6, dotted, line width=3.0pt, mark size=2.5pt, mark=+, mark options={solid, mycolor6}]
  table[row sep=crcr]{%
1	1.555\\
2	1.568\\
4	1.58\\
8	1.58\\
16	1.587\\
32	1.594\\
64	1.593\\
128	1.597\\
256	1.603\\
512	1.6\\
1024	1.601\\
2048	1.607\\
4096	1.603\\
8192	1.604\\
};
\addlegendentry{VA5S1CS1}

\end{axis}

\begin{axis}[%
width=0.41\textwidth,
height=2.1in,
at={(2.4in,0in)},
scale only axis,
xmode=log,
log basis x=2,
xmin=1,
xmax=8192,
xtick={    1,     2,     4,     8,    16,    32,    64,   128,   256,   512,  1024,  2048,  4096,  8192},
xticklabels={    1,     2,     4,     8,    16,    32,    64,   128,   256,   512,  1024,  2048,  4096,  8192},
x tick label style={rotate=60,anchor=east,font=\footnotesize},
xminorticks=true,
log basis x=2,
xminorticks=true,
xlabel style={font=\color{white!15!black},at={(0.5,-0.12)}},
xlabel={number of MPI cores},
ymin=0,
ymax=45,
yminorticks=true,
ylabel style={font=\color{white!15!black}},
title={Iterations},
axis background/.style={fill=white},
]
\addplot [color=mycolor1, dashdotted, line width=4.0pt]
  table[row sep=crcr]{%
1	6\\
2	7\\
4	7\\
8	8\\
16	8\\
32	8\\
64	9\\
128	10\\
256	10\\
512	11\\
1024	12\\
2048	13\\
4096	13\\
8192	14\\
};

\addplot [color=mycolor2, dashed, line width=3.0pt]
  table[row sep=crcr]{%
1	6\\
2	9\\
4	12\\
8	13\\
16	14\\
32	14\\
64	16\\
128	18\\
256	20\\
512	22\\
1024	25\\
2048	29\\
4096	35\\
8192	37\\
};

\addplot [color=mycolor3, dotted, line width=3.0pt, mark size=2.5pt, mark=+, mark options={solid, mycolor3}]
  table[row sep=crcr]{%
1	12\\
2	15\\
4	17\\
8	17\\
16	20\\
32	20\\
64	22\\
128	25\\
256	27\\
512	29\\
1024	32\\
2048	37\\
4096	40\\
8192	44\\
};

\addplot [color=mycolor4, dotted, line width=3.0pt, mark size=2.5pt, mark=+, mark options={solid, mycolor4}]
  table[row sep=crcr]{%
1	12\\
2	12\\
4	12\\
8	12\\
16	12\\
32	13\\
64	12\\
128	12\\
256	13\\
512	14\\
1024	15\\
2048	16\\
4096	15\\
8192	16\\
};

\addplot [color=mycolor5, dotted, line width=3.0pt, mark size=2.5pt, mark=+, mark options={solid, mycolor5}]
  table[row sep=crcr]{%
1	7\\
2	7\\
4	7\\
8	7\\
16	8\\
32	9\\
64	8\\
128	8\\
256	9\\
512	8\\
1024	8\\
2048	12\\
4096	8\\
8192	9\\
};

\addplot [color=mycolor6, dotted, line width=3.0pt, mark size=2.5pt, mark=+, mark options={solid, mycolor6}]
  table[row sep=crcr]{%
1	11\\
2	12\\
4	13\\
8	14\\
16	14\\
32	15\\
64	15\\
128	15\\
256	16\\
512	18\\
1024	17\\
2048	21\\
4096	26\\
8192	27\\
};

\end{axis}

\end{tikzpicture}%
    \caption{Comparison with \texttt{Hypre} 256k dofs per core. Operator complexity and number of iterations for different preconditioners.}
    \label{fig:hyprevsamg4psblas_cost}
\end{figure}
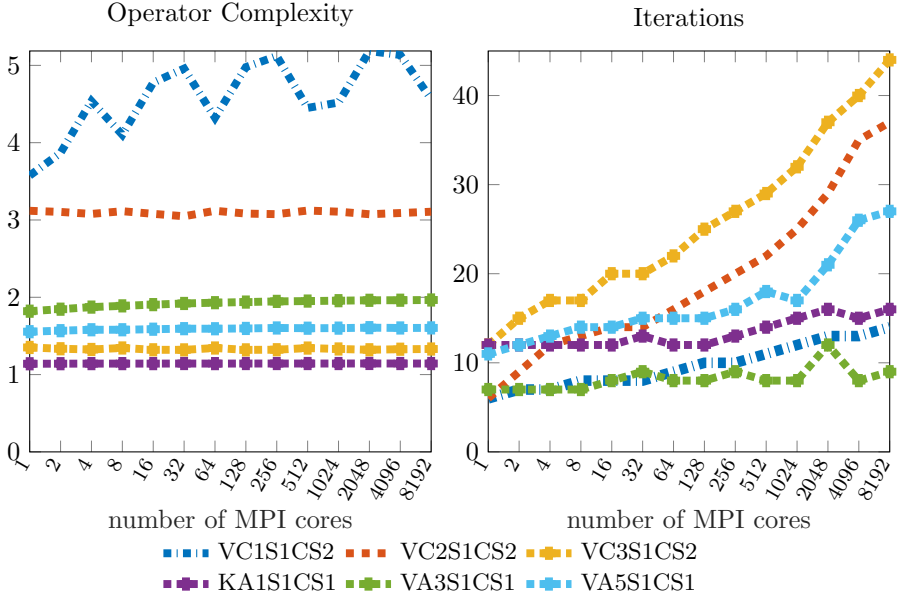
{We detail first, in Fig.~\ref{fig:hyprevsamg4psblas_cost}, the operator complexity obtained by the different coarsening approaches; we can observe that the new strategy based on the A1 aggregation consistently produces the smallest operator complexity, thus reducing the memory footprint for the preconditioner. On the other hand, the smoothed version A3 produces a small increase in the operator complexity, whilst still  remaining much smaller than VC1S1CS2 and VC2S1CS2 and smaller than $2$, therefore ensuring a linear cost in the application of the preconditioner. Note that our VA3S1CS1 shows the best algorithmic scalability, requiring a number of iterations which is almost constant (between $7$ and $9$) for increasing number of MPI cores except on $2048$ where a small increase to $12$ is observed. This is essentially due to a small decrease in the ability of obtaining a matching of maximum cardinality,  hence producing a small deterioration of the quality of the final aggregates at the coarser levels during the hierarchy setup. This aspect is currently under investigation: we are studying  new parallel approximate matching algorithms aiming to realize the best trade-off between maximum weight and cardinality in a way that is independent of both problem size and data partitioning.  
VC1S1CS2 shows a comparable number of iterations ranging from $6$ to $14$: this is the benefit accrued by the large operator complexity of the corresponding AMG hierarchy. On the other hand, KA1S1CS1 also shows a very good algorithmic scalability, with number of iterations ranging from $12$ to $16$, especially considering the fact that both  operator complexity and memory requirements are about $5$ times smaller than those of VC1S1CS2. All the other variants show a more significant increase in the number of iterations with respect to the problem size.} 
{To complement this information we look also at the setup time for all the preconditioners. 
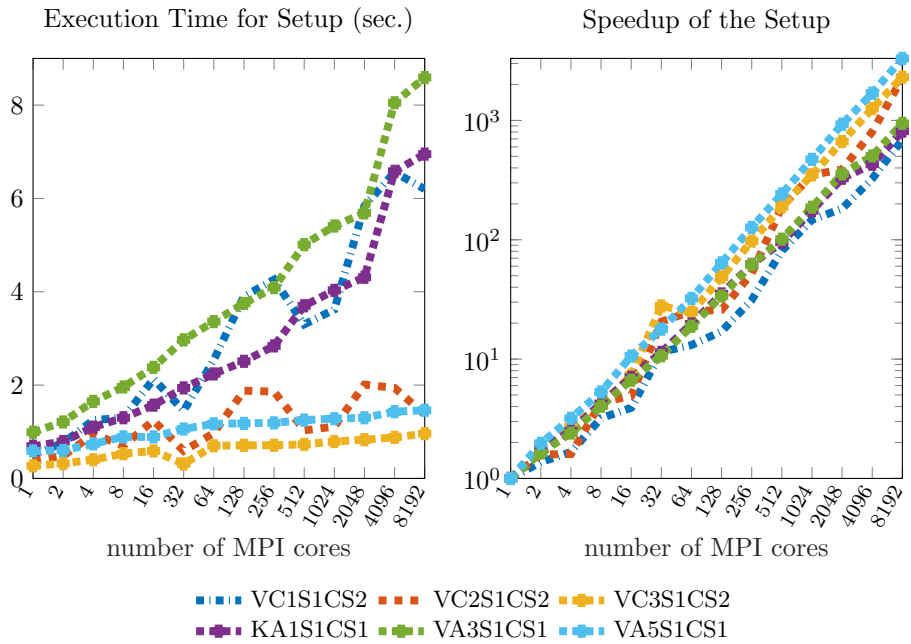
\begin{figure}[htbp]
    \centering
%
%
\definecolor{mycolor1}{rgb}{0.00000,0.44700,0.74100}%
\definecolor{mycolor2}{rgb}{0.85000,0.32500,0.09800}%
\definecolor{mycolor3}{rgb}{0.92900,0.69400,0.12500}%
\definecolor{mycolor4}{rgb}{0.49400,0.18400,0.55600}%
\definecolor{mycolor5}{rgb}{0.46600,0.67400,0.18800}%
\definecolor{mycolor6}{rgb}{0.30100,0.74500,0.93300}%
\begin{tikzpicture}

\begin{axis}[%
width=0.40\textwidth,
height=2.20in,
at={(0in,0in)},
scale only axis,
xmode=log,
xmin=1,
xmax=8192,
xtick={    1,     2,     4,     8,    16,    32,    64,   128,   256,   512,  1024,  2048,  4096,  8192},
xticklabels={    1,     2,     4,     8,    16,    32,    64,   128,   256,   512,  1024,  2048,  4096,  8192},
x tick label style={rotate=60,anchor=east,font=\footnotesize},
xminorticks=true,
log basis x=2,
xminorticks=true,
xlabel style={font=\color{white!15!black},at={(0.5,-0.12)}},
xlabel={number of MPI cores},
log basis x=2,
ymin=0,
ymax=9,
yminorticks=true,
axis background/.style={fill=white},
title={Execution Time for Setup (sec.)},
legend columns=3, 
legend style={at={(1.1,-0.24)}, anchor=north, legend cell align=left, align=left,font=\small,draw=none}
]
\addplot [color=mycolor1, dashdotted, line width=3.0pt]
  table[row sep=crcr]{%
1	0.52\\
2	0.76\\
4	1.25\\
8	1.3\\
16	2.13\\
32	1.45\\
64	2.54\\
128	3.9\\
256	4.26\\
512	3.3\\
1024	3.62\\
2048	5.84\\
4096	6.56\\
8192	6.2\\
};
\addlegendentry{VC1S1CS2}

\addplot [color=mycolor2, dashed, line width=3.0pt]
  table[row sep=crcr]{%
1	0.38\\
2	0.48\\
4	0.95\\
8	0.7\\
16	1.24\\
32	0.59\\
64	0.97\\
128	1.88\\
256	1.86\\
512	1.03\\
1024	1.11\\
2048	2.01\\
4096	1.94\\
8192	1.37\\
};
\addlegendentry{VC2S1CS2}

\addplot [color=mycolor3, dotted, line width=3.0pt, mark size=2.5pt, mark=+, mark options={solid, mycolor3}]
  table[row sep=crcr]{%
1	0.27\\
2	0.32\\
4	0.4\\
8	0.53\\
16	0.59\\
32	0.31\\
64	0.7\\
128	0.71\\
256	0.71\\
512	0.73\\
1024	0.79\\
2048	0.83\\
4096	0.88\\
8192	0.96\\
};
\addlegendentry{VC3S1CS2}

\addplot [color=mycolor4, dotted, line width=3.0pt, mark size=2.5pt, mark=+, mark options={solid, mycolor4}]
  table[row sep=crcr]{%
1	0.69\\
2	0.79\\
4	1.1\\
8	1.3\\
16	1.57\\
32	1.94\\
64	2.24\\
128	2.51\\
256	2.84\\
512	3.7\\
1024	4.03\\
2048	4.31\\
4096	6.58\\
8192	6.95\\
};
\addlegendentry{KA1S1CS1}

\addplot [color=mycolor5, dotted, line width=3.0pt, mark size=2.5pt, mark=+, mark options={solid, mycolor5}]
  table[row sep=crcr]{%
1	0.99\\
2	1.21\\
4	1.65\\
8	1.96\\
16	2.38\\
32	2.97\\
64	3.36\\
128	3.75\\
256	4.09\\
512	5.01\\
1024	5.41\\
2048	5.69\\
4096	8.05\\
8192	8.59\\
};
\addlegendentry{VA3S1CS1}

\addplot [color=mycolor6, dotted, line width=3.0pt, mark size=2.5pt, mark=+, mark options={solid, mycolor6}]
  table[row sep=crcr]{%
1	0.59\\
2	0.6\\
4	0.74\\
8	0.89\\
16	0.89\\
32	1.06\\
64	1.17\\
128	1.18\\
256	1.19\\
512	1.25\\
1024	1.28\\
2048	1.3\\
4096	1.43\\
8192	1.46\\
};
\addlegendentry{VA5S1CS1}

\end{axis}
\begin{axis}[%
width=0.40\textwidth,
height=2.20in,
at={(2.5in,0in)},
scale only axis,
xmode=log,
xmin=1,
xmax=8192,
xtick={    1,     2,     4,     8,    16,    32,    64,   128,   256,   512,  1024,  2048,  4096,  8192},
xticklabels={    1,     2,     4,     8,    16,    32,    64,   128,   256,   512,  1024,  2048,  4096,  8192},
x tick label style={rotate=60,anchor=east,font=\footnotesize},
xminorticks=true,
log basis x=2,
xminorticks=true,
xlabel style={font=\color{white!15!black},at={(0.5,-0.12)}},
xlabel={number of MPI cores},
ymode=log,
ymin=1,
ymax=3310.46575342466,
yminorticks=true,
axis background/.style={fill=white},
title={Speedup of the Setup},
]
\addplot [color=mycolor1, dashdotted, line width=3.0pt]
  table[row sep=crcr]{%
1	1\\
2	1.36842105263158\\
4	1.664\\
8	3.2\\
16	3.90610328638498\\
32	11.4758620689655\\
64	13.1023622047244\\
128	17.0666666666667\\
256	31.2488262910798\\
512	80.6787878787879\\
1024	147.093922651934\\
2048	182.356164383562\\
4096	324.682926829268\\
8192	687.070967741936\\
};

\addplot [color=mycolor2, dashed, line width=3.0pt]
  table[row sep=crcr]{%
1	1\\
2	1.58333333333333\\
4	1.6\\
8	4.34285714285714\\
16	4.90322580645161\\
32	20.6101694915254\\
64	25.0721649484536\\
128	25.8723404255319\\
256	52.3010752688172\\
512	188.893203883495\\
1024	350.558558558559\\
2048	387.18407960199\\
4096	802.309278350515\\
8192	2272.23357664234\\
};

\addplot [color=mycolor3, dotted, line width=3.0pt, mark size=2.5pt, mark=+, mark options={solid, mycolor3}]
  table[row sep=crcr]{%
1	1\\
2	1.6875\\
4	2.7\\
8	4.07547169811321\\
16	7.32203389830509\\
32	27.8709677419355\\
64	24.6857142857143\\
128	48.6760563380282\\
256	97.3521126760564\\
512	189.369863013699\\
1024	349.974683544304\\
2048	666.21686746988\\
4096	1256.72727272727\\
8192	2304\\
};

\addplot [color=mycolor4, dotted, line width=3.0pt, mark size=2.5pt, mark=+, mark options={solid, mycolor4}]
  table[row sep=crcr]{%
1	1\\
2	1.74683544303797\\
4	2.50909090909091\\
8	4.24615384615385\\
16	7.03184713375796\\
32	11.3814432989691\\
64	19.7142857142857\\
128	35.1872509960159\\
256	62.1971830985915\\
512	95.4810810810811\\
1024	175.325062034739\\
2048	327.870069605568\\
4096	429.519756838906\\
8192	813.306474820144\\
};

\addplot [color=mycolor5, dotted, line width=3.0pt, mark size=2.5pt, mark=+, mark options={solid, mycolor5}]
  table[row sep=crcr]{%
1	1\\
2	1.63636363636364\\
4	2.4\\
8	4.04081632653061\\
16	6.65546218487395\\
32	10.6666666666667\\
64	18.8571428571429\\
128	33.792\\
256	61.9657701711491\\
512	101.173652694611\\
1024	187.386321626617\\
2048	356.330404217926\\
4096	503.731677018633\\
8192	944.130384167637\\
};

\addplot [color=mycolor6, dotted, line width=3.0pt, mark size=2.5pt, mark=+, mark options={solid, mycolor6}]
  table[row sep=crcr]{%
1	1\\
2	1.96666666666667\\
4	3.18918918918919\\
8	5.30337078651685\\
16	10.6067415730337\\
32	17.811320754717\\
64	32.2735042735043\\
128	64\\
256	126.924369747899\\
512	241.664\\
1024	472\\
2048	929.476923076923\\
4096	1689.95804195804\\
8192	3310.46575342466\\
};

\end{axis}
\end{tikzpicture}%
    \caption{Comparison with \texttt{Hypre} 256k dofs per core. Preconditioners setup: execution time (left), speedup (right).}
    \label{fig:hyprevsamg4psblas_setup}
\end{figure}
From the results in Fig.~\ref{fig:hyprevsamg4psblas_setup} we observe that using KA1S1CS1 tends to have a similar cost to that of VC1S1CS2, while VA3S1CS1 shows a small increase
to the setup time due to the application of $1$ step of the weighted Jacobi smoother
to the hierarchy prolongators, as described in Sec.~\ref{sec:parmatch}. VA5S1C1, based on a decoupled smoothed aggregation, has a clear advantage in the setup cost due to the absence of communication in the aggregation algorithm, and, as expected, obtains the best speedup for the setup phase. Very good speedups are also shown by KA1S1CS1 and VA3S1CS1, which confirms the effectiveness of the parallel implementations of all the computational kernels described in sections~\ref{sec:parmatch}-\ref{sec:smoothers}.
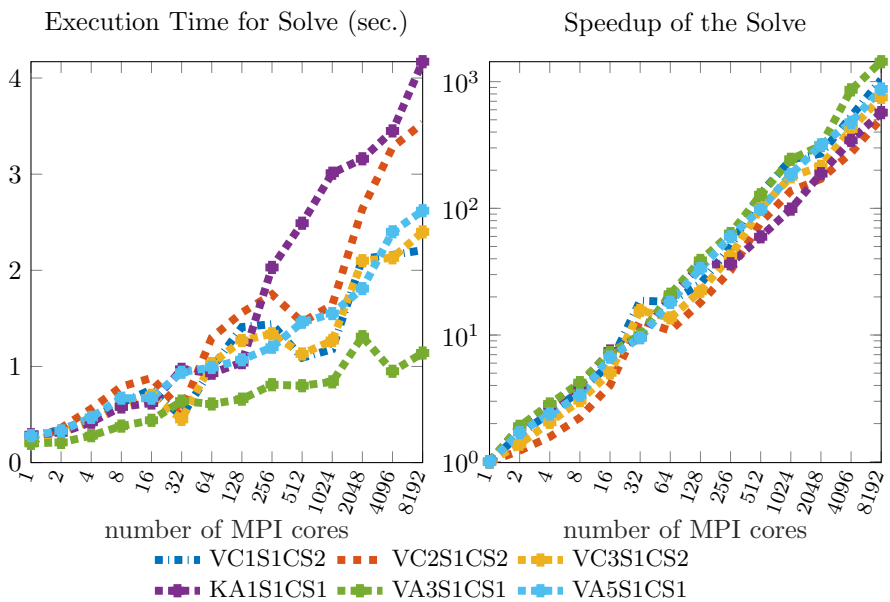
\begin{figure}[htbp]
    \centering
%
%
\definecolor{mycolor1}{rgb}{0.00000,0.44700,0.74100}%
\definecolor{mycolor2}{rgb}{0.85000,0.32500,0.09800}%
\definecolor{mycolor3}{rgb}{0.92900,0.69400,0.12500}%
\definecolor{mycolor4}{rgb}{0.49400,0.18400,0.55600}%
\definecolor{mycolor5}{rgb}{0.46600,0.67400,0.18800}%
\definecolor{mycolor6}{rgb}{0.30100,0.74500,0.93300}%
\begin{tikzpicture}

\begin{axis}[%
width=0.40\textwidth,
height=2.1in,
at={(0in,0in)},
scale only axis,
xmode=log,
xmin=1,
xmax=8192,
xtick={   1,    2,    4,    8,   16,   32,   64,  128,  256,  512, 1024, 2048, 4096, 8192},
xticklabels={   1,    2,    4,    8,   16,   32,   64,  128,  256,  512, 1024, 2048, 4096, 8192},
x tick label style={rotate=70,anchor=east,font=\footnotesize},
log basis x=2,
xminorticks=false,
xlabel style={font=\color{white!15!black},at={(0.5,-0.12)}},
xlabel={number of MPI cores},
ymin=0,
ymax=4.17,
yminorticks=true,
axis background/.style={fill=white},
title={Execution Time for Solve (sec.)},
legend columns=3, 
legend style={at={(1.0,-0.19)}, anchor=north, legend cell align=left, align=left,font=\small,draw=none}
]
\addplot [color=mycolor1, dashdotted, line width=3.0pt]
  table[row sep=crcr]{%
1	0.28\\
2	0.34\\
4	0.46\\
8	0.62\\
16	0.76\\
32	0.48\\
64	0.98\\
128	1.41\\
256	1.44\\
512	1.09\\
1024	1.18\\
2048	2.12\\
4096	2.16\\
8192	2.21\\
};
\addlegendentry{VC1S1CS2}

\addplot [color=mycolor2, dashed, line width=3.0pt]
  table[row sep=crcr]{%
1	0.22\\
2	0.36\\
4	0.56\\
8	0.79\\
16	0.88\\
32	0.54\\
64	1.3\\
128	1.56\\
256	1.75\\
512	1.46\\
1024	1.64\\
2048	2.64\\
4096	3.28\\
8192	3.53\\
};
\addlegendentry{VC2S1CS2}

\addplot [color=mycolor3, dotted, line width=3.0pt, mark size=2.5pt, mark=+, mark options={solid, mycolor3}]
  table[row sep=crcr]{%
1	0.22\\
2	0.32\\
4	0.43\\
8	0.58\\
16	0.7\\
32	0.45\\
64	1.03\\
128	1.27\\
256	1.34\\
512	1.13\\
1024	1.27\\
2048	2.1\\
4096	2.13\\
8192	2.4\\
};
\addlegendentry{VC3S1CS2}

\addplot [color=mycolor4, dotted, line width=3.0pt, mark size=2.5pt, mark=+, mark options={solid, mycolor4}]
  table[row sep=crcr]{%
1	0.29\\
2	0.31\\
4	0.42\\
8	0.58\\
16	0.62\\
32	0.97\\
64	0.93\\
128	1.04\\
256	2.03\\
512	2.49\\
1024	3.01\\
2048	3.16\\
4096	3.45\\
8192	4.17\\
};
\addlegendentry{KA1S1CS1}

\addplot [color=mycolor5, dotted, line width=3.0pt, mark size=2.5pt, mark=+, mark options={solid, mycolor5}]
  table[row sep=crcr]{%
1	0.2\\
2	0.21\\
4	0.28\\
8	0.38\\
16	0.44\\
32	0.64\\
64	0.61\\
128	0.66\\
256	0.81\\
512	0.8\\
1024	0.84\\
2048	1.31\\
4096	0.95\\
8192	1.14\\
};
\addlegendentry{VA3S1CS1}

\addplot [color=mycolor6, dotted, line width=3.0pt, mark size=2.5pt, mark=+, mark options={solid, mycolor6}]
  table[row sep=crcr]{%
1	0.28\\
2	0.33\\
4	0.47\\
8	0.67\\
16	0.67\\
32	0.94\\
64	0.99\\
128	1.07\\
256	1.2\\
512	1.46\\
1024	1.55\\
2048	1.81\\
4096	2.4\\
8192	2.62\\
};
\addlegendentry{VA5S1CS1}

\end{axis}

\begin{axis}[%
width=0.40\textwidth,
height=2.1in,
at={(2.4in,0in)},
scale only axis,
xmode=log,
xmin=1,
xmax=8192,
xtick={   1,    2,    4,    8,   16,   32,   64,  128,  256,  512, 1024, 2048, 4096, 8192},
xticklabels={   1,    2,    4,    8,   16,   32,   64,  128,  256,  512, 1024, 2048, 4096, 8192},
x tick label style={rotate=70,anchor=east,font=\footnotesize},
log basis x=2,
xminorticks=false,
xlabel style={font=\color{white!15!black},at={(0.5,-0.12)}},
xlabel={number of MPI cores},
ymode=log,
ymin=0.984126984126984,
ymax=1437.19298245614,
yminorticks=true,
axis background/.style={fill=white},
title={Speedup of the Solve},
]
\addplot [color=mycolor1, dashdotted, line width=3.0pt]
  table[row sep=crcr]{%
1	1\\
2	1.64705882352941\\
4	2.43478260869565\\
8	3.61290322580645\\
16	5.89473684210526\\
32	18.6666666666667\\
64	18.2857142857143\\
128	25.4184397163121\\
256	49.7777777777778\\
512	131.522935779817\\
1024	242.983050847458\\
2048	270.490566037736\\
4096	530.962962962963\\
8192	1037.90045248869\\
};

\addplot [color=mycolor2, dashed, line width=3.0pt]
  table[row sep=crcr]{%
1	1\\
2	1.22222222222222\\
4	1.57142857142857\\
8	2.22784810126582\\
16	4\\
32	13.037037037037\\
64	10.8307692307692\\
128	18.0512820512821\\
256	32.1828571428571\\
512	77.1506849315068\\
1024	137.365853658537\\
2048	170.666666666667\\
4096	274.731707317073\\
8192	510.549575070822\\
};

\addplot [color=mycolor3, dotted, line width=3.0pt, mark size=2.5pt, mark=+, mark options={solid, mycolor3}]
  table[row sep=crcr]{%
1	1\\
2	1.375\\
4	2.04651162790698\\
8	3.03448275862069\\
16	5.02857142857143\\
32	15.6444444444444\\
64	13.6699029126214\\
128	22.1732283464567\\
256	42.0298507462687\\
512	99.6814159292036\\
1024	177.385826771654\\
2048	214.552380952381\\
4096	423.06103286385\\
8192	750.933333333333\\
};

\addplot [color=mycolor4, dotted, line width=3.0pt, mark size=2.5pt, mark=+, mark options={solid, mycolor4}]
  table[row sep=crcr]{%
1	1\\
2	1.87096774193548\\
4	2.76190476190476\\
8	4\\
16	7.48387096774194\\
32	9.56701030927835\\
64	19.9569892473118\\
128	35.6923076923077\\
256	36.5714285714286\\
512	59.6305220883534\\
1024	98.6578073089701\\
2048	187.949367088608\\
4096	344.301449275362\\
8192	569.707434052758\\
};

\addplot [color=mycolor5, dotted, line width=3.0pt, mark size=2.5pt, mark=+, mark options={solid, mycolor5}]
  table[row sep=crcr]{%
1	1\\
2	1.9047619047619\\
4	2.85714285714286\\
8	4.21052631578947\\
16	7.27272727272727\\
32	10\\
64	20.9836065573771\\
128	38.7878787878788\\
256	63.2098765432099\\
512	128\\
1024	243.809523809524\\
2048	312.671755725191\\
4096	862.315789473684\\
8192	1437.19298245614\\
};

\addplot [color=mycolor6, dotted, line width=3.0pt, mark size=2.5pt, mark=+, mark options={solid, mycolor6}]
  table[row sep=crcr]{%
1	1\\
2	1.6969696969697\\
4	2.38297872340426\\
8	3.34328358208955\\
16	6.68656716417911\\
32	9.53191489361702\\
64	18.1010101010101\\
128	33.4953271028037\\
256	59.7333333333333\\
512	98.1917808219178\\
1024	184.98064516129\\
2048	316.817679558011\\
4096	477.866666666667\\
8192	875.480916030534\\
};
\end{axis}
\end{tikzpicture}%
    \caption{Comparison with \texttt{Hypre} 256k dofs per core. Solve: execution time (left), speedup (right).}
    \label{fig:hyprevsamg4psblas_solve}
\end{figure}
If we look at the solve phase in Fig.\ref{fig:hyprevsamg4psblas_solve}, we can see that VA3S1CS1 has generally the best solve time with respect to the other preconditioners, with a very small increase with increasing number of cores. This efficiency  compensates the small increase in the setup time with respect to the \texttt{Hypre} preconditioners when applied to time-dependent or non-linear problems, where the same preconditioner can be reused  for multiple  external iterations. We observe that our KA1S1CS1, although showing a very good algorithmic scalability, has a rapid increase in the solve time for increasing number of cores, due to the K-cycle application. Indeed this cycle requires $2^{nl-2}+1$ coarsest level visit and solutions per iteration, leading to a worse ratio between computation and communication with respect to the simpler V-cycle. Therefore, when the number of levels grows with the problem size so as to keep the coarsest matrix reasonably small, the above approach is poorly scalable. The speedups of the solve phase for all preconditioners are broadly comparable.}

{We conclude this comparison by claiming that our new coarsening based on weighted matching coupled with a smoothed version of the transfer operators{, and with a robust iterative solver at the coarsest level} shows very good scalability results and outperforms the preconditioners included in \texttt{Hypre} both in terms of memory requirements and in terms of solve times. Algorithmic scalability for Poisson-type problems of VA3S1CS1 appears optimal and comparable with the standard geometric multigrid approach.}

\subsection{Performance results towards extreme scale}
\label{res:pizdaintcpu-gpu}

In this section we discuss scalability results obtained on Piz Daint, running tests with $512 \times 10^3$ dofs per core up to $27000$ cores, i.e., we reach an overall number of $\sim 1.4 \times 10^{10}$ dofs; in the same vein, we also analyze results obtained with  GPU accelerators (in the solve phase). In the latter case we run tests with $12 \times 512 \times 10^3 \sim 6.2 \times 10^6$ dofs per GPU, and up to $2048$ GPUs, i.e, we reach an overall number dofs of more than $1.2 \times 10^{10}$; in these experiments on GPUs we keep the same amount of memory per node, hence each GPU will handle the same number of dofs as 12 CPU cores. In the same vein, we stop the coarsening process for the setup of the multilevel hierarchy when the maximum size of the coarsest matrix is $12\times 200 \times np$, where $np$ is the number of GPUs, i.e. with the same size of coarsest matrix per node.

{We begin with results on pure MPI; for the sake of space, we limit our discussion to preconditioners using the HGS smoother. We compare the KA1S1CS1 and KA2S1CS1 preconditioners (KA1/A2-types), using the K-cycle in the application of an AMG hierarchy where unsmoothed prolongators are employed, with the VA3S1CS1 and VA4S1CS1 preconditioners (VA3/A4-types), which use the V-cycle coupled with the smoothed version of the prolongators. 
In Fig.~\ref{fig:pizdaintcpu_opcmplx} we show the operator complexity of the multilevel hierarchies corresponding to all preconditioners. As expected, the operator complexity of the VA3/A4-type preconditioners, with smoothed prolongators, is larger than that of the corresponding KA1/A2-type preconditioners, since the coarse matrices are denser than those built with the unsmoothed prolongators. Nevertheless, when aggregates of size $8$ are built, the operator complexity is about $1.9$, while when aggregates of size  $16$ are employed, the operator complexity is about $1.3$. This indicates that even for VA3/A4-type preconditioners, the memory requirements for the AMG hierarchies are less than double the memory needed for the system matrix. 
Moreover, despite the small operator complexity, the numerical scalability of all the preconditioners is very satisfactory, even optimal for some of them; in Fig.~\ref{fig:pizdaintcpu_it} we see that KA1S1CS1 requires a number of iterations ranging from $12$ to $17$, while KA2S1CS1 requires a number of iterations ranging from $17$ to $24$, showing a small increase in iterations despite a reduction in operator complexity, and preserving numerical scalability for increasing number of cores. 
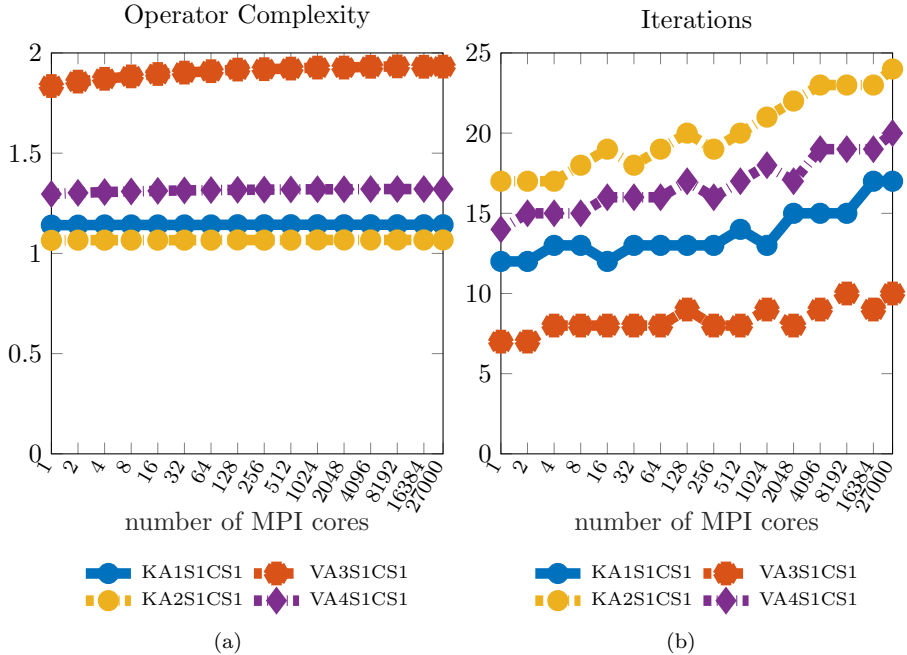
\begin{figure}[htbp]
    \centering
    \subfloat[\label{fig:pizdaintcpu_opcmplx}]{
%
%
\definecolor{mycolor1}{rgb}{0.00000,0.44700,0.74100}%
\definecolor{mycolor2}{rgb}{0.85000,0.32500,0.09800}%
\definecolor{mycolor3}{rgb}{0.92900,0.69400,0.12500}%
\definecolor{mycolor4}{rgb}{0.49400,0.18400,0.55600}%
\begin{tikzpicture}

\begin{axis}[%
width=0.40\textwidth,
height=2.1in,
at={(2.509in,2.689in)},
scale only axis,
xmode=log,
xmin=1,
xmax=27000,
xtick={    1,     2,     4,     8,    16,    32,    64,   128,   256,   512,  1024,  2048,  4096,  8192, 16384, 27000},
xticklabels={    1,     2,     4,     8,    16,    32,    64,   128,   256,   512,  1024,  2048,  4096,  8192, 16384, 27000},
x tick label style={rotate=60,anchor=east,font=\footnotesize},
xminorticks=true,
log basis x=2,
xminorticks=true,
xlabel style={font=\color{white!15!black},at={(0.5,-0.12)}},
xlabel={number of MPI cores},
ymin=0,
ymax=2,
yminorticks=true,
title={Operator Complexity},
axis background/.style={fill=white},
legend columns = 2,
legend style={at={(0.5,-0.25)}, anchor=north, legend cell align=left, align=left, draw=none,font=\scriptsize}
]
\addplot [color=mycolor1, line width=4.0pt, mark=o, mark options={solid, mycolor1}]
  table[row sep=crcr]{%
1	1.141\\
2	1.141\\
4	1.142\\
8	1.142\\
16	1.142\\
32	1.142\\
64	1.142\\
128	1.143\\
256	1.143\\
512	1.143\\
1024	1.143\\
2048	1.143\\
4096	1.143\\
8192	1.143\\
16384	1.143\\
27000	1.143\\
};
\addlegendentry{KA1S1CS1}

\addplot [color=mycolor2, dashed, line width=4.0pt,mark=asterisk, mark size=4.5pt, mark options={solid, mycolor2}]
  table[row sep=crcr]{%
1	1.835\\
2	1.857\\
4	1.871\\
8	1.883\\
16	1.895\\
32	1.902\\
64	1.909\\
128	1.917\\
256	1.919\\
512	1.922\\
1024	1.927\\
2048	1.927\\
4096	1.931\\
8192	1.933\\
16384	1.931\\
27000	1.932\\
};
\addlegendentry{VA3S1CS1}

\addplot [color=mycolor3, dashdotted, line width=4.0pt, mark=o, mark options={solid, mycolor3}]
  table[row sep=crcr]{%
1	1.065\\
2	1.066\\
4	1.066\\
8	1.066\\
16	1.066\\
32	1.066\\
64	1.066\\
128	1.067\\
256	1.066\\
512	1.066\\
1024	1.067\\
2048	1.067\\
4096	1.067\\
8192	1.067\\
16384	1.067\\
27000	1.067\\
};
\addlegendentry{KA2S1CS1}

\addplot [color=mycolor4, dashdotted, line width=4.0pt, mark=diamond, mark options={solid, mycolor4}]
  table[row sep=crcr]{%
1	1.296\\
2	1.302\\
4	1.306\\
8	1.309\\
16	1.312\\
32	1.314\\
64	1.316\\
128	1.317\\
256	1.318\\
512	1.319\\
1024	1.32\\
2048	1.32\\
4096	1.321\\
8192	1.321\\
16384	1.321\\
27000	1.321\\
};
\addlegendentry{VA4S1CS1}

\end{axis}
\end{tikzpicture}%
}
    \subfloat[\label{fig:pizdaintcpu_it}]{
%
%
\definecolor{mycolor1}{rgb}{0.00000,0.44700,0.74100}%
\definecolor{mycolor2}{rgb}{0.85000,0.32500,0.09800}%
\definecolor{mycolor3}{rgb}{0.92900,0.69400,0.12500}%
\definecolor{mycolor4}{rgb}{0.49400,0.18400,0.55600}%
\begin{tikzpicture}

\begin{axis}[%
width=0.40\textwidth,
height=2.1in,
at={(2.509in,2.689in)},
scale only axis,
xmode=log,
xmin=1,
xmax=27000,
xtick={    1,     2,     4,     8,    16,    32,    64,   128,   256,   512,  1024,  2048,  4096,  8192, 16384, 27000},
xticklabels={    1,     2,     4,     8,    16,    32,    64,   128,   256,   512,  1024,  2048,  4096,  8192, 16384, 27000},
x tick label style={rotate=60,anchor=east,font=\footnotesize},
xminorticks=true,
log basis x=2,
xminorticks=true,
xlabel style={font=\color{white!15!black},at={(0.5,-0.12)}},
xlabel={number of MPI cores},
ymin=0,
ymax=25,
yminorticks=true,
title={Iterations},
axis background/.style={fill=white},
legend columns = 2,
legend style={at={(0.5,-0.25)}, anchor=north, legend cell align=left, align=left, draw=none,font=\scriptsize}
]
\addplot [color=mycolor1, line width=4.0pt, mark=o, mark options={solid, mycolor1}]
  table[row sep=crcr]{%
1	12\\
2	12\\
4	13\\
8	13\\
16	12\\
32	13\\
64	13\\
128	13\\
256	13\\
512	14\\
1024	13\\
2048	15\\
4096	15\\
8192	15\\
16384	17\\
27000	17\\
};
\addlegendentry{KA1S1CS1}

\addplot [color=mycolor2, dashed, line width=4.0pt,mark=asterisk, mark size=4.5pt, mark options={solid, mycolor2}]
  table[row sep=crcr]{%
1	7\\
2	7\\
4	8\\
8	8\\
16	8\\
32	8\\
64	8\\
128	9\\
256	8\\
512	8\\
1024	9\\
2048	8\\
4096	9\\
8192	10\\
16384	9\\
27000	10\\
};
\addlegendentry{VA3S1CS1}

\addplot [color=mycolor3, dashdotted, line width=4.0pt, mark=o, mark options={solid, mycolor3}]
  table[row sep=crcr]{%
1	17\\
2	17\\
4	17\\
8	18\\
16	19\\
32	18\\
64	19\\
128	20\\
256	19\\
512	20\\
1024	21\\
2048	22\\
4096	23\\
8192	23\\
16384	23\\
27000	24\\
};
\addlegendentry{KA2S1CS1}

\addplot [color=mycolor4, dashdotted, line width=4.0pt, mark=diamond, mark options={solid, mycolor4}]
  table[row sep=crcr]{%
1	14\\
2	15\\
4	15\\
8	15\\
16	16\\
32	16\\
64	16\\
128	17\\
256	16\\
512	17\\
1024	18\\
2048	17\\
4096	19\\
8192	19\\
16384	19\\
27000	20\\
};
\addlegendentry{VA4S1CS1}

\end{axis}
\end{tikzpicture}%
}
    \caption{Weak Scaling results for 512k dofs per core. Operator complexity and number of iterations.}
\end{figure}
The VA3/A4-type preconditioners, employing the smoothed version of the prolongators and a less expensive cycle, require a smaller number of iterations than the corresponding KA1/A2-type preconditioners. In more details, VA3S1CS1 continues to show an almost perfect algorithmic scalability with
a number of iterations ranging from $7$ to $10$ going from $1$ core to $27700$. A very good algorithmic scalability is also observed for VA4S1CS1, where the number of iterations ranges from $14$
to $20$. This behavior confirms that the new parallel coarsening based on weighted matching is able
to detect, in a completely automatic way, good quality aggregates of variable size and, exploiting smoothed operators allows to obtain optimality also using V-cycle at low operator complexity.
VA3S1CS1 always obtains the best solve time per each number of cores, with a very slow increase for increasing number of cores; it solves the biggest size problem on $27000$ cores in less than $2$ seconds. In Fig.~\ref{fig:pizdaintcpu_solvetime} the application of the more expensive K-cycle shows its effect on the solve time, especially when increasing the number of cores. Better solve times than the KA1/A2-type preconditioners are also obtained by VA4S1CS1 for increasing number of cores, showing that the use of V-cycle coupled with A3/A4 aggregation types is the best choice for extreme scalability, albeit at a small increase in the setup cost (see Fig.~\ref{fig:pizdaintcpu_setuptime}). 
For the sake of completeness we report in Figs.~\ref{fig:pizdaintcpu_setupsp}-\ref{fig:pizdaintcpu_solvesp} the speedups obtained; we observe a very similar behavior for all methods, displaying a smooth increase with increasing number of cores, and demonstrating the good implementation scalability of all computational kernels.}
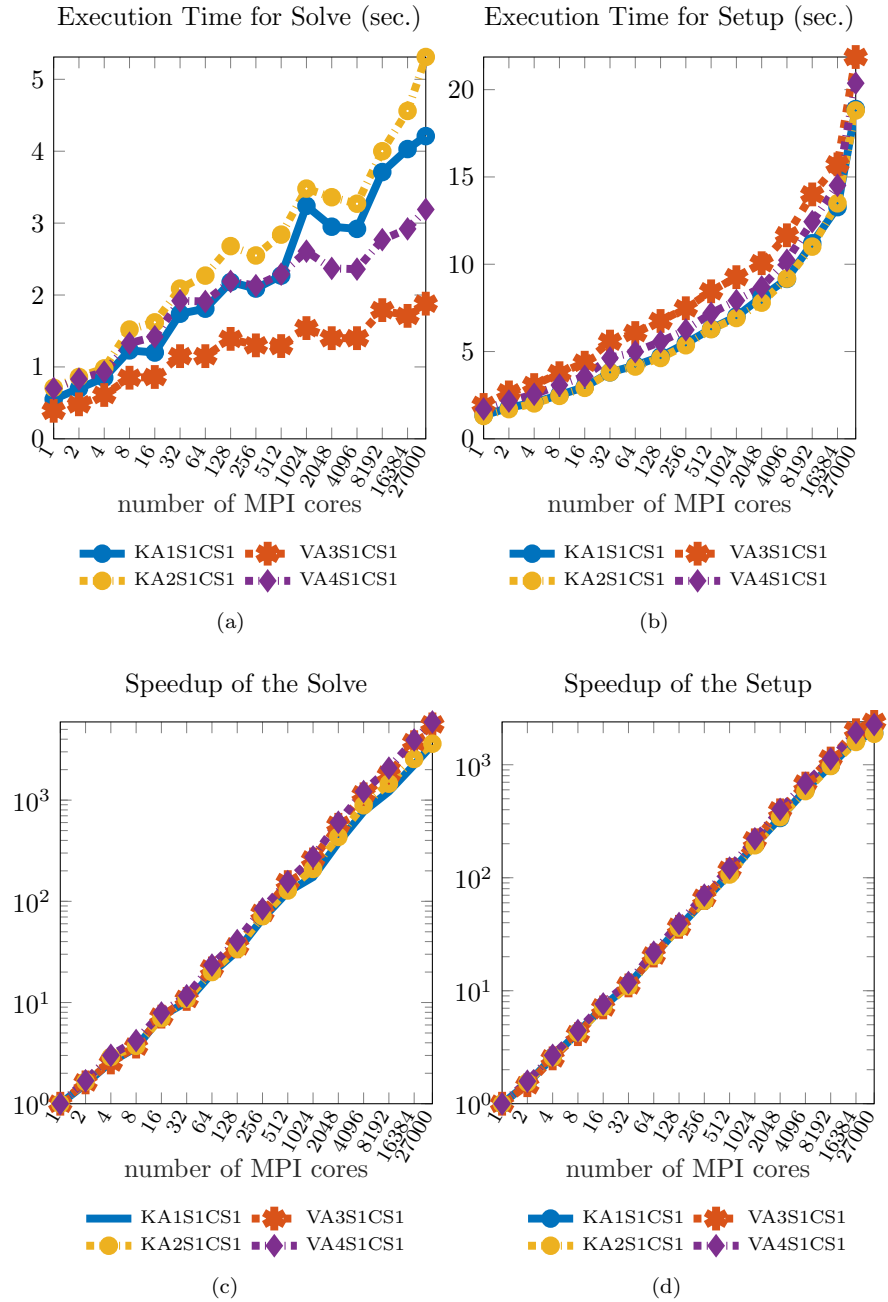
\begin{figure}[htbp]
    \centering
    \subfloat[\label{fig:pizdaintcpu_solvetime}]{
%
%
\definecolor{mycolor1}{rgb}{0.00000,0.44700,0.74100}%
\definecolor{mycolor2}{rgb}{0.85000,0.32500,0.09800}%
\definecolor{mycolor3}{rgb}{0.92900,0.69400,0.12500}%
\definecolor{mycolor4}{rgb}{0.49400,0.18400,0.55600}%
\begin{tikzpicture}

\begin{axis}[%
width=0.38\textwidth,
height=2in,
at={(2.509in,2.689in)},
scale only axis,
xmode=log,
xmin=1,
xmax=27000,
xtick={    1,     2,     4,     8,    16,    32,    64,   128,   256,   512,  1024,  2048,  4096,  8192, 16384, 27000},
xticklabels={    1,     2,     4,     8,    16,    32,    64,   128,   256,   512,  1024,  2048,  4096,  8192, 16384, 27000},
x tick label style={rotate=60,anchor=east,font=\footnotesize},
xminorticks=true,
log basis x=2,
xminorticks=true,
xlabel style={font=\color{white!15!black},at={(0.5,-0.12)}},
xlabel={number of MPI cores},
ymin=0,
ymax=5.31,
title={Execution Time for Solve (sec.)},
axis background/.style={fill=white},
legend columns = 2,
legend style={at={(0.5,-0.25)}, anchor=north, legend cell align=left, align=left, draw=none,font=\scriptsize}
]
\addplot [color=mycolor1, line width=3.0pt, mark=o, mark options={solid, mycolor1}]
  table[row sep=crcr]{%
1	0.55\\
2	0.7\\
4	0.85\\
8	1.23\\
16	1.2\\
32	1.74\\
64	1.81\\
128	2.18\\
256	2.09\\
512	2.27\\
1024	3.24\\
2048	2.95\\
4096	2.92\\
8192	3.71\\
16384	4.03\\
27000	4.21\\
};
\addlegendentry{KA1S1CS1}

\addplot [color=mycolor2, dashed, line width=3.0pt, mark=asterisk, mark size=4.5pt, mark options={solid, mycolor2}]
  table[row sep=crcr]{%
1	0.39\\
2	0.48\\
4	0.61\\
8	0.85\\
16	0.86\\
32	1.15\\
64	1.15\\
128	1.39\\
256	1.3\\
512	1.29\\
1024	1.54\\
2048	1.4\\
4096	1.4\\
8192	1.79\\
16384	1.71\\
27000	1.88\\
};
\addlegendentry{VA3S1CS1}

\addplot [color=mycolor3, dashdotted, line width=3.0pt, mark=o, mark options={solid, mycolor3}]
  table[row sep=crcr]{%
1	0.71\\
2	0.86\\
4	0.98\\
8	1.52\\
16	1.62\\
32	2.09\\
64	2.27\\
128	2.68\\
256	2.55\\
512	2.84\\
1024	3.48\\
2048	3.36\\
4096	3.27\\
8192	4\\
16384	4.56\\
27000	5.31\\
};
\addlegendentry{KA2S1CS1}

\addplot [color=mycolor4, dashdotted, line width=3.0pt, mark=diamond, mark options={solid, mycolor4}]
  table[row sep=crcr]{%
1	0.7\\
2	0.83\\
4	0.93\\
8	1.33\\
16	1.42\\
32	1.92\\
64	1.91\\
128	2.19\\
256	2.13\\
512	2.29\\
1024	2.61\\
2048	2.37\\
4096	2.36\\
8192	2.77\\
16384	2.92\\
27000	3.19\\
};
\addlegendentry{VA4S1CS1}

\end{axis}
\end{tikzpicture}%
}
    \subfloat[\label{fig:pizdaintcpu_setuptime}]{
%
%
\definecolor{mycolor1}{rgb}{0.00000,0.44700,0.74100}%
\definecolor{mycolor2}{rgb}{0.85000,0.32500,0.09800}%
\definecolor{mycolor3}{rgb}{0.92900,0.69400,0.12500}%
\definecolor{mycolor4}{rgb}{0.49400,0.18400,0.55600}%
\begin{tikzpicture}

\begin{axis}[%
width=0.38\textwidth,
height=2in,
at={(2.509in,2.689in)},
scale only axis,
xmode=log,
xmin=1,
xmax=27000,
xtick={    1,     2,     4,     8,    16,    32,    64,   128,   256,   512,  1024,  2048,  4096,  8192, 16384, 27000},
xticklabels={    1,     2,     4,     8,    16,    32,    64,   128,   256,   512,  1024,  2048,  4096,  8192, 16384, 27000},
x tick label style={rotate=60,anchor=east,font=\footnotesize},
xminorticks=true,
log basis x=2,
xminorticks=true,
xlabel style={font=\color{white!15!black},at={(0.5,-0.12)}},
xlabel={number of MPI cores},
ymin=0,
ymax=21.87,
yminorticks=true,
title={Execution Time for Setup (sec.)},
axis background/.style={fill=white},
legend columns = 2,
legend style={at={(0.5,-0.25)}, anchor=north, legend cell align=left, align=left, draw=none,font=\scriptsize}
]
\addplot [color=mycolor1, line width=3.0pt, mark=o, mark options={solid, mycolor1}]
  table[row sep=crcr]{%
1	1.34\\
2	1.78\\
4	2.08\\
8	2.51\\
16	2.97\\
32	3.81\\
64	4.2\\
128	4.71\\
256	5.47\\
512	6.3\\
1024	7.01\\
2048	8.12\\
4096	9.16\\
8192	11.2\\
16384	13.27\\
27000	18.89\\
};
\addlegendentry{KA1S1CS1}

\addplot [color=mycolor2, dashed, line width=3.0pt, mark=asterisk, mark size=4.5pt, mark options={solid, mycolor2}]
  table[row sep=crcr]{%
1	1.94\\
2	2.66\\
4	3.08\\
8	3.76\\
16	4.36\\
32	5.57\\
64	6.06\\
128	6.75\\
256	7.48\\
512	8.45\\
1024	9.26\\
2048	10.06\\
4096	11.66\\
8192	14\\
16384	15.72\\
27000	21.87\\
};
\addlegendentry{VA3S1CS1}

\addplot [color=mycolor3, dashdotted, line width=3.0pt, mark=o, mark options={solid, mycolor3}]
  table[row sep=crcr]{%
1	1.31\\
2	1.7\\
4	2.03\\
8	2.46\\
16	2.92\\
32	3.81\\
64	4.14\\
128	4.63\\
256	5.36\\
512	6.28\\
1024	6.92\\
2048	7.8\\
4096	9.18\\
8192	11.01\\
16384	13.49\\
27000	18.81\\
};
\addlegendentry{KA2S1CS1}

\addplot [color=mycolor4, dashdotted, line width=3.0pt, mark=diamond, mark options={solid, mycolor4}]
  table[row sep=crcr]{%
1	1.71\\
2	2.16\\
4	2.55\\
8	3.08\\
16	3.58\\
32	4.63\\
64	4.99\\
128	5.56\\
256	6.24\\
512	7.19\\
1024	7.92\\
2048	8.7\\
4096	10.23\\
8192	12.44\\
16384	14.52\\
27000	20.37\\
};
\addlegendentry{VA4S1CS1}

\end{axis}
\end{tikzpicture}%
}
    
    \subfloat[\label{fig:pizdaintcpu_solvesp}]{
%
%
\definecolor{mycolor1}{rgb}{0.00000,0.44700,0.74100}%
\definecolor{mycolor2}{rgb}{0.85000,0.32500,0.09800}%
\definecolor{mycolor3}{rgb}{0.92900,0.69400,0.12500}%
\definecolor{mycolor4}{rgb}{0.49400,0.18400,0.55600}%
\begin{tikzpicture}

\begin{axis}[%
width=0.38\textwidth,
height=2in,
at={(2.509in,2.689in)},
scale only axis,
xmode=log,
xmin=1,
xmax=27000,
xtick={    1,     2,     4,     8,    16,    32,    64,   128,   256,   512,  1024,  2048,  4096,  8192, 16384, 27000},
xticklabels={    1,     2,     4,     8,    16,    32,    64,   128,   256,   512,  1024,  2048,  4096,  8192, 16384, 27000},
x tick label style={rotate=60,anchor=east,font=\footnotesize},
xminorticks=true,
log basis x=2,
xminorticks=true,
xlabel style={font=\color{white!15!black},at={(0.5,-0.12)}},
xlabel={number of MPI cores},
ymode=log,
ymin=1,
ymax=5924.8,
yminorticks=true,
title={Speedup of the Solve},
axis background/.style={fill=white},
legend columns = 2,
legend style={at={(0.5,-0.25)}, anchor=north, legend cell align=left, align=left, draw=none,font=\scriptsize}
]
\addplot [color=mycolor1, line width=3.0pt]
  table[row sep=crcr]{%
1	1\\
2	1.57142857142857\\
4	2.58823529411765\\
8	3.57723577235772\\
16	7.33333333333333\\
32	10.1149425287356\\
64	19.4475138121547\\
128	32.2935779816514\\
256	67.3684210526316\\
512	124.052863436123\\
1024	173.827160493827\\
2048	381.830508474576\\
4096	771.506849315069\\
8192	1214.4474393531\\
16384	2236.02977667494\\
27000	3527.31591448931\\
};
\addlegendentry{KA1S1CS1}

\addplot [color=mycolor2, dashed, line width=3.0pt, mark=asterisk, mark size=4.5pt, mark options={solid, mycolor2}]
  table[row sep=crcr]{%
1	1\\
2	1.625\\
4	2.55737704918033\\
8	3.67058823529412\\
16	7.25581395348837\\
32	10.8521739130435\\
64	21.704347826087\\
128	35.9136690647482\\
256	76.8\\
512	154.790697674419\\
1024	259.324675324675\\
2048	570.514285714286\\
4096	1141.02857142857\\
8192	1784.84916201117\\
16384	3736.70175438597\\
27000	5601.06382978723\\
};
\addlegendentry{VA3S1CS1}

\addplot [color=mycolor3, dashdotted, line width=3.0pt, mark=o, mark options={solid, mycolor3}]
  table[row sep=crcr]{%
1	1\\
2	1.65116279069767\\
4	2.89795918367347\\
8	3.73684210526316\\
16	7.01234567901234\\
32	10.8708133971292\\
64	20.0176211453744\\
128	33.910447761194\\
256	71.278431372549\\
512	128\\
1024	208.919540229885\\
2048	432.761904761905\\
4096	889.345565749235\\
8192	1454.08\\
16384	2551.01754385965\\
27000	3610.16949152542\\
};
\addlegendentry{KA2S1CS1}

\addplot [color=mycolor4, dashdotted, line width=3.0pt, mark=diamond, mark options={solid, mycolor4}]
  table[row sep=crcr]{%
1	1\\
2	1.68674698795181\\
4	3.01075268817204\\
8	4.21052631578947\\
16	7.88732394366197\\
32	11.6666666666667\\
64	23.455497382199\\
128	40.9132420091324\\
256	84.131455399061\\
512	156.506550218341\\
1024	274.63601532567\\
2048	604.894514767932\\
4096	1214.91525423729\\
8192	2070.18050541516\\
16384	3927.67123287671\\
27000	5924.76489028213\\
};
\addlegendentry{VA4S1CS1}

\end{axis}
\end{tikzpicture}%
}
    \subfloat[\label{fig:pizdaintcpu_setupsp}]{
%
%
\definecolor{mycolor1}{rgb}{0.00000,0.44700,0.74100}%
\definecolor{mycolor2}{rgb}{0.85000,0.32500,0.09800}%
\definecolor{mycolor3}{rgb}{0.92900,0.69400,0.12500}%
\definecolor{mycolor4}{rgb}{0.49400,0.18400,0.55600}%
\begin{tikzpicture}

\begin{axis}[%
width=0.38\textwidth,
height=2in,
at={(2.509in,2.689in)},
scale only axis,
xmode=log,
xmin=1,
xmax=27000,
xtick={    1,     2,     4,     8,    16,    32,    64,   128,   256,   512,  1024,  2048,  4096,  8192, 16384, 27000},
xticklabels={    1,     2,     4,     8,    16,    32,    64,   128,   256,   512,  1024,  2048,  4096,  8192, 16384, 27000},
x tick label style={rotate=60,anchor=east,font=\footnotesize},
xminorticks=true,
log basis x=2,
xminorticks=true,
xlabel style={font=\color{white!15!black},at={(0.5,-0.12)}},
xlabel={number of MPI cores},
ymode=log,
ymin=1,
ymax=2395.06172839506,
yminorticks=true,
title={Speedup of the Setup},
axis background/.style={fill=white},
legend columns = 2,
legend style={at={(0.5,-0.25)}, anchor=north, legend cell align=left, align=left, draw=none,font=\scriptsize}
]
\addplot [color=mycolor1, line width=3.0pt, mark=o, mark options={solid, mycolor1}]
  table[row sep=crcr]{%
1	1\\
2	1.50561797752809\\
4	2.57692307692308\\
8	4.27091633466136\\
16	7.21885521885522\\
32	11.254593175853\\
64	20.4190476190476\\
128	36.416135881104\\
256	62.7129798903108\\
512	108.901587301587\\
1024	195.743223965763\\
2048	337.970443349754\\
4096	599.196506550218\\
8192	980.114285714286\\
16384	1654.45064054258\\
27000	1915.29910005294\\
};
\addlegendentry{KA1S1CS1}

\addplot [color=mycolor2, dashed, line width=3.0pt, mark=asterisk, mark size=4.5pt, mark options={solid, mycolor2}]
  table[row sep=crcr]{%
1	1\\
2	1.45864661654135\\
4	2.51948051948052\\
8	4.12765957446809\\
16	7.11926605504587\\
32	11.1454219030521\\
64	20.4884488448845\\
128	36.7881481481481\\
256	66.3957219251337\\
512	117.547928994083\\
1024	214.5313174946\\
2048	394.942345924453\\
4096	681.495711835334\\
8192	1135.17714285714\\
16384	2021.94402035623\\
27000	2395.06172839506\\
};
\addlegendentry{VA3S1CS1}

\addplot [color=mycolor3, dashdotted, line width=3.0pt, mark=o, mark options={solid, mycolor3}]
  table[row sep=crcr]{%
1	1\\
2	1.54117647058824\\
4	2.58128078817734\\
8	4.26016260162602\\
16	7.17808219178082\\
32	11.002624671916\\
64	20.2512077294686\\
128	36.2159827213823\\
256	62.5671641791045\\
512	106.802547770701\\
1024	193.849710982659\\
2048	343.958974358974\\
4096	584.505446623094\\
8192	974.706630336058\\
16384	1591.03335804299\\
27000	1880.38277511962\\
};
\addlegendentry{KA2S1CS1}

\addplot [color=mycolor4, dashdotted, line width=3.0pt, mark=diamond, mark options={solid, mycolor4}]
  table[row sep=crcr]{%
1	1\\
2	1.58333333333333\\
4	2.68235294117647\\
8	4.44155844155844\\
16	7.64245810055866\\
32	11.8185745140389\\
64	21.9318637274549\\
128	39.3669064748201\\
256	70.1538461538461\\
512	121.769123783032\\
1024	221.090909090909\\
2048	402.537931034483\\
4096	684.66862170088\\
8192	1126.07073954984\\
16384	1929.52066115702\\
27000	2266.56848306333\\
};
\addlegendentry{VA4S1CS1}

\end{axis}
\end{tikzpicture}%
}
    \caption{Weak scaling results 512k dofs per core. Setup and solve time with the relative speed-up.}\label{fig:pizdaintcpu}
\end{figure}

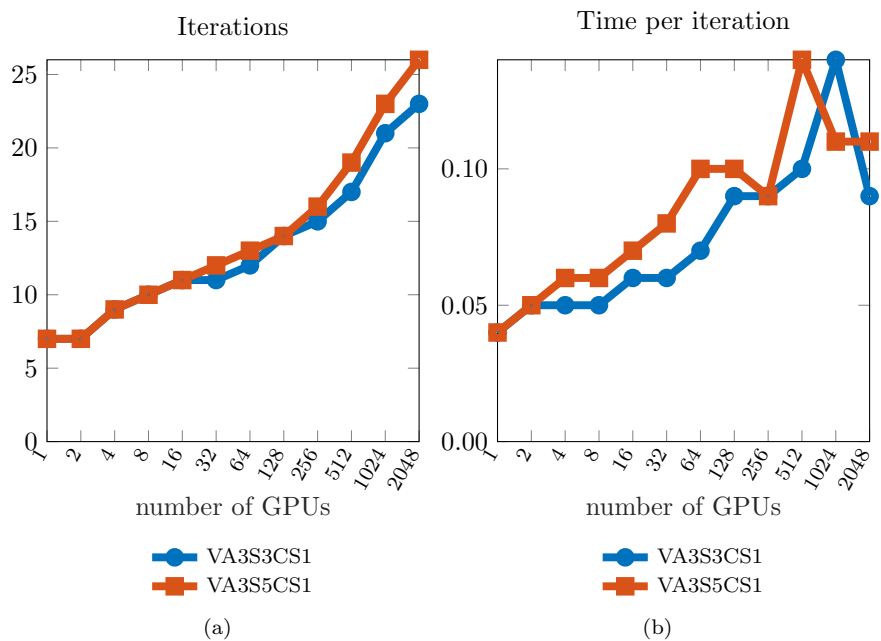
\begin{figure}[htbp]
    \centering
    \subfloat[\label{fig:vspmc3it_pizdaint_gpu}]{
%
%
\definecolor{mycolor1}{rgb}{0.00000,0.44700,0.74100}%
\definecolor{mycolor2}{rgb}{0.85000,0.32500,0.09800}%
\begin{tikzpicture}

\begin{axis}[%
width=0.38\textwidth,
height=2in,
at={(2.509in,2.133in)},
scale only axis,
xmode=log,
xmin=1,
xmax=2048,
xtick={   1,    2,    4,    8,   16,   32,   64,  128,  256,  512, 1024, 2048},
xticklabels={    1,     2,     4,     8,    16,    32,    64,   128,   256,   512,  1024,  2048,  4096,  8192, 16384},
x tick label style={rotate=60,anchor=east,font=\footnotesize},
xminorticks=true,
log basis x=2,
xminorticks=true,
xlabel style={font=\color{white!15!black},at={(0.5,-0.12)}},
xlabel={number of GPUs},
ymin=0,
ymax=26,
yminorticks=true,
title={Iterations},
axis background/.style={fill=white},
legend columns = 1,
legend style={at={(0.5,-0.25)}, anchor=north, legend cell align=left, align=left, draw=none,font=\footnotesize}
]
\addplot [color=mycolor1, line width=3.0pt, mark=o, mark options={solid, mycolor1}]
  table[row sep=crcr]{%
1	7\\
2	7\\
4	9\\
8	10\\
16	11\\
32	11\\
64	12\\
128	14\\
256	15\\
512	17\\
1024	21\\
2048	23\\
};
\addlegendentry{VA3S3CS1}

\addplot [color=mycolor2, line width=3.0pt, mark=square, mark options={solid, mycolor2}]
  table[row sep=crcr]{%
1	7\\
2	7\\
4	9\\
8	10\\
16	11\\
32	12\\
64	13\\
128	14\\
256	16\\
512	19\\
1024	23\\
2048	26\\
};
\addlegendentry{VA3S5CS1}

\end{axis}
\end{tikzpicture}%
}
    \subfloat[\label{fig:vspmc3timexiter_pizdaint_gpu}]{
%
%
\definecolor{mycolor1}{rgb}{0.00000,0.44700,0.74100}%
\definecolor{mycolor2}{rgb}{0.85000,0.32500,0.09800}%
\begin{tikzpicture}

\begin{axis}[%
width=0.38\textwidth,
height=2in,
at={(2.509in,2.133in)},
scale only axis,
xmode=log,
xmin=1,
xmax=2048,
xtick={   1,    2,    4,    8,   16,   32,   64,  128,  256,  512, 1024, 2048},
xticklabels={    1,     2,     4,     8,    16,    32,    64,   128,   256,   512,  1024,  2048,  4096,  8192, 16384},
x tick label style={rotate=60,anchor=east,font=\footnotesize},
xminorticks=true,
log basis x=2,
xminorticks=true,
xlabel style={font=\color{white!15!black},at={(0.5,-0.12)}},
xlabel={number of GPUs},
ymin=0,
ymax=0.14,
yminorticks=true,
y tick label style={
    /pgf/number format/.cd,
        fixed,
        fixed zerofill,
        precision=2,
    /tikz/.cd
},
title={Time per iteration},
axis background/.style={fill=white},
legend columns = 1,
legend style={at={(0.5,-0.25)}, anchor=north, legend cell align=left, align=left, draw=none,font=\footnotesize}
]
\addplot [color=mycolor1, line width=3.0pt, mark=o, mark options={solid, mycolor1}]
  table[row sep=crcr]{%
1	0.04\\
2	0.05\\
4	0.05\\
8	0.05\\
16	0.06\\
32	0.06\\
64	0.07\\
128	0.09\\
256	0.09\\
512	0.1\\
1024	0.14\\
2048	0.09\\
};
\addlegendentry{VA3S3CS1}

\addplot [color=mycolor2, line width=3.0pt, mark=square, mark options={solid, mycolor2}]
  table[row sep=crcr]{%
1	0.04\\
2	0.05\\
4	0.06\\
8	0.06\\
16	0.07\\
32	0.08\\
64	0.1\\
128	0.1\\
256	0.09\\
512	0.14\\
1024	0.11\\
2048	0.11\\
};
\addlegendentry{VA3S5CS1}

\end{axis}
\end{tikzpicture}%
}
    \caption{Weak scaling results 6M dofs per GPU. Number of iterations and time per iteration on GPUs.}
\end{figure}

The GPU plugin of \texttt{PSBLAS} allows us to run the solve phase on a cluster of GPUs, implementing the preconditioned FCG method coupled with some of the preconditioners included in \texttt{AMG4PSBLAS}. The plugin implements efficient  GPU versions of  the sparse matrix-vector products and vector-vector operations, such as vector updates and scalar products, including the necessary MPI communications~\cite{CaFiRo:14,Filippone:2017}.
In the following we discuss results obtained by employing VA3-type preconditioners coupled with both the HINVK and $\ell_1-$Jacobi smoothers, whose implementation on GPU users  the sparse matrix-vector kernel. We compare VA3S3CS1 and VA3S5CS1, when they are coupled with the FCG iterative solver; in the case of VA3S3CS1, we apply HINVK also as preconditioner of the CG method at the coarsest level, while for VA3S5CS1 we use $\ell_1-$Jacobi both as smoother and as preconditioner for CG at the coarsest level. Note that, due to the larger constant value of the $\ell_1-$Jacobi smoother (see Tab.~\ref{tab:smoothingfactor}), we apply $4$ pre/post-smoothing sweeps of the method at each V-cycle application, whereas for HINVK we only apply 1 sweep.
From Fig.~\ref{fig:vspmc3it_pizdaint_gpu} we see that both methods have a very similar behavior, showing a number of iterations ranging from $7$ to $23$ for VA3S3CS1 and from $7$ to $26$ for VA3S5CS1. The cost per iteration of VA3S3CS1 is generally better, as shown in Fig.~\ref{fig:vspmc3timexiter_pizdaint_gpu}: we can observe a cost per iteration ranging from $0.04$ to $0.14$ seconds, corresponding to a solve time per dof in the range $[10^{-12}:10^{-9}]$.  
\begin{figure}[h]
    \centering
    \subfloat[\label{fig:vspmc3solvetime_pizdaint_gpu}]{
%
%
\definecolor{mycolor1}{rgb}{0.00000,0.44700,0.74100}%
\definecolor{mycolor2}{rgb}{0.85000,0.32500,0.09800}%
\begin{tikzpicture}

\begin{axis}[%
width=0.38\textwidth,
height=2in,
at={(2.509in,2.133in)},
scale only axis,
xmode=log,
xmin=1,
xmax=2048,
xtick={   1,    2,    4,    8,   16,   32,   64,  128,  256,  512, 1024, 2048},
xticklabels={    1,     2,     4,     8,    16,    32,    64,   128,   256,   512,  1024,  2048,  4096,  8192, 16384},
x tick label style={rotate=60,anchor=east,font=\footnotesize},
xminorticks=true,
log basis x=2,
xminorticks=true,
xlabel style={font=\color{white!15!black},at={(0.5,-0.12)}},
xlabel={number of GPUs},
ymin=0,
ymax=2.87,
title={Execution Time for Solve (sec.)},
axis background/.style={fill=white},
legend columns = 1,
legend style={at={(0.5,-0.25)}, anchor=north, legend cell align=left, align=left, draw=none,font=\footnotesize}
]
\addplot [color=mycolor1, line width=3.0pt, mark=o, mark options={solid, mycolor1}]
  table[row sep=crcr]{%
1	0.26\\
2	0.36\\
4	0.46\\
8	0.54\\
16	0.66\\
32	0.7\\
64	0.87\\
128	1.22\\
256	1.3\\
512	1.63\\
1024	2.87\\
2048	2.1\\
};
\addlegendentry{VA3S3CS1}

\addplot [color=mycolor2, line width=3.0pt, mark=square, mark options={solid, mycolor2}]
  table[row sep=crcr]{%
1	0.28\\
2	0.38\\
4	0.52\\
8	0.65\\
16	0.79\\
32	0.92\\
64	1.29\\
128	1.37\\
256	1.48\\
512	2.73\\
1024	2.62\\
2048	2.85\\
};
\addlegendentry{VA3S5CS1}

\end{axis}
\end{tikzpicture}%
}
    \subfloat[\label{fig:vspmc3setuptime_pizdaint_gpu}]{
%
%
\definecolor{mycolor1}{rgb}{0.00000,0.44700,0.74100}%
\definecolor{mycolor2}{rgb}{0.85000,0.32500,0.09800}%
\begin{tikzpicture}

\begin{axis}[%
width=0.38\textwidth,
height=2in,
at={(2.509in,2.133in)},
scale only axis,
xmode=log,
xmin=1,
xmax=2048,
xtick={   1,    2,    4,    8,   16,   32,   64,  128,  256,  512, 1024, 2048},
xticklabels={    1,     2,     4,     8,    16,    32,    64,   128,   256,   512,  1024,  2048,  4096,  8192, 16384},
x tick label style={rotate=60,anchor=east,font=\footnotesize},
xminorticks=true,
log basis x=2,
xminorticks=true,
xlabel style={font=\color{white!15!black},at={(0.5,-0.12)}},
xlabel={number of CPUs},
ymin=0,
ymax=163.3,
yminorticks=true,
title={Execution Time for Setup (sec.)},
axis background/.style={fill=white},
legend columns = 1,
legend style={at={(0.5,-0.25)}, anchor=north, legend cell align=left, align=left, draw=none,font=\footnotesize}
]
\addplot [color=mycolor1, line width=3.0pt, mark=o, mark options={solid, mycolor1}]
  table[row sep=crcr]{%
1	113.07\\
2	116\\
4	118.31\\
8	116.94\\
16	121.91\\
32	129.67\\
64	129.91\\
128	134.46\\
256	143.21\\
512	149.31\\
1024	150.1\\
2048	163.3\\
};
\addlegendentry{VA3S3CS1}

\addplot [color=mycolor2, line width=3.0pt, mark=square, mark options={solid, mycolor2}]
  table[row sep=crcr]{%
1	27.64\\
2	30.85\\
4	33.4\\
8	37.15\\
16	39.45\\
32	42.62\\
64	45.24\\
128	51.51\\
256	54.07\\
512	63.47\\
1024	66.48\\
2048	75.31\\
};
\addlegendentry{VA3S5CS1}

\end{axis}
\end{tikzpicture}%
}
    \caption{Weak scaling results 6M dofs per GPU. Solve time on GPUs and Setup time on CPU.}\label{fig:pizdaint_gpu}
\end{figure}
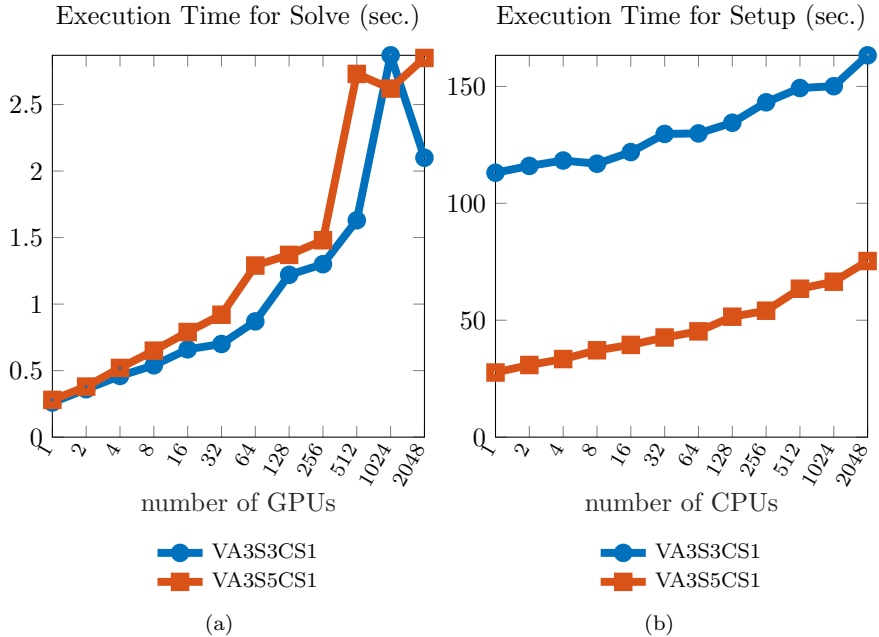
For the total solve time reported in~Fig.\ref{fig:vspmc3solvetime_pizdaint_gpu}, we observe that the two preconditioners have a similar behavior but VA3S3CS1 is generally better. In Fig.~\ref{fig:vspmc3setuptime_pizdaint_gpu} we also show the setup time for both preconditioners. In the current version of the library the setup of the preconditioners is not yet implemented on the GPU: indeed, it is carried out on a single core of the CPU host device. We see that, as expected, VA3S3CS1 shows larger setup costs with respect to VA3S5CS1 due to the larger cost for the setup of the HINVK smoother. The smoother setup  cost, as expected, is about constant for all numbers of CPU cores, with the $\ell_1-$Jacobi smoother being less expensive by about 2 orders of magnitude. Further development activities planned for future releases of the library include the implementation of hybrid OpenMP-MPI and possibly CUDA versions of the HINVK setup phase to make better use of hybrid computing nodes. Let us observe that in time-dependent problems, such as in most of the test cases from the EoCoE-II project, large setup times are generally well-tolerated by the computational procedure if the solve phase is very efficient, since the same preconditioner is applied in a very large number of time steps.
An extended comparison of the behavior of  the various preconditioners in Figures~\ref{fig:pizdaintcpu}, and~\ref{fig:pizdaint_gpu}, is available in the supplementary materials; in particular we have illustrated in Fig.~\ref{fig:global_behavior} the solve time for CPU and GPU variants for weak scaling up to $10^{10}$ dofs.

\section{Conclusions and future work}
\label{sec:end}
We presented a new version of a parallel preconditioners package implementing AMG methods, addressing scalability, flexibility and robustness for high-performance scientific computing at extreme scale. 
We focused on the design and implementation of a new parallel coarsening algorithm to be used in conjunction with highly parallel smoothers within AMG cycles for iterative solution of s.p.d. linear systems with size larger than billions on a current pre-exascale computer. Scalability results both on CPU and on GPU accelerators are very promising and a comparison with available software demonstrates the validity of our approaches both in terms of algorithms and in terms of software development.  

{Currently we are  investigating  alternative OpenMP and CUDA implementations of the kernels involved in the AMG setup, including graph matching, sparse matrix by sparse matrix product and approximate sparse matrix inversion, as well as other kernels employed in various phases of the matrix and communication  setup. We plan to include these alternative OpenMP and/or CUDA versions in future releases of the package.}

\section*{Acknowledgments}
We gratefully acknowledge  PRACE for awarding us access to 
Piz Daint at ETH Zurich/CSCS, Switzerland. We thank Mahantesh Halappanavar from PPNL-USA for useful discussion on the algorithm and implementation details of \texttt{MatchBox-P}. Finally we thank the anonymous reviewers for their comments and remarks which helped us in  improving the paper.

\bibliographystyle{siamplain}
\bibliography{bibliography}

\pagebreak
\setcounter{page}{0}
\pagenumbering{arabic}
\setcounter{page}{1}
\renewcommand{\thepage}{SM\arabic{page}}
\appendix
\renewcommand{\thesection}{\arabic{section}}
\setcounter{section}{0}

\makeatletter
\null
\@slugger
\begin{center}%
	\@supertitle%
	{\bfseries\MakeUppercase{\siampretitle}\@ucnothanks{SUPPLEMENTARY MATERIALS: AMG preconditioners for Linear Solvers towards Extreme Scale}}
	\vskip .075in%
	{\footnotesize
		\spaceskip=1.5\fontdimen2\the\and{\ \scriptsize AND }
		\@tempswafalse\expandafter\@doauth\@author\and\relax}%
	\@subtitle%
\end{center}%
\par
\vskip .11in
\def\@ucnothanks#1{\expandafter\@ucnt#1\thanks\relax}
\makeatother

This supplementary document provides further information on the experiments described in the paper. The sections to which the information refer have the same name of the corresponding sections in the paper.

\section{Parallel Smoothers and Coarsest Solvers} To describe the effect of using different smoothers in the AMG algorithm, we have used the two-level convergence analysis based on the computation of the constant $K$ to bound the spectral radius of the error propagation matrix of the method, i.e: 
\begin{equation*}
\|I-BA\|_A^2 = \| (I - P(P^TAP)^{-1}P^TA)(I-M^{-1}A) \|_A^2 \leq 1 - \frac{1}{K},
\end{equation*}
for which it holds
\begin{equation*}
K = \sup_{e} \frac{\|(I-PR) e \|_{\overline{M}}^2}{\|e\|_A^2} \geq 1, \quad \overline{M} = M^T(M^T + M - A)^{-1} M,
\end{equation*}
for a fixed projector $P$, and $R$ any matrix such that $PR$ is a projection onto $\Range(P)$. As we mentioned in section~\ref{sec:smoothers}, to improve the convergence properties of a single V-cycle, one can consider a smoothed version of such a prolongator, i.e., a more accurate interpolation operator obtained by applying one step of a weighted-Jacobi smoother to the $P$ obtained from the matching strategy. The actual prolongator $\overline{P}$ is obtained from $P$ as $\overline{P}=(I-\omega D^{-1}A)P$, where $D=\diag(A)$ and $\omega = 1/\|D^{-1}A\|_{\infty} \approx 1/\rho(D^{-1}A)$, with $\rho(D^{-1}A)$ the spectral radius of $D^{-1}A$. In this case we no longer have an explicit expression for the $P$ and $R$ needed for computing the constant $K$ inherited from the matching construction in section~\ref{sec:parmatch}. Nevertheless, we can analyze the behavior on the same test case by numerically computing an $R$ such that $PR$ is a projection onto $\Range(P)$. 
\begin{table}[htbp]
	\centering
	\caption{Convergence constant $K$ for various smoothers for the homogeneous 3D Laplacian problem with parallel coarsening based on compatible weighted matching. The global size of the problem is $m = 2^4 \times 2^4 \times 2^4$ distributed over $np$ processes using a 3D block distribution. By INVK we denote here the approximate inverse computed from the incomplete LU factorization with $0$ levels of fill-in, and admitting a single level of fill-in in the inversion procedure.}
	\begin{tabular}{rrcccccccc}
		&   & \multicolumn{6}{c}{Smoothed  prolongator 3 sweeps of matching} \\
		\cmidrule{3-8}
		m    & np & HGS  & L1-HGS & HINVK & L1-INVK & HINVK & L1-INVK\\
		&   &      &        &       \multicolumn{2}{c}{$\omega = 1$}         & \multicolumn{2}{c}{$\omega_{\text{opt}}$}\\
		\midrule
		4096 & 1 & 1.3686 & 1.3686  & 1.5611 & 1.5611 & 1.0056 & 1.0056\\
		2048 & 2 & 1.4108 & 1.5166 & 1.5293 & 1.5150 & 1.0390 & 1.0751 \\
		1024 & 4 & 1.4496 & 1.6500 & 1.5681 & 1.6318 & 1.2541 & 1.3553 \\
		512  & 8 & 1.4632 & 1.7443 & 1.7362 & 1.7450 & 1.4971 & 1.5554 \\
		256 & 16 & 1.4832 & 1.7938 & 1.7796 & 1.7966 & 1.5310 & 1.5400 \\
		128 & 32 & 1.5035 & 1.8435 & 1.8172 & 1.8413 & 1.5663  & 1.5371\\
		64  & 64 & 1.5218 & 1.8948 & 1.8533 & 1.8836 & 1.6034 & 1.5523 \\
		32 & 128 & 1.5757 & 2.0123 & 1.9392  & 2.0133 & 1.7287 & 1.6242\\
		16 & 256 & 1.6278 & 2.1369 & 2.0326 & 2.1491 & 1.8903 & 1.7451 \\
		8 & 512  & 1.6523 & 2.2833 & 2.1858 & 2.2921 & 2.1073 & 1.8999 \\
		\bottomrule
	\end{tabular}
	\label{tab:convergence_smoothed_aggregation}
\end{table}
In Tab.~\ref{tab:convergence_smoothed_aggregation} we report the results of these  computations; if we compare them  with those  obtained for the unsmoothed aggregation in Tab.~\ref{tab:smoothingfactor}, we observe the same behavior with respect to the smoothers, but we achieve better constants in absolute terms. This confirms  that increasing the regularity of the projection operator leads to better convergence properties  in this framework.

\section{Comparison with \texttt{Hypre}}\label{sec:hyprecoarsening}
We compare our preconditioners with some preconditioners available in the \texttt{Hypre} library~\cite{FJY2006}, using default algorithmic parameters for best practices. For our choice of  coarsening algorithms we consider:
\begin{description}
	\item[Falgout] This is a combination of the classic Ruge-St\"uben and the parallel independent set, named CLJP, coarsening algorithms. An initial RS0 coarsening requiring no communication is built; this generates a set of coarse variables (C-indices) for each processor, and the the variables associated with  the internal points of each process are subsequently used as the first independent set in a CLJP coarsening sweep (see~\cite{MR1908755}, and~\cite[section~3.2]{MR2267933}).
	\item[HMIS] The Hybrid Maximal Independent Set algorithm uses one pass of the Ruge-St\"uben coarsening algorithm on each processor independently. The interior coarse variables (C-indices) are then used as first independent sets for a parallel coarsening algorithm~\cite{MR2205610}. 
	\item[HMIS1] This is a variant of the HMIS algorithm in which $1$ level of \emph{aggressive} coarsening is used, i.e., instead of the classical Ruge-St\"uben coarsening algorithm a stricter definition of connectivity is enforced resulting in sparser coarse grids; see~\cite[section~3]{MR2650221} for the details.
\end{description} 
{We also include here the comparison between the \texttt{Hypre} preconditioners included in the paper and the same versions of the preconditioners, when the direct solution at the coarsest level is substituted by an iterative solution. \texttt{Hypre} makes available at the coarsest level the same (stationary) smoothers available at the intermediate levels; we selected the Symmetrized Gauss-Seidel method, which is the most robust among the available methods, and applied 10 iterations of this method to solve the coarsest system. All the other parameters are left as in the default.
	We can observe that there is no significant difference between the use of the iterative coarsest solver (see the preconditioner's names with the CS3 suffix) and the direct one (whose preconditioner's names have the CS2 suffix) both in terms of iterations and in terms of time for the setup and the solve phases.}
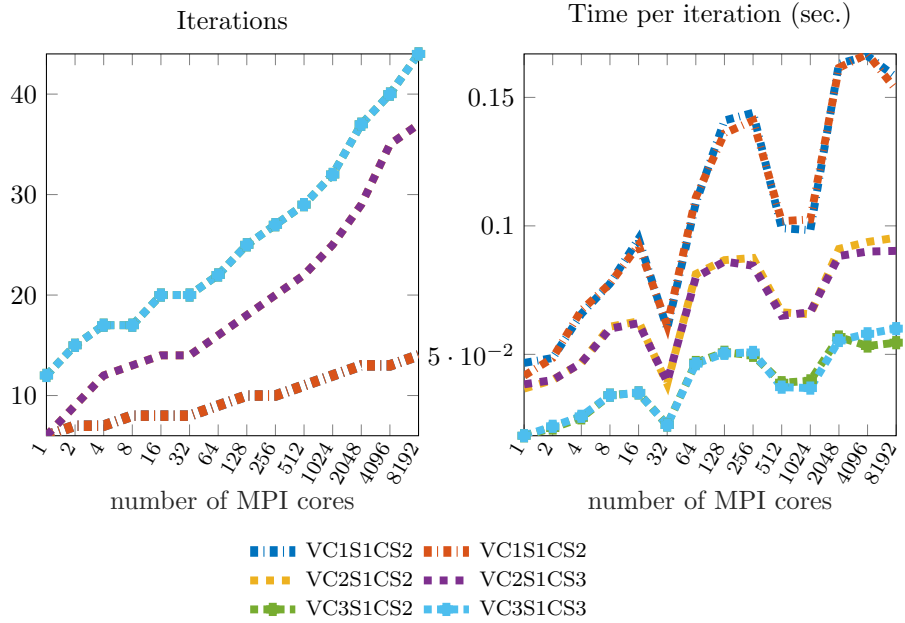
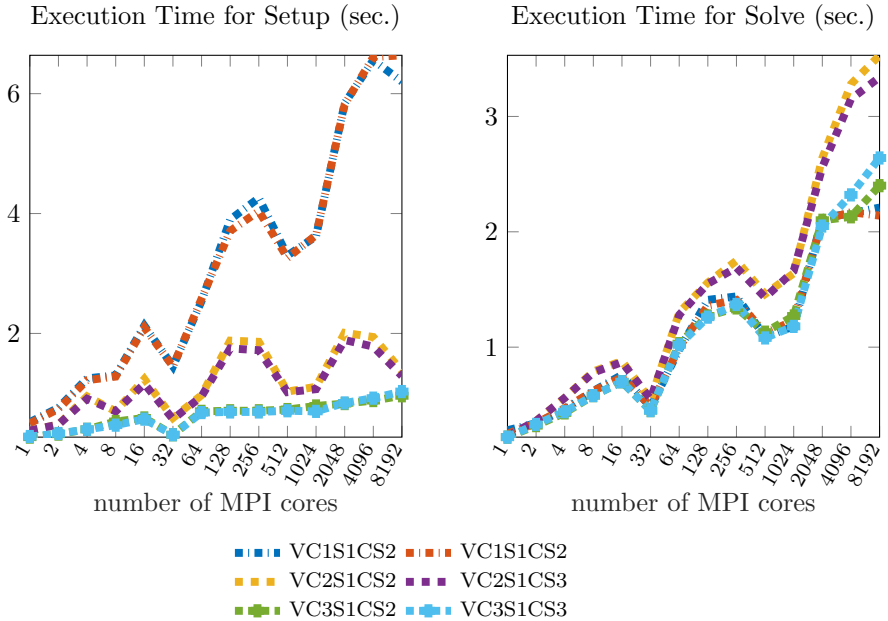
\begin{figure}[htbp]
	\centering
	\subfloat[\label{fig:HypreIt} Number of iterations and Time per iteration of the \texttt{Hypre} preconditioners ]{
%
%
\definecolor{mycolor1}{rgb}{0.00000,0.44700,0.74100}%
\definecolor{mycolor2}{rgb}{0.85000,0.32500,0.09800}%
\definecolor{mycolor3}{rgb}{0.92900,0.69400,0.12500}%
\definecolor{mycolor4}{rgb}{0.49400,0.18400,0.55600}%
\definecolor{mycolor5}{rgb}{0.46600,0.67400,0.18800}%
\definecolor{mycolor6}{rgb}{0.30100,0.74500,0.93300}%
\begin{tikzpicture}

\begin{axis}[%
width=0.38\textwidth,
height=2in,
at={(0in,0in)},
scale only axis,
xmode=log,
xmin=1,
xmax=8192,
xtick={    1,     2,     4,     8,    16,    32,    64,   128,   256,   512,  1024,  2048,  4096,  8192, 16384},
xticklabels={    1,     2,     4,     8,    16,    32,    64,   128,   256,   512,  1024,  2048,  4096,  8192, 16384},
x tick label style={rotate=60,anchor=east,font=\footnotesize},
xminorticks=true,
log basis x=2,
xminorticks=true,
xlabel style={font=\color{white!15!black},at={(0.5,-0.12)}},
xlabel={number of MPI cores},
ymin=6,
ymax=44,
title={Iterations},
axis background/.style={fill=white},
legend columns = 2,
legend style={at={(1.0,-0.25)}, anchor=north, legend cell align=left, align=left, draw=none,font=\footnotesize}
]
\addplot [color=mycolor1, dashdotted, line width=4.0pt]
  table[row sep=crcr]{%
1	6\\
2	7\\
4	7\\
8	8\\
16	8\\
32	8\\
64	9\\
128	10\\
256	10\\
512	11\\
1024	12\\
2048	13\\
4096	13\\
8192	14\\
};
\addlegendentry{VC1S1CS2}

\addplot [color=mycolor2, dashdotted, line width=4.0pt]
  table[row sep=crcr]{%
1	6\\
2	7\\
4	7\\
8	8\\
16	8\\
32	8\\
64	9\\
128	10\\
256	10\\
512	11\\
1024	12\\
2048	13\\
4096	13\\
8192	14\\
};
\addlegendentry{VC1S1CS2}

\addplot [color=mycolor3, dashed, line width=3.0pt]
  table[row sep=crcr]{%
1	6\\
2	9\\
4	12\\
8	13\\
16	14\\
32	14\\
64	16\\
128	18\\
256	20\\
512	22\\
1024	25\\
2048	29\\
4096	35\\
8192	37\\
};
\addlegendentry{VC2S1CS2}

\addplot [color=mycolor4, dashed, line width=3.0pt]
  table[row sep=crcr]{%
1	6\\
2	9\\
4	12\\
8	13\\
16	14\\
32	14\\
64	16\\
128	18\\
256	20\\
512	22\\
1024	25\\
2048	29\\
4096	35\\
8192	37\\
};
\addlegendentry{VC2S1CS3}

\addplot [color=mycolor5, dotted, line width=3.0pt, mark size=2.5pt, mark=+, mark options={solid, mycolor5}]
  table[row sep=crcr]{%
1	12\\
2	15\\
4	17\\
8	17\\
16	20\\
32	20\\
64	22\\
128	25\\
256	27\\
512	29\\
1024	32\\
2048	37\\
4096	40\\
8192	44\\
};
\addlegendentry{VC3S1CS2}

\addplot [color=mycolor6, dotted, line width=3.0pt, mark size=2.5pt, mark=+, mark options={solid, mycolor6}]
  table[row sep=crcr]{%
1	12\\
2	15\\
4	17\\
8	17\\
16	20\\
32	20\\
64	22\\
128	25\\
256	27\\
512	29\\
1024	32\\
2048	37\\
4096	40\\
8192	44\\
};
\addlegendentry{VC3S1CS3}

\end{axis}

\begin{axis}[%
width=0.38\textwidth,
height=2in,
at={(2.5in,0in)},
scale only axis,
xmode=log,
xmin=1,
xmax=8192,
xtick={    1,     2,     4,     8,    16,    32,    64,   128,   256,   512,  1024,  2048,  4096,  8192, 16384},
xticklabels={    1,     2,     4,     8,    16,    32,    64,   128,   256,   512,  1024,  2048,  4096,  8192, 16384},
x tick label style={rotate=60,anchor=east,font=\footnotesize},
xminorticks=true,
log basis x=2,
xminorticks=true,
xlabel style={font=\color{white!15!black},at={(0.5,-0.12)}},
xlabel={number of MPI cores},
ymin=0.0183333333333333,
ymax=0.166923076923077,
title={Time per iteration (sec.)},
axis background/.style={fill=white},
]
\addplot [color=mycolor1, dashdotted, line width=3.0pt]
  table[row sep=crcr]{%
1	0.0466666666666667\\
2	0.0485714285714286\\
4	0.0657142857142857\\
8	0.0775\\
16	0.095\\
32	0.06\\
64	0.108888888888889\\
128	0.141\\
256	0.144\\
512	0.0990909090909091\\
1024	0.0983333333333333\\
2048	0.163076923076923\\
4096	0.166153846153846\\
8192	0.157857142857143\\
};

\addplot [color=mycolor2, dashdotted, line width=3.0pt]
  table[row sep=crcr]{%
1	0.0416666666666667\\
2	0.0485714285714286\\
4	0.0671428571428571\\
8	0.0775\\
16	0.0925\\
32	0.06125\\
64	0.111111111111111\\
128	0.136\\
256	0.141\\
512	0.101818181818182\\
1024	0.1025\\
2048	0.161538461538462\\
4096	0.166923076923077\\
8192	0.152857142857143\\
};

\addplot [color=mycolor3, dashed, line width=3.0pt]
  table[row sep=crcr]{%
1	0.0366666666666667\\
2	0.04\\
4	0.0466666666666667\\
8	0.0607692307692308\\
16	0.0628571428571429\\
32	0.0385714285714286\\
64	0.08125\\
128	0.0866666666666667\\
256	0.0875\\
512	0.0663636363636364\\
1024	0.0656\\
2048	0.0910344827586207\\
4096	0.0937142857142857\\
8192	0.0954054054054054\\
};

\addplot [color=mycolor4, dashed, line width=3.0pt]
  table[row sep=crcr]{%
1	0.0383333333333333\\
2	0.04\\
4	0.0466666666666667\\
8	0.06\\
16	0.0621428571428571\\
32	0.0392857142857143\\
64	0.08\\
128	0.0861111111111111\\
256	0.0845\\
512	0.065\\
1024	0.0664\\
2048	0.0882758620689655\\
4096	0.09\\
8192	0.0902702702702703\\
};

\addplot [color=mycolor5, dotted, line width=3.0pt, mark size=2.5pt, mark=+, mark options={solid, mycolor5}]
  table[row sep=crcr]{%
1	0.0183333333333333\\
2	0.0213333333333333\\
4	0.0252941176470588\\
8	0.0341176470588235\\
16	0.035\\
32	0.0225\\
64	0.0468181818181818\\
128	0.0508\\
256	0.0496296296296296\\
512	0.0389655172413793\\
1024	0.0396875\\
2048	0.0567567567567568\\
4096	0.05325\\
8192	0.0545454545454545\\
};

\addplot [color=mycolor6, dotted, line width=3.0pt, mark size=2.5pt, mark=+, mark options={solid, mycolor6}]
  table[row sep=crcr]{%
1	0.0183333333333333\\
2	0.022\\
4	0.0258823529411765\\
8	0.0341176470588235\\
16	0.035\\
32	0.0225\\
64	0.0463636363636364\\
128	0.0504\\
256	0.0507407407407407\\
512	0.0372413793103448\\
1024	0.036875\\
2048	0.0554054054054054\\
4096	0.058\\
8192	0.06\\
};

\end{axis}

\end{tikzpicture}%
}
	
	\subfloat[\label{fig:HypreTime} Setup Time and Solve Time of the \texttt{Hypre} Preconditioners ]{
%
%
\definecolor{mycolor1}{rgb}{0.00000,0.44700,0.74100}%
\definecolor{mycolor2}{rgb}{0.85000,0.32500,0.09800}%
\definecolor{mycolor3}{rgb}{0.92900,0.69400,0.12500}%
\definecolor{mycolor4}{rgb}{0.49400,0.18400,0.55600}%
\definecolor{mycolor5}{rgb}{0.46600,0.67400,0.18800}%
\definecolor{mycolor6}{rgb}{0.30100,0.74500,0.93300}%
\begin{tikzpicture}

\begin{axis}[%
width=0.38\textwidth,
height=2in,
at={(0in,0in)},
scale only axis,
xmode=log,
xmin=1,
xmax=8192,
xtick={    1,     2,     4,     8,    16,    32,    64,   128,   256,   512,  1024,  2048,  4096,  8192, 16384},
xticklabels={    1,     2,     4,     8,    16,    32,    64,   128,   256,   512,  1024,  2048,  4096,  8192, 16384},
x tick label style={rotate=60,anchor=east,font=\footnotesize},
xminorticks=true,
log basis x=2,
xminorticks=true,
xlabel style={font=\color{white!15!black},at={(0.5,-0.12)}},
xlabel={number of MPI cores},
ymin=0.27,
ymax=6.64,
yminorticks=true,
title={Execution Time for Setup (sec.)},
axis background/.style={fill=white},
legend columns = 2,
legend style={at={(1.0,-0.25)}, anchor=north, legend cell align=left, align=left, draw=none,font=\footnotesize}
]
\addplot [color=mycolor1, dashdotted, line width=3.0pt]
  table[row sep=crcr]{%
1	0.52\\
2	0.76\\
4	1.25\\
8	1.3\\
16	2.13\\
32	1.45\\
64	2.54\\
128	3.9\\
256	4.26\\
512	3.3\\
1024	3.62\\
2048	5.84\\
4096	6.56\\
8192	6.2\\
};
\addlegendentry{VC1S1CS2}

\addplot [color=mycolor2, dashdotted, line width=3.0pt]
  table[row sep=crcr]{%
1	0.5\\
2	0.74\\
4	1.22\\
8	1.28\\
16	2.1\\
32	1.45\\
64	2.58\\
128	3.71\\
256	4.04\\
512	3.26\\
1024	3.64\\
2048	5.82\\
4096	6.61\\
8192	6.64\\
};
\addlegendentry{VC1S1CS2}

\addplot [color=mycolor3, dashed, line width=3.0pt]
  table[row sep=crcr]{%
1	0.38\\
2	0.48\\
4	0.95\\
8	0.7\\
16	1.24\\
32	0.59\\
64	0.97\\
128	1.88\\
256	1.86\\
512	1.03\\
1024	1.11\\
2048	2.01\\
4096	1.94\\
8192	1.37\\
};
\addlegendentry{VC2S1CS2}

\addplot [color=mycolor4, dashed, line width=3.0pt]
  table[row sep=crcr]{%
1	0.38\\
2	0.48\\
4	0.91\\
8	0.7\\
16	1.16\\
32	0.57\\
64	0.95\\
128	1.75\\
256	1.72\\
512	1.02\\
1024	1.07\\
2048	1.89\\
4096	1.78\\
8192	1.29\\
};
\addlegendentry{VC2S1CS3}

\addplot [color=mycolor5, dotted, line width=3.0pt, mark size=2.5pt, mark=+, mark options={solid, mycolor5}]
  table[row sep=crcr]{%
1	0.27\\
2	0.32\\
4	0.4\\
8	0.53\\
16	0.59\\
32	0.31\\
64	0.7\\
128	0.71\\
256	0.71\\
512	0.73\\
1024	0.79\\
2048	0.83\\
4096	0.88\\
8192	0.96\\
};
\addlegendentry{VC3S1CS2}

\addplot [color=mycolor6, dotted, line width=3.0pt, mark size=2.5pt, mark=+, mark options={solid, mycolor6}]
  table[row sep=crcr]{%
1	0.28\\
2	0.33\\
4	0.4\\
8	0.47\\
16	0.57\\
32	0.31\\
64	0.68\\
128	0.69\\
256	0.69\\
512	0.71\\
1024	0.7\\
2048	0.84\\
4096	0.92\\
8192	1.03\\
};
\addlegendentry{VC3S1CS3}

\end{axis}
\begin{axis}[%
width=0.38\textwidth,
height=2in,
at={(2.5in,0in)},
scale only axis,
xmode=log,
xmin=1,
xmax=8192,
xtick={    1,     2,     4,     8,    16,    32,    64,   128,   256,   512,  1024,  2048,  4096,  8192, 16384},
xticklabels={    1,     2,     4,     8,    16,    32,    64,   128,   256,   512,  1024,  2048,  4096,  8192, 16384},
x tick label style={rotate=60,anchor=east,font=\footnotesize},
xminorticks=true,
log basis x=2,
xminorticks=true,
xlabel style={font=\color{white!15!black},at={(0.5,-0.12)}},
xlabel={number of MPI cores},
ymin=0.22,
ymax=3.53,
title={Execution Time for Solve (sec.)},
axis background/.style={fill=white},
]
\addplot [color=mycolor1, dashdotted, line width=3.0pt]
  table[row sep=crcr]{%
1	0.28\\
2	0.34\\
4	0.46\\
8	0.62\\
16	0.76\\
32	0.48\\
64	0.98\\
128	1.41\\
256	1.44\\
512	1.09\\
1024	1.18\\
2048	2.12\\
4096	2.16\\
8192	2.21\\
};

\addplot [color=mycolor2, dashdotted, line width=3.0pt]
  table[row sep=crcr]{%
1	0.25\\
2	0.34\\
4	0.47\\
8	0.62\\
16	0.74\\
32	0.49\\
64	1\\
128	1.36\\
256	1.41\\
512	1.12\\
1024	1.23\\
2048	2.1\\
4096	2.17\\
8192	2.14\\
};

\addplot [color=mycolor3, dashed, line width=3.0pt]
  table[row sep=crcr]{%
1	0.22\\
2	0.36\\
4	0.56\\
8	0.79\\
16	0.88\\
32	0.54\\
64	1.3\\
128	1.56\\
256	1.75\\
512	1.46\\
1024	1.64\\
2048	2.64\\
4096	3.28\\
8192	3.53\\
};

\addplot [color=mycolor4, dashed, line width=3.0pt]
  table[row sep=crcr]{%
1	0.23\\
2	0.36\\
4	0.56\\
8	0.78\\
16	0.87\\
32	0.55\\
64	1.28\\
128	1.55\\
256	1.69\\
512	1.43\\
1024	1.66\\
2048	2.56\\
4096	3.15\\
8192	3.34\\
};

\addplot [color=mycolor5, dotted, line width=3.0pt, mark size=2.5pt, mark=+, mark options={solid, mycolor5}]
  table[row sep=crcr]{%
1	0.22\\
2	0.32\\
4	0.43\\
8	0.58\\
16	0.7\\
32	0.45\\
64	1.03\\
128	1.27\\
256	1.34\\
512	1.13\\
1024	1.27\\
2048	2.1\\
4096	2.13\\
8192	2.4\\
};

\addplot [color=mycolor6, dotted, line width=3.0pt, mark size=2.5pt, mark=+, mark options={solid, mycolor6}]
  table[row sep=crcr]{%
1	0.22\\
2	0.33\\
4	0.44\\
8	0.58\\
16	0.7\\
32	0.45\\
64	1.02\\
128	1.26\\
256	1.37\\
512	1.08\\
1024	1.18\\
2048	2.05\\
4096	2.32\\
8192	2.64\\
};

\end{axis}
\end{tikzpicture}%
} 
	
	\caption{Weak scaling results for  256k  dofs  per  core. Comparison of the convergence behavior and execution time of the \texttt{Hypre} preconditioners, when direct or iterative coarsest solver is applied.}
	\label{fig:hypreitsvsds}
\end{figure}

\subsection{{Anisotropic Poisson problem}} {To complement the experiment in Section~\ref{sec:testcase}, we consider here also the case of an anisotropic Poisson problem on a structured grid, i.e., we replace \eqref{eq:bench} with 
	\begin{equation*}
	\sum_{i=1}^{3} k_i \frac{\partial^2 u}{\partial x_i^2}= \mathbf{f}, \quad\text{ in } [0,1]^3.
	\end{equation*}
}{
	Let us observe that in general the weighted matching algorithm is able of recovering the aggregates one would expect for the different type of anisotropies, consider, e.g., the case of z-aligned anisotropy, in which $k_3 = 100$ and $k_1=k_2=1$ represented in Figure~\ref{fig:aggregate_anisotropy}.
	\begin{figure}[htbp]
		\centering
		\includegraphics[width=\columnwidth]{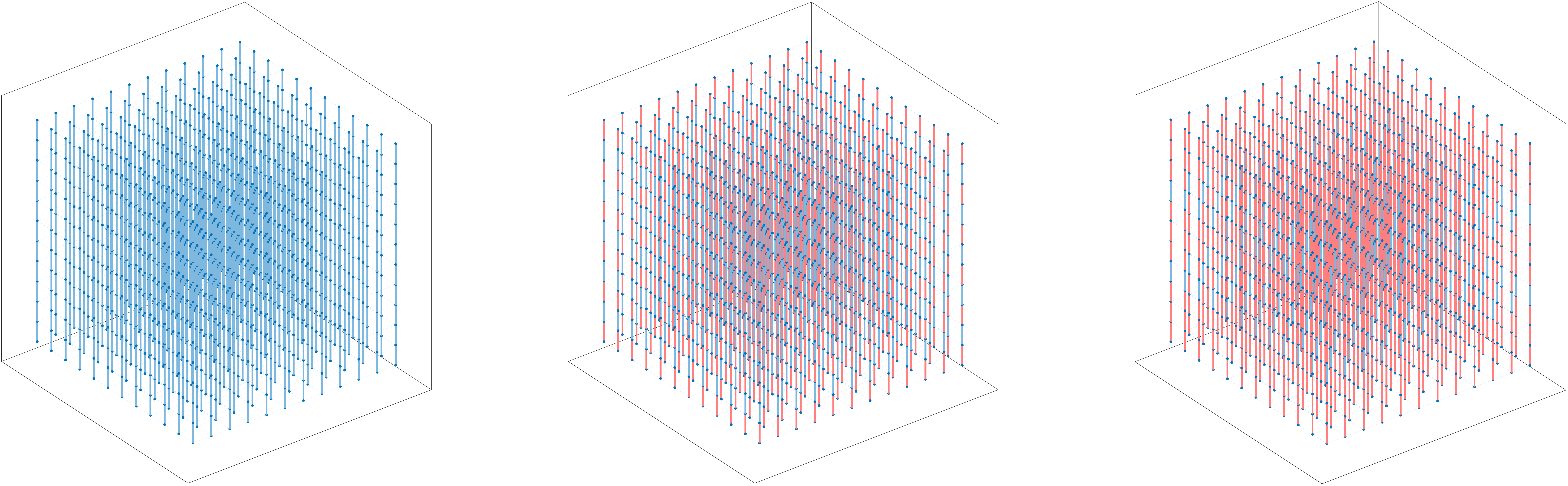}
		\caption{Anisotropic Poisson problem with $z$-oriented axial anisotropy. The width of the edges is proportional to the weight of the corresponding edge. Edges that are part of the matching are highlighted in red, while unmatched edges are reported in cyan. The graph in the central panel corresponds to the case of a single sweep of matching, while the one in the last panel corresponds to the case of two matching sweeps.}
		\label{fig:aggregate_anisotropy}
	\end{figure}
	We investigate then the behavior of the preconditioners on the test case for $k_1 = k_2 = 100$, and $k_3 = 1$ while using again $512k$ dofs per node. In Fig.~\ref{fig:aniso_opcplx} we show the operator complexity of the multilevel hierarchies corresponding to all the considered preconditioners.
	\begin{figure}[htbp]
		\centering
		\subfloat[\label{fig:aniso_opcplx}]{
%
%
\definecolor{mycolor1}{rgb}{0.00000,0.44700,0.74100}%
\definecolor{mycolor2}{rgb}{0.85000,0.32500,0.09800}%
\definecolor{mycolor3}{rgb}{0.92900,0.69400,0.12500}%
\definecolor{mycolor4}{rgb}{0.49400,0.18400,0.55600}%
\begin{tikzpicture}

\begin{axis}[%
width=0.40\textwidth,
height=2.1in,
at={(2.509in,2.689in)},
scale only axis,
xmode=log,
xmin=1,
xmax=4096,
xtick={    1,     2,     4,     8,    16,    32,    64,   128,   256,   512,  1024,  2048,  4096,  8192, 16384, 28800},
xticklabels={    1,     2,     4,     8,    16,    32,    64,   128,   256,   512,  1024,  2048,  4096,  8192, 16384, 28800},
x tick label style={rotate=60,anchor=east,font=\footnotesize},
xminorticks=true,
xlabel style={font=\color{white!15!black}},
xlabel={number of MPI cores},
ymin=0,
ymax=4.284957,
title={Operator Complexity},
axis background/.style={fill=white},
legend columns = 2,
legend style={at={(0.5,-0.25)}, anchor=north, legend cell align=left, align=left, draw=none,font=\scriptsize}
]
\addplot [color=mycolor1, dashdotted, line width=3.0pt]
  table[row sep=crcr]{%
1	3.825403\\
2	3.950092\\
4	4.011222\\
8	4.028796\\
16	4.102072\\
32	4.151046\\
64	4.147568\\
128	4.202475\\
256	4.238566\\
512	4.215979\\
1024	4.256825\\
2048	4.284957\\
4096	4.252871\\
};
\addlegendentry{VC1S1CS2}

\addplot [color=mycolor2, dashed, line width=3.0pt]
  table[row sep=crcr]{%
1	3.526049\\
2	3.554471\\
4	3.577735\\
8	3.600243\\
16	3.611048\\
32	3.629618\\
64	3.649799\\
128	3.647961\\
256	3.657762\\
512	3.676675\\
1024	3.666846\\
2048	3.673876\\
4096	3.691714\\
};
\addlegendentry{VC2S1CS2}

\addplot [color=mycolor3, dotted, line width=3.0pt, mark size=2.5pt, mark=+, mark options={solid, mycolor3}]
  table[row sep=crcr]{%
1	1.14\\
2	1.141\\
4	1.141\\
8	1.141\\
16	1.142\\
32	1.142\\
64	1.142\\
128	1.143\\
256	1.142\\
512	1.143\\
1024	1.144\\
2048	1.143\\
4096	1.143\\
};
\addlegendentry{KA1S1CS1}

\addplot [color=mycolor4, dotted, line width=3.0pt, mark size=2.5pt, mark=+, mark options={solid, mycolor4}]
  table[row sep=crcr]{%
1	1.704\\
2	1.726\\
4	1.739\\
8	1.752\\
16	1.767\\
32	1.772\\
64	1.78\\
128	1.792\\
256	1.79\\
512	1.796\\
1024	1.812\\
2048	1.799\\
4096	1.81\\
};
\addlegendentry{VA3S1CS1}

\end{axis}
\end{tikzpicture}%
}
		\subfloat[\label{fig:aniso_iters}]{
%
%
\definecolor{mycolor1}{rgb}{0.00000,0.44700,0.74100}%
\definecolor{mycolor2}{rgb}{0.85000,0.32500,0.09800}%
\definecolor{mycolor3}{rgb}{0.92900,0.69400,0.12500}%
\definecolor{mycolor4}{rgb}{0.49400,0.18400,0.55600}%
\begin{tikzpicture}

\begin{axis}[%
width=0.40\textwidth,
height=2.1in,
at={(2.509in,2.689in)},
scale only axis,
xmode=log,
xmin=1,
xmax=4096,
xtick={   1,    2,    4,    8,   16,   32,   64,  128,  256,  512, 1024, 2048, 4096, 8192},
xticklabels={   1,    2,    4,    8,   16,   32,   64,  128,  256,  512, 1024, 2048, 4096, 8192},
x tick label style={rotate=70,anchor=east,font=\footnotesize},
xlabel style={font=\color{white!15!black}},
xlabel={number of MPI cores},
ymin=0,
ymax=27,
title={Iterations},
axis background/.style={fill=white},
legend columns = 2,
legend style={at={(0.5,-0.25)}, anchor=north, legend cell align=left, align=left, draw=none,font=\scriptsize}
]
\addplot [color=mycolor1, dashdotted, line width=4.0pt]
  table[row sep=crcr]{%
1	6\\
2	6\\
4	7\\
8	7\\
16	7\\
32	8\\
64	8\\
128	9\\
256	9\\
512	10\\
1024	10\\
2048	11\\
4096	12\\
};
\addlegendentry{VC1S1CS2}

\addplot [color=mycolor2, dashed, line width=3.0pt]
  table[row sep=crcr]{%
1	6\\
2	9\\
4	11\\
8	11\\
16	13\\
32	14\\
64	14\\
128	17\\
256	18\\
512	20\\
1024	22\\
2048	25\\
4096	27\\
};
\addlegendentry{VC2S1CS2}

\addplot [color=mycolor3, dotted, line width=3.0pt, mark size=2.5pt, mark=+, mark options={solid, mycolor3}]
  table[row sep=crcr]{%
1	16\\
2	17\\
4	17\\
8	18\\
16	19\\
32	18\\
64	19\\
128	22\\
256	20\\
512	20\\
1024	22\\
2048	23\\
4096	23\\
};
\addlegendentry{KA1S1CS1}

\addplot [color=mycolor4, dotted, line width=3.0pt, mark size=2.5pt, mark=+, mark options={solid, mycolor4}]
  table[row sep=crcr]{%
1	13\\
2	14\\
4	15\\
8	15\\
16	18\\
32	17\\
64	17\\
128	20\\
256	18\\
512	18\\
1024	20\\
2048	18\\
4096	19\\
};
\addlegendentry{VA3S1CS1}

\end{axis}
\end{tikzpicture}%
}
		\caption{Weak Scaling results for 512k dofs per core. Operator complexity and number of iterations.}
	\end{figure}
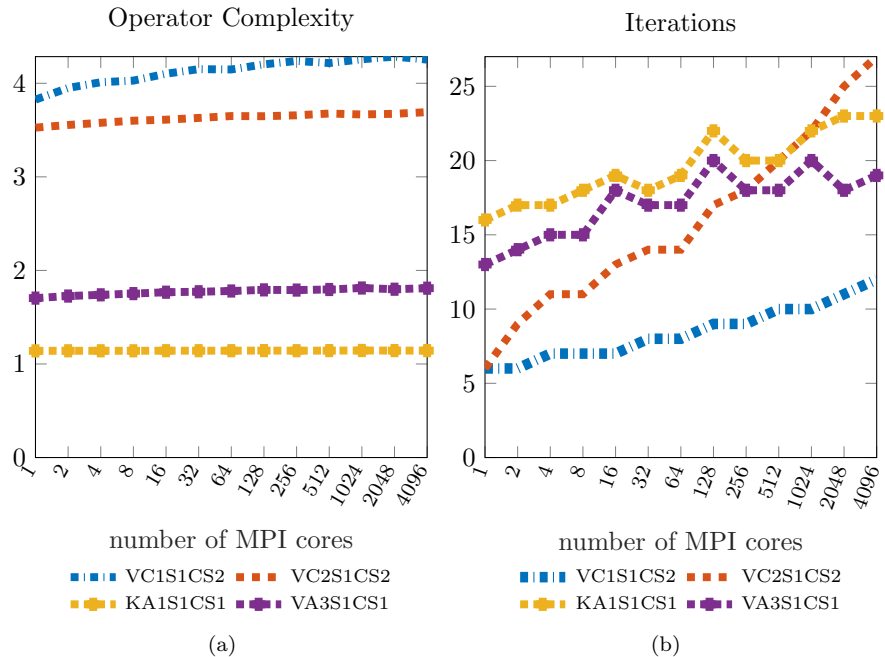
	What we observe is that, consistently with what happens with the isotropic case, both the HMIS and Falgout coarsening schemes produce hierarchies with a significantly higher operator complexity. As seen in Fig.~\ref{fig:aniso_iters}, for Falgout this corresponds to the smaller number of overall iterations for the preconditioned system; however the number of iterations in the HMIS case tends to increase with the number of cores. Let us also stress that we have stopped the analysis at 4096 cores for this case since both \texttt{Hypre} AMG schemes produce out of memory errors for the next step. On the other hand, both of our schemes have a smaller operator complexity and a better behavior with respect to the number of iterations if compared with the HMIS case in \texttt{Hypre}. 
}
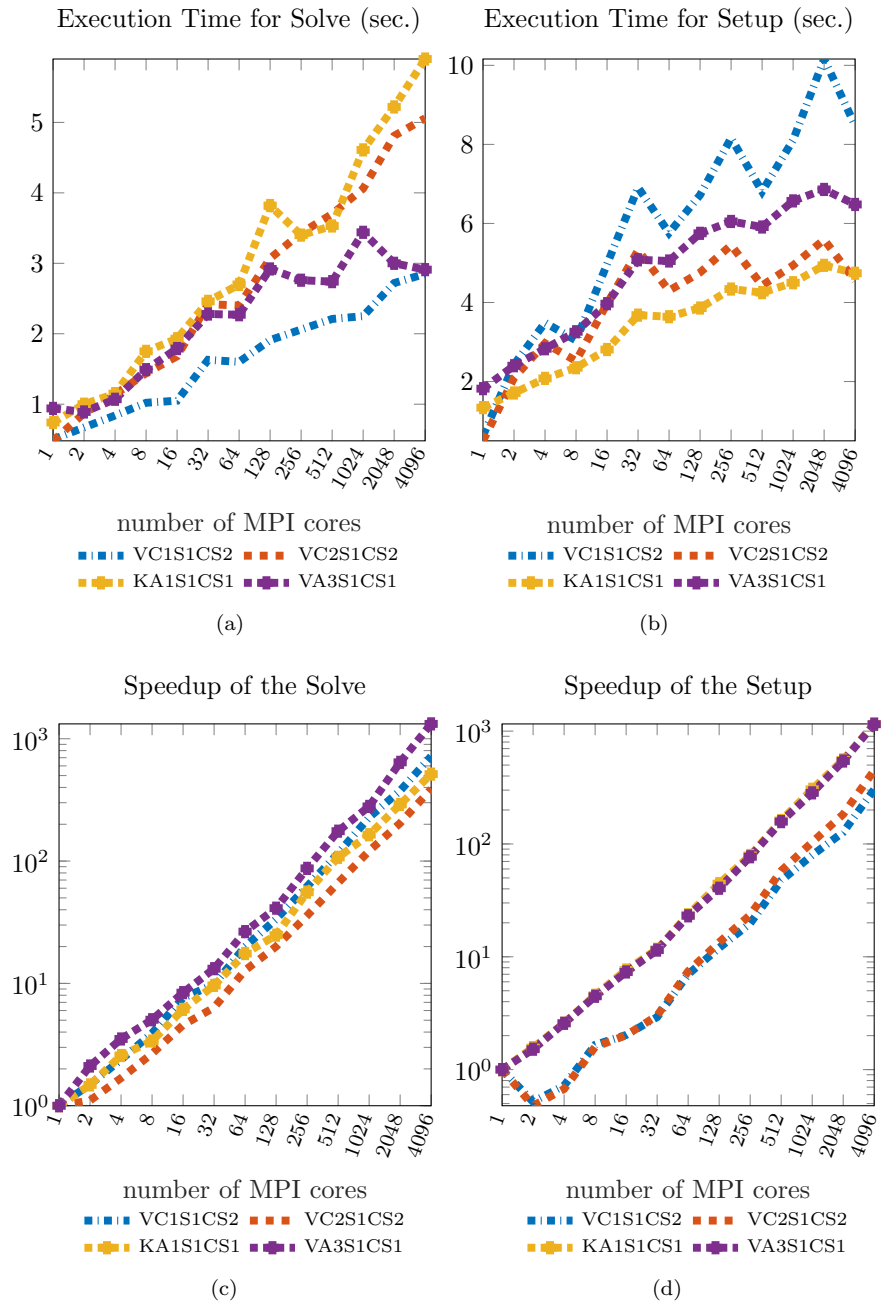
\begin{figure}[htbp]
	\centering
	\subfloat[\label{fig:aniso_solvetime}]{
%
%
\definecolor{mycolor1}{rgb}{0.00000,0.44700,0.74100}%
\definecolor{mycolor2}{rgb}{0.85000,0.32500,0.09800}%
\definecolor{mycolor3}{rgb}{0.92900,0.69400,0.12500}%
\definecolor{mycolor4}{rgb}{0.49400,0.18400,0.55600}%
\begin{tikzpicture}

\begin{axis}[%
width=0.38\textwidth,
height=2in,
at={(2.509in,2.689in)},
scale only axis,
xmode=log,
xmin=1,
xmax=4096,
xtick={    1,     2,     4,     8,    16,    32,    64,   128,   256,   512,  1024,  2048,  4096,  8192, 16384, 28800},
xticklabels={    1,     2,     4,     8,    16,    32,    64,   128,   256,   512,  1024,  2048,  4096,  8192, 16384, 28800},
x tick label style={rotate=60,anchor=east,font=\footnotesize},
xlabel style={font=\color{white!15!black}},
xlabel={number of MPI cores},
ymin=0.48,
ymax=5.9,
title={Execution Time for Solve (sec.)},
axis background/.style={fill=white},
legend columns = 2,
legend style={at={(0.5,-0.25)}, anchor=north, legend cell align=left, align=left, draw=none,font=\scriptsize}
]
\addplot [color=mycolor1, dashdotted, line width=3.0pt]
  table[row sep=crcr]{%
1	0.5\\
2	0.67\\
4	0.84\\
8	1.02\\
16	1.05\\
32	1.63\\
64	1.6\\
128	1.91\\
256	2.06\\
512	2.21\\
1024	2.25\\
2048	2.72\\
4096	2.85\\
};
\addlegendentry{VC1S1CS2}

\addplot [color=mycolor2, dashed, line width=3.0pt]
  table[row sep=crcr]{%
1	0.48\\
2	0.87\\
4	1.14\\
8	1.44\\
16	1.68\\
32	2.42\\
64	2.4\\
128	3.07\\
256	3.43\\
512	3.7\\
1024	4.06\\
2048	4.81\\
4096	5.06\\
};
\addlegendentry{VC2S1CS2}

\addplot [color=mycolor3, dotted, line width=3.0pt, mark size=2.5pt, mark=+, mark options={solid, mycolor3}]
  table[row sep=crcr]{%
1	0.74\\
2	1\\
4	1.15\\
8	1.75\\
16	1.93\\
32	2.46\\
64	2.7\\
128	3.82\\
256	3.4\\
512	3.53\\
1024	4.61\\
2048	5.22\\
4096	5.9\\
};
\addlegendentry{KA1S1CS1}

\addplot [color=mycolor4, dotted, line width=3.0pt, mark size=2.5pt, mark=+, mark options={solid, mycolor4}]
  table[row sep=crcr]{%
1	0.94\\
2	0.89\\
4	1.07\\
8	1.49\\
16	1.79\\
32	2.28\\
64	2.27\\
128	2.92\\
256	2.76\\
512	2.74\\
1024	3.44\\
2048	3\\
4096	2.91\\
};
\addlegendentry{VA3S1CS1}

\end{axis}
\end{tikzpicture}%
}
	\subfloat[\label{fig:aniso_setuptime}]{
%
%
\definecolor{mycolor1}{rgb}{0.00000,0.44700,0.74100}%
\definecolor{mycolor2}{rgb}{0.85000,0.32500,0.09800}%
\definecolor{mycolor3}{rgb}{0.92900,0.69400,0.12500}%
\definecolor{mycolor4}{rgb}{0.49400,0.18400,0.55600}%
\begin{tikzpicture}

\begin{axis}[%
width=0.38\textwidth,
height=2in,
at={(2.509in,2.689in)},
scale only axis,
xmode=log,
xmin=1,
xmax=4096,
xtick={   1,    2,    4,    8,   16,   32,   64,  128,  256,  512, 1024, 2048, 4096, 8192},
xticklabels={   1,    2,    4,    8,   16,   32,   64,  128,  256,  512, 1024, 2048, 4096, 8192},
x tick label style={rotate=70,anchor=east,font=\footnotesize},
xminorticks=true,
xlabel style={font=\color{white!15!black}},
xlabel={number of MPI cores},
ymin=0.5,
ymax=10.16,
title={Execution Time for Setup (sec.)},
axis background/.style={fill=white},
legend columns = 2,
legend style={at={(0.5,-0.25)}, anchor=north, legend cell align=left, align=left, draw=none,font=\scriptsize}
]
\addplot [color=mycolor1, dashdotted, line width=3.0pt]
  table[row sep=crcr]{%
1	0.63\\
2	2.43\\
4	3.49\\
8	3.04\\
16	4.97\\
32	6.91\\
64	5.76\\
128	6.72\\
256	8.16\\
512	6.81\\
1024	8.11\\
2048	10.16\\
4096	8.48\\
};
\addlegendentry{VC1S1CS2}

\addplot [color=mycolor2, dashed, line width=3.0pt]
  table[row sep=crcr]{%
1	0.5\\
2	2.09\\
4	2.97\\
8	2.52\\
16	3.97\\
32	5.32\\
64	4.3\\
128	4.76\\
256	5.45\\
512	4.43\\
1024	4.94\\
2048	5.59\\
4096	4.55\\
};
\addlegendentry{VC2S1CS2}

\addplot [color=mycolor3, dotted, line width=3.0pt, mark size=2.5pt, mark=+, mark options={solid, mycolor3}]
  table[row sep=crcr]{%
1	1.34\\
2	1.71\\
4	2.08\\
8	2.35\\
16	2.81\\
32	3.68\\
64	3.64\\
128	3.86\\
256	4.34\\
512	4.25\\
1024	4.5\\
2048	4.94\\
4096	4.74\\
};
\addlegendentry{KA1S1CS1}

\addplot [color=mycolor4, dotted, line width=3.0pt, mark size=2.5pt, mark=+, mark options={solid, mycolor4}]
  table[row sep=crcr]{%
1	1.82\\
2	2.4\\
4	2.83\\
8	3.26\\
16	3.97\\
32	5.08\\
64	5.05\\
128	5.75\\
256	6.05\\
512	5.91\\
1024	6.57\\
2048	6.86\\
4096	6.48\\
};
\addlegendentry{VA3S1CS1}

\end{axis}
\end{tikzpicture}%
}
	
	\subfloat[\label{fig:aniso_solvesp}]{
%
%
\definecolor{mycolor1}{rgb}{0.00000,0.44700,0.74100}%
\definecolor{mycolor2}{rgb}{0.85000,0.32500,0.09800}%
\definecolor{mycolor3}{rgb}{0.92900,0.69400,0.12500}%
\definecolor{mycolor4}{rgb}{0.49400,0.18400,0.55600}%
\begin{tikzpicture}

\begin{axis}[%
width=0.38\textwidth,
height=2in,
at={(2.509in,2.689in)},
scale only axis,
xmode=log,
xmin=1,
xmax=4096,
xtick={   1,    2,    4,    8,   16,   32,   64,  128,  256,  512, 1024, 2048, 4096, 8192},
xticklabels={   1,    2,    4,    8,   16,   32,   64,  128,  256,  512, 1024, 2048, 4096, 8192},
x tick label style={rotate=70,anchor=east,font=\footnotesize},
xminorticks=true,
xlabel style={font=\color{white!15!black}},
xlabel={number of MPI cores},
ymode=log,
ymin=1,
ymax=1323.10652920962,
yminorticks=true,
title={Speedup of the Solve},
axis background/.style={fill=white},
legend columns = 2,
legend style={at={(0.5,-0.25)}, anchor=north, legend cell align=left, align=left, draw=none,font=\scriptsize}
]
\addplot [color=mycolor1, dashdotted, line width=3.0pt]
  table[row sep=crcr]{%
1	1\\
2	1.49253731343284\\
4	2.38095238095238\\
8	3.92156862745098\\
16	7.61904761904762\\
32	9.8159509202454\\
64	20\\
128	33.5078534031414\\
256	62.1359223300971\\
512	115.837104072398\\
1024	227.555555555556\\
2048	376.470588235294\\
4096	718.59649122807\\
};
\addlegendentry{VC1S1CS2}

\addplot [color=mycolor2, dashed, line width=3.0pt]
  table[row sep=crcr]{%
1	1\\
2	1.10344827586207\\
4	1.68421052631579\\
8	2.66666666666667\\
16	4.57142857142857\\
32	6.34710743801653\\
64	12.8\\
128	20.0130293159609\\
256	35.8250728862974\\
512	66.4216216216216\\
1024	121.064039408867\\
2048	204.37422037422\\
4096	388.553359683794\\
};
\addlegendentry{VC2S1CS2}

\addplot [color=mycolor3, dotted, line width=3.0pt, mark size=2.5pt, mark=+, mark options={solid, mycolor3}]
  table[row sep=crcr]{%
1	1\\
2	1.48\\
4	2.57391304347826\\
8	3.38285714285714\\
16	6.13471502590674\\
32	9.6260162601626\\
64	17.5407407407407\\
128	24.7958115183246\\
256	55.7176470588235\\
512	107.331444759207\\
1024	164.373101952278\\
2048	290.329501915709\\
4096	513.735593220339\\
};
\addlegendentry{KA1S1CS1}

\addplot [color=mycolor4, dotted, line width=3.0pt, mark size=2.5pt, mark=+, mark options={solid, mycolor4}]
  table[row sep=crcr]{%
1	1\\
2	2.1123595505618\\
4	3.51401869158878\\
8	5.04697986577181\\
16	8.40223463687151\\
32	13.1929824561404\\
64	26.5022026431718\\
128	41.2054794520548\\
256	87.1884057971015\\
512	175.649635036496\\
1024	279.813953488372\\
2048	641.706666666667\\
4096	1323.10652920962\\
};
\addlegendentry{VA3S1CS1}

\end{axis}
\end{tikzpicture}%
}
	\subfloat[\label{fig:aniso_setupsp}]{
%
%
\definecolor{mycolor1}{rgb}{0.00000,0.44700,0.74100}%
\definecolor{mycolor2}{rgb}{0.85000,0.32500,0.09800}%
\definecolor{mycolor3}{rgb}{0.92900,0.69400,0.12500}%
\definecolor{mycolor4}{rgb}{0.49400,0.18400,0.55600}%
\begin{tikzpicture}

\begin{axis}[%
width=0.38\textwidth,
height=2in,
at={(2.509in,2.689in)},
scale only axis,
xmode=log,
xmin=1,
xmax=4096,
xtick={   1,    2,    4,    8,   16,   32,   64,  128,  256,  512, 1024, 2048, 4096, 8192},
xticklabels={   1,    2,    4,    8,   16,   32,   64,  128,  256,  512, 1024, 2048, 4096, 8192},
x tick label style={rotate=70,anchor=east,font=\footnotesize},
xminorticks=true,
xlabel style={font=\color{white!15!black}},
xlabel={number of MPI cores},
ymode=log,
ymin=0.478468899521531,
ymax=1157.94092827004,
yminorticks=true,
title={Speedup of the Setup},
axis background/.style={fill=white},
legend columns = 2,
legend style={at={(0.5,-0.25)}, anchor=north, legend cell align=left, align=left, draw=none,font=\scriptsize}
]
\addplot [color=mycolor1, dashdotted, line width=3.0pt]
  table[row sep=crcr]{%
1	1\\
2	0.518518518518518\\
4	0.722063037249284\\
8	1.65789473684211\\
16	2.02816901408451\\
32	2.91751085383502\\
64	7\\
128	12\\
256	19.7647058823529\\
512	47.3656387665198\\
1024	79.5462392108508\\
2048	126.992125984252\\
4096	304.301886792453\\
};
\addlegendentry{VC1S1CS2}

\addplot [color=mycolor2, dashed, line width=3.0pt]
  table[row sep=crcr]{%
1	1\\
2	0.478468899521531\\
4	0.673400673400673\\
8	1.58730158730159\\
16	2.01511335012594\\
32	3.00751879699248\\
64	7.44186046511628\\
128	13.4453781512605\\
256	23.4862385321101\\
512	57.7878103837472\\
1024	103.643724696356\\
2048	183.184257602862\\
4096	450.10989010989\\
};
\addlegendentry{VC2S1CS2}

\addplot [color=mycolor3, dotted, line width=3.0pt, mark size=2.5pt, mark=+, mark options={solid, mycolor3}]
  table[row sep=crcr]{%
1	1\\
2	1.5672514619883\\
4	2.57692307692308\\
8	4.56170212765957\\
16	7.62989323843416\\
32	11.6521739130435\\
64	23.5604395604396\\
128	44.4352331606218\\
256	79.0414746543779\\
512	161.430588235294\\
1024	304.924444444444\\
2048	555.53036437247\\
4096	1157.94092827004\\
};
\addlegendentry{KA1S1CS1}

\addplot [color=mycolor4, dotted, line width=3.0pt, mark size=2.5pt, mark=+, mark options={solid, mycolor4}]
  table[row sep=crcr]{%
1	1\\
2	1.51666666666667\\
4	2.57243816254417\\
8	4.46625766871166\\
16	7.33501259445844\\
32	11.4645669291339\\
64	23.0653465346535\\
128	40.5147826086957\\
256	77.0115702479339\\
512	157.671742808799\\
1024	283.665144596651\\
2048	543.34693877551\\
4096	1150.41975308642\\
};
\addlegendentry{VA3S1CS1}

\end{axis}
\end{tikzpicture}%
}
	\caption{Weak scaling results 512k dofs per core. Setup and solve time with the relative speed-up.}\label{fig:anisocpu}
\end{figure}
{
	When we turn to the solve times in Fig.~\ref{fig:aniso_solvetime}, we observe that the KA1S1CS1 and VC2S1CS2 have comparable solve time; indeed what is happening in this case is that the interplay between the number of iterations and the cost of applying the K-Cycle in our KA1S1CS1 balances the higher operator complexity and the cost of using a V-Cycle in VC2S1CS2.  On the other hand, as it already seen  in the isotropic case, our VA3S1CS1 achieves a solve time that is comparable with the \texttt{Hypre} scheme based on the Falgout coarsening strategy VC1S1CS2, especially for increasing number of cores. Clearly, the smaller solve time of the latter for a reduced number of cores comes at the cost of a higher setup time, as it can be seen by comparing the two relative lines in Fig.~\ref{fig:aniso_setuptime}. The KA1S1CS1 and VC2S1CS2 have again comparable setup time, with the first one being marginally smaller. It is interesting to note that the speedups obtained in the solve phase are very similar for all methods, and behave consistently with the  number of iterations and application cost; see Fig.~\ref{fig:aniso_solvesp}. On the other hand, in Fig.~\ref{fig:aniso_setupsp}, the speedup for the setup phase displays a smooth increase with  increasing number of cores for all  methods, with a slight better behavior for  our proposed method, demonstrating again the good implementation scalability of all computational kernels.
}

\section{Performance results towards extreme scale} The resulting performance of the multigrid preconditioner in term of implementation scalability depends also on how effective the coarsening procedure is, and on how well balanced is the distribution of the coarsest matrix. We report in Tab.~\ref{tab:averagecr} the average coarsening ratio $\overline{{cr}}$ obtained for the Weak Scaling test with 512k dofs per node, that is, the average of the coarsening ratio ${cr}_l$ at each level $l$ defined as $\frac{1}{\text{nl}}\sum_{l=0}^{\text{nl}-1}{cr}_l = \frac{1}{\text{nl}} \sum_{l=0}^{\text{nl}-1} n_l/n_{l+1}$, with $n_l$ the matrix size at the level $l$.
\begin{table}[htbp]
	\centering
	\caption{Average coarsening ratio $\overline{{cr}}$ for the weak scalability test with 512k dofs per node while using the Parallel Matching with 3 and 4 sweeps coarsening strategy, and information on the block-distribution of the coarsest matrix.}
	
	\subfloat[average coarsening ratio\label{tab:averagecr}]{
		\begin{tabular}{rllll}
			\toprule
			Aggregation & \multicolumn{2}{c}{PMC3} & \multicolumn{2}{c}{PMC4} \\
			\cmidrule{2-3} \cmidrule{4-5}
			np & $\text{nl}$ & $\overline{{cr}}$ & $\text{nl}$ & $\overline{{cr}}$  \\
			\midrule
			1 & 5 & 8.00 & 4 & 16.00 \\
			2 & 5 & 7.98 & 4 & 15.95 \\
			4 & 5 & 8.00 & 4 & 16.00 \\
			8 & 5 & 8.00 & 4 & 16.00 \\
			16 & 5 & 7.95 & 4 & 15.86 \\
			32 & 5 & 8.00 & 4 & 16.00 \\
			64 & 5 & 8.00 & 4 & 16.00 \\
			128 & 5 & 7.93 & 4 & 15.81 \\
			256 & 5 & 8.00 & 4 & 16.00 \\
			512 & 5 & 8.00 & 4 & 16.00 \\
			1024 & 5 & 7.96 & 4 & 15.90 \\
			2048 & 5 & 8.00 & 4 & 16.00 \\
			4096 & 5 & 7.94 & 4 & 15.84 \\
			8192 & 5 & 7.96 & 4 & 15.90 \\
			16384 & 5 & 8.00 & 4 & 16.00 \\
			27000 & 5 & 8.00 & 4 & 16.00 \\
			\bottomrule
	\end{tabular}}\hfill
	\subfloat[Coarsest matrix block distribution\label{tab:coarsest_blocks}]{    
		\begin{tabular}{rcccccc}%
			\toprule
			&  \multicolumn{2}{c}{block size} & coarsest \\
			\cmidrule{2-3}%
			np &  min & max & size \\
			\midrule
			1     & 125 & 125 & 125 \\
			2     & 138 & 139 & 277 \\
			4     & 128 & 128 & 512 \\
			8     & 125 & 125 & 1000 \\
			16    & 121 & 129 & 2004 \\
			32    & 128 & 128 & 4096 \\
			64    & 125 & 125 & 8000 \\
			128   & 119 & 139 & 16206 \\
			256   & 128 & 128 & 32768 \\
			512   & 125 & 125 & 64000 \\
			1024  & 123 & 162 & 127466 \\
			2048  & 128 & 128 & 262144 \\
			4096  & 106 & 117 & 452335 \\
			8192  & 115 & 133 & 1019452 \\
			16384 & 128 & 128 & 2097152 \\
			27000 & 125 & 125 & 3375000 \\
			\bottomrule
	\end{tabular}}
\end{table}
Moreover, in Tab.~\ref{tab:coarsest_blocks} we report the size of the coarsest matrix together with the minimum and maximum block size. The size of the blocks and the minimum/maximum block sizes for the coarse matrix are the same for both the parallel matching with three (PMC3) and four sweeps (PMC4) aggregation strategies. The difference is that we reach this target dimension with a hierarchy with $nl=4$ instead of $nl=5$, and this  is reflected in the average coarsening ratios in Tab.~\ref{tab:averagecr}. As we can observe, we achieve a good load-balancing among the different cores, i.e., there is a small variation between the size of the largest and smallest block size on each processor. {Observe moreover that we can exactly compute an upper bound of the cost of the application of 1 K-cycle, so that, as also described in~\cite{NV2008}, the complexity for our K-cycle can be defined in terms of its operator complexity: $\text{K-opc}=\frac{\sum_{l=0}^{\text{nl}-1} 2^{l} \operatorname{nnz}(A_{l})}{\operatorname{nnz}(A_0)}$.  This takes into account the number of iterations ($it=2$) at each level and the maximum number of times each level is visited in the recursion ($2^{nl-l-1}$). This complexity is then generally larger than that of the V-cycle although a still linearly increasing function in the non-zero entries of the system matrix (while the increase factor of the V-cycle is $1/(1-q)$, where $q$ is the inverse of the coarsening ratio, the increase factor of the K-cycle is then $1/(1-2q)$.}

\begin{figure}
	\centering
	\subfloat[\label{fig:VSCycle} VA3S1CS1 preconditioner on MPI and VA3S5CS1 on MPI-CUDA ]{
%
%
\definecolor{mycolor1}{rgb}{0.00000,0.44700,0.74100}%
\definecolor{mycolor2}{rgb}{0.85000,0.32500,0.09800}%
\definecolor{mycolor3}{rgb}{0.92900,0.69400,0.12500}%
\definecolor{mycolor4}{rgb}{0.49400,0.18400,0.55600}%
\definecolor{mycolor5}{rgb}{0.46600,0.67400,0.18800}%
\begin{tikzpicture}

\begin{axis}[%
width=0.87\textwidth,
height=2.8in,
at={(2.509in,2.768in)},
scale only axis,
xmode=log,
xmin=100000,
xmax=100000000000,
xminorticks=true,
xlabel style={font=\color{white!15!black}},
xlabel={dofs},
ymode=log,
ymin=0.1,
ymax=10,
yminorticks=true,
title={Execution Time for Solve (sec.)},
axis background/.style={fill=white},
ymajorgrids,
yminorgrids,
legend columns=3,
legend style={at={(0.5,-0.15)}, anchor=north, legend cell align=left, align=left, draw=none,font=\footnotesize}
]
\addplot [color=mycolor1, line width=3.0pt, mark=o, mark options={solid, mycolor1}]
  table[row sep=crcr]{%
6144000	0.28\\
12288000	0.38\\
24576000	0.52\\
49152000	0.65\\
98304000	0.79\\
196608000	0.92\\
393216000	1.29\\
786432000	1.37\\
1572864000	1.48\\
3145728000	2.73\\
6291456000	2.62\\
12582912000	2.85\\
};
\addlegendentry{6M dof x GPU}

\addplot [color=mycolor2, line width=3.0pt, mark=o, mark options={solid, mycolor2}]
  table[row sep=crcr]{%
3072000	0.16\\
6144000	0.22\\
12288000	0.28\\
24576000	0.49\\
49152000	0.56\\
98304000	0.83\\
196608000	0.95\\
393216000	1.22\\
786432000	1.4\\
1572864000	1.44\\
3145728000	2.25\\
6291456000	2.21\\
};
\addlegendentry{3M dof x GPU}

\addplot [color=mycolor3, line width=3.0pt, mark=triangle, mark options={solid, rotate=180, mycolor3}]
  table[row sep=crcr]{%
256000	0.2\\
500000	0.21\\
1099000	0.28\\
2048000	0.38\\
4000000	0.44\\
7813000	0.64\\
16384000	0.61\\
32000000	0.66\\
62500000	0.81\\
131072000	0.8\\
256000000	0.84\\
500000000	1.31\\
1048576000	0.95\\
2048000000	1.14\\
4000000000	2.27\\
5695313000	1.45\\
};
\addlegendentry{256k dof x MPI core}

\addplot [color=mycolor4, line width=3.0pt, mark=triangle, mark options={solid, rotate=180, mycolor4}]
  table[row sep=crcr]{%
500000	0.39\\
1099000	0.48\\
2048000	0.61\\
4000000	0.85\\
7813000	0.86\\
16384000	1.15\\
32000000	1.15\\
62500000	1.39\\
131072000	1.3\\
256000000	1.29\\
500000000	1.54\\
1048576000	1.4\\
1755904000	1.4\\
4000000000	1.79\\
8388608000	1.71\\
13500000000	1.88\\
};
\addlegendentry{512k dof x MPI core}

\addplot [color=mycolor5, line width=3.0pt, mark=triangle, mark options={solid, rotate=180, mycolor5}]
  table[row sep=crcr]{%
977000	0.79\\
1953000	1.13\\
4000000	1.21\\
7813000	1.66\\
16384000	1.83\\
32000000	2.26\\
62500000	2.24\\
131072000	2.59\\
256000000	2.55\\
500000000	2.89\\
1048576000	2.75\\
2048000000	2.77\\
4000000000	4.06\\
8388608000	3.7\\
16384000000	3.71\\
26367188000	5.55\\
};
\addlegendentry{1M dof x MPI core}

\node[above left, align=right, font=\itshape]
at (axis cs:977000,0.79) {    1};
\node[above left, align=right, font=\itshape]
at (axis cs:1953000,1.13) {    2};
\node[above left, align=right, font=\itshape]
at (axis cs:4000000,1.21) {    4};
\node[above left, align=right, font=\itshape]
at (axis cs:7813000,1.66) {    8};
\node[above left, align=right, font=\itshape]
at (axis cs:16384000,1.83) {   16};
\node[above left, align=right, font=\itshape]
at (axis cs:32000000,2.26) {   32};
\node[above left, align=right, font=\itshape]
at (axis cs:62500000,2.24) {   64};
\node[above left, align=right, font=\itshape]
at (axis cs:131072000,2.59) {  128};
\node[above left, align=right, font=\itshape]
at (axis cs:256000000,2.55) {  256};
\node[above left, align=right, font=\itshape]
at (axis cs:500000000,2.89) {  512};
\node[below left, align=right, font=\itshape]
at (axis cs:1048576000,2.75) { 1024};
\node[above left, align=right, font=\itshape]
at (axis cs:2048000000,2.77) { 2048};
\node[above left, align=right, font=\itshape]
at (axis cs:4000000000,4.06) { 4096};
\node[below left, align=right, font=\itshape]
at (axis cs:8388608000,3.7) { 8192};
\node[above left, align=right, font=\itshape]
at (axis cs:16384000000,3.71) {16384};
\node[above left, align=right, font=\itshape]
at (axis cs:26367188000,5.55) {27000};
\node[below left, align=right]
at (axis cs:3072000,0.16) {   1};
\node[below left, align=right]
at (axis cs:6144000,0.22) {   2};
\node[below left, align=right]
at (axis cs:12288000,0.28) {   4};
\node[below right, align=right]
at (axis cs:24576000,0.49) {   8};
\node[below right, align=right]
at (axis cs:49152000,0.56) {  16};
\node[below right, align=right]
at (axis cs:98304000,0.83) {  32};
\node[below right, align=right]
at (axis cs:196608000,0.95) {  64};
\node[below right, align=right]
at (axis cs:393216000,1.22) { 128};
\node[below right, align=right]
at (axis cs:786432000,1.4) { 256};
\node[below right, align=right]
at (axis cs:1572864000,1.44) { 512};
\node[below right, align=right]
at (axis cs:3145728000,2.25) {1024};
\node[above right, align=right]
at (axis cs:6291456000,2.21) {2048};
\end{axis}
\end{tikzpicture}%
}
	
	\subfloat[\label{fig:KCycle} KA1S1CS1 preconditioner on MPI and VA3S5CS1 on MPI-CUDA ]{
%
%
\definecolor{mycolor1}{rgb}{0.00000,0.44700,0.74100}%
\definecolor{mycolor2}{rgb}{0.85000,0.32500,0.09800}%
\definecolor{mycolor3}{rgb}{0.92900,0.69400,0.12500}%
\definecolor{mycolor4}{rgb}{0.49400,0.18400,0.55600}%
\definecolor{mycolor5}{rgb}{0.46600,0.67400,0.18800}%
\begin{tikzpicture}

\begin{axis}[%
width=0.87\textwidth,
height=2.8in,
at={(2.509in,2.768in)},
scale only axis,
xmode=log,
xmin=100000,
xmax=100000000000,
xminorticks=true,
xlabel style={font=\color{white!15!black}},
xlabel={dofs},
ymode=log,
ymin=0.1,
ymax=21,
yminorticks=true,
title={Execution Time for Solve (sec.)},
axis background/.style={fill=white},
ymajorgrids,
yminorgrids,
legend columns=3,
legend style={at={(0.5,-0.15)}, anchor=north, legend cell align=left, align=left, draw=none,font=\footnotesize}
]
\addplot [color=mycolor1, line width=3.0pt, mark=o, mark options={solid, mycolor1}]
  table[row sep=crcr]{%
6144000	0.28\\
12288000	0.38\\
24576000	0.52\\
49152000	0.65\\
98304000	0.79\\
196608000	0.92\\
393216000	1.29\\
786432000	1.37\\
1572864000	1.48\\
3145728000	2.73\\
6291456000	2.62\\
12582912000	2.85\\
};
\addlegendentry{6M dof x GPU}

\addplot [color=mycolor2, line width=3.0pt, mark=o, mark options={solid, mycolor2}]
  table[row sep=crcr]{%
3072000	0.16\\
6144000	0.22\\
12288000	0.28\\
24576000	0.49\\
49152000	0.56\\
98304000	0.83\\
196608000	0.95\\
393216000	1.22\\
786432000	1.4\\
1572864000	1.44\\
3145728000	2.25\\
6291456000	2.21\\
};
\addlegendentry{3M dof x GPU}

\addplot [color=mycolor3, line width=3.0pt, mark=triangle, mark options={solid, rotate=180, mycolor3}]
  table[row sep=crcr]{%
256000	0.29\\
500000	0.31\\
1099000	0.42\\
2048000	0.58\\
4000000	0.62\\
7813000	0.97\\
16384000	0.93\\
32000000	1.04\\
62500000	2.03\\
131072000	2.49\\
256000000	3.01\\
500000000	3.16\\
1048576000	3.45\\
2048000000	4.17\\
4000000000	5.41\\
5695313000	7.58\\
};
\addlegendentry{256k dof x MPI core}

\addplot [color=mycolor4, line width=3.0pt, mark=triangle, mark options={solid, rotate=180, mycolor4}]
  table[row sep=crcr]{%
500000	0.55\\
1099000	0.7\\
2048000	0.85\\
4000000	1.23\\
7813000	1.2\\
16384000	1.74\\
32000000	1.81\\
62500000	2.18\\
131072000	2.09\\
256000000	2.27\\
500000000	3.24\\
1048576000	2.95\\
1755904000	2.92\\
4000000000	3.71\\
8388608000	4.03\\
13500000000	4.21\\
};
\addlegendentry{512k dof x MPI core}

\addplot [color=mycolor5, line width=3.0pt, mark=triangle, mark options={solid, rotate=180, mycolor5}]
  table[row sep=crcr]{%
977000	1.13\\
1953000	1.4\\
4000000	1.83\\
7813000	2.31\\
16384000	2.77\\
32000000	3.49\\
62500000	3.65\\
131072000	4.07\\
256000000	4.05\\
500000000	4.71\\
1048576000	5.4\\
2048000000	5.23\\
4000000000	9.41\\
8388608000	11.47\\
16384000000	12.44\\
26367188000	15.76\\
};
\addlegendentry{1M dof x MPI core}

\node[above left, align=right, font=\itshape]
at (axis cs:977000,1.13) {    1};
\node[above left, align=right, font=\itshape]
at (axis cs:1953000,1.4) {    2};
\node[above left, align=right, font=\itshape]
at (axis cs:4000000,1.83) {    4};
\node[above left, align=right, font=\itshape]
at (axis cs:7813000,2.31) {    8};
\node[above left, align=right, font=\itshape]
at (axis cs:16384000,2.77) {   16};
\node[above left, align=right, font=\itshape]
at (axis cs:32000000,3.49) {   32};
\node[above left, align=right, font=\itshape]
at (axis cs:62500000,3.65) {   64};
\node[above left, align=right, font=\itshape]
at (axis cs:131072000,4.07) {  128};
\node[above left, align=right, font=\itshape]
at (axis cs:256000000,4.05) {  256};
\node[above left, align=right, font=\itshape]
at (axis cs:500000000,4.71) {  512};
\node[above left, align=right, font=\itshape]
at (axis cs:1048576000,5.4) { 1024};
\node[below left, align=right, font=\itshape]
at (axis cs:2048000000,5.23) { 2048};
\node[above left, align=right, font=\itshape]
at (axis cs:4000000000,9.41) { 4096};
\node[above left, align=right, font=\itshape]
at (axis cs:8388608000,11.47) { 8192};
\node[above, align=right, font=\itshape]
at (axis cs:16384000000,12.44) {16384};
\node[below right, align=right, font=\itshape]
at (axis cs:26367188000,15.76) {27000};
\node[below left, align=right]
at (axis cs:3072000,0.16) {   1};
\node[below left, align=right]
at (axis cs:6144000,0.22) {   2};
\node[below left, align=right]
at (axis cs:12288000,0.28) {   4};
\node[below left, align=right]
at (axis cs:24576000,0.49) {   8};
\node[below left, align=right]
at (axis cs:49152000,0.56) {  16};
\node[below left, align=right]
at (axis cs:98304000,0.83) {  32};
\node[below left, align=right]
at (axis cs:196608000,0.95) {  64};
\node[below left, align=right]
at (axis cs:393216000,1.22) { 128};
\node[above left, align=right]
at (axis cs:786432000,1.4) { 256};
\node[below left, align=right]
at (axis cs:1572864000,1.44) { 512};
\node[below left, align=right]
at (axis cs:3145728000,2.25) {1024};
\node[below right, align=right]
at (axis cs:6291456000,2.21) {2048};
\end{axis}

\end{tikzpicture}%
}
	
	\caption{Weak scaling results. Comparison of the execution time for solve between the pure MPI implementation of the VA3S1CS1 and KA1S1CS1 preconditioners, and the hybrid MPI-CUDA implementation of the VA3S5CS1 preconditioner.}
	\label{fig:global_behavior}
\end{figure}
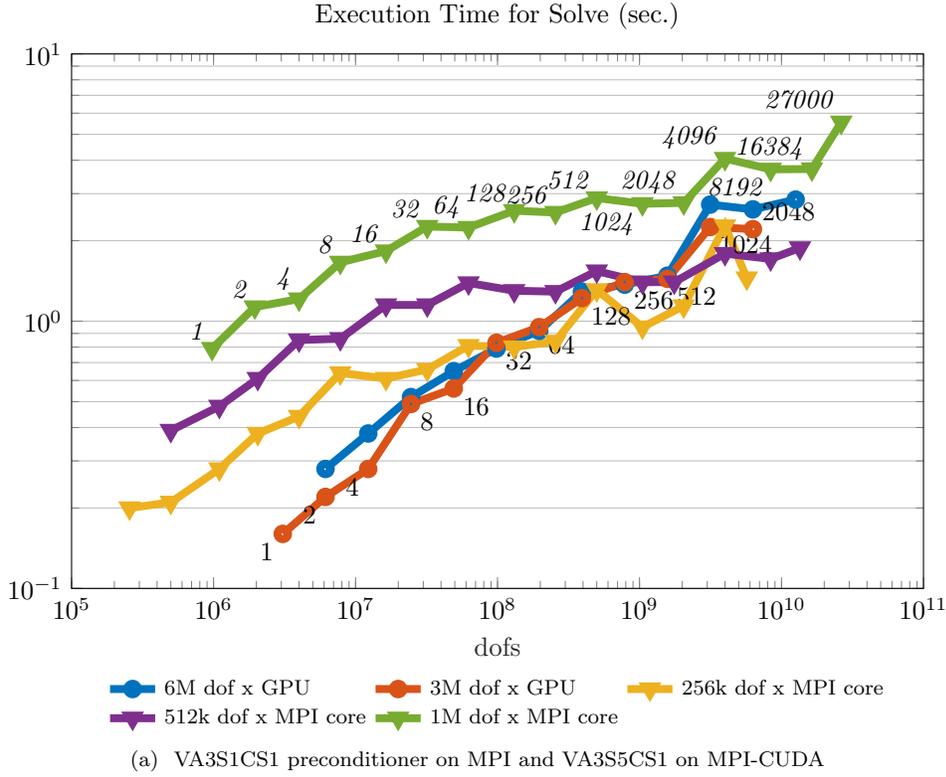
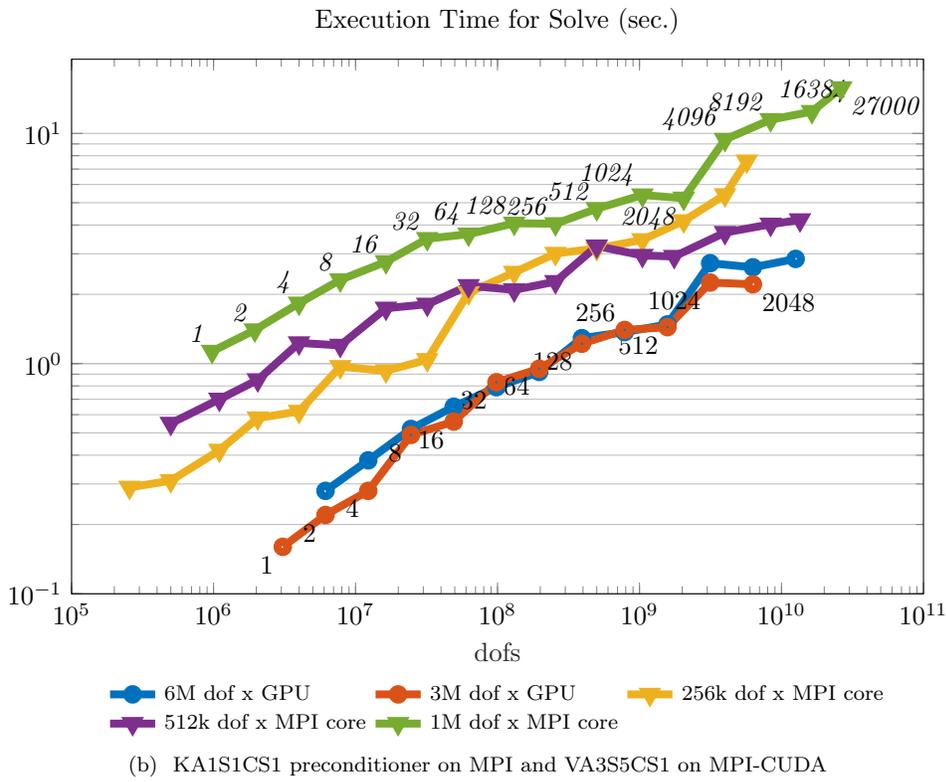
Finally, we consider Fig.~\ref{fig:VSCycle} and Fig.~\ref{fig:KCycle} in which we summarize the result that are discussed in detail in section~\ref{res:pizdaintcpu-gpu}. We can observe again how the smoothed version of the preconditioner using a V-cycle consistently achieves faster solution times with respect to the corresponding unsmoothed version. Moreover, when we compare these results with the GPU version of the VA3S5CS1 preconditioner we observe that the hybrid approach permits savings in the solve time and, possibly even more importantly, large savings in the energy consumption, since to reach the same number of dofs and very similar execution time per linear iteration we need to use less CPU nodes.

\section{Software} The development versions of the \texttt{PSBLAS}, \texttt{PSBLAS-GPU} and \texttt{AMG4PSBLAS}, can be obtained from the website \url{https://psctoolkit.github.io/}. The new extension, containing the aggregation strategy based on graph matching will be made available by the end of the EoCoE-II project at the same address, together with the releases of the new versions of the libraries composing the whole software toolkit.

\end{document}